\documentclass[final]{ws-m3as}

\usepackage{color}  
\definecolor{darkgreen}{rgb}{0, 0.7, 0}
\usepackage{soul}

\usepackage{setspace}
\usepackage{amsmath,amsfonts,amssymb,amsbsy,amscd}
\usepackage{mathrsfs,array}
\usepackage{tikz}
\usepackage{upgreek}
\usetikzlibrary{intersections, arrows, automata, backgrounds, calendar, chains, matrix, mindmap, patterns, petri, shadows, shapes.geometric, shapes.misc, spy, trees, shapes}
\usepackage{subcaption}
\usepackage[numbers]{natbib}
\usepackage{tikz}
\usetikzlibrary{arrows.meta}
\usetikzlibrary{3d}
\usetikzlibrary{shadings}
\usetikzlibrary{math}
\usepackage{upgreek}
\usetikzlibrary{intersections, arrows, automata, backgrounds, calendar, chains, matrix, mindmap, patterns, petri, shadows, shapes.geometric, shapes.misc, spy, trees, shapes}

\usepackage{hyperref}
\hypersetup{
    colorlinks=true,
    linkcolor=blue,
    urlcolor=red,
    linktoc=all
}

\newcommand{\jump}[1]{\mbox{$[\![ #1 ]\!]$}}
\newcommand{\avg}[1]{\mbox{$\{\!\>\!\>\!\!\!\{ #1 \}\!\>\!\>\!\!\!\}$}}

\newcommand{\nn}{\nonumber}

\newcommand{\bx}{\mathbf{x}}

\newcommand{\HH}{\mathcal{H}}

\def\be{\begin{equation}}
\def\ee{\end{equation}}
\def\ba{\begin{array}}
\def\ea{\end{array}}
\def\bea{\begin{eqnarray}}
\def\eea{\end{eqnarray}}
\def\beas{\begin{eqnarray*}}
\def\eeas{\end{eqnarray*}}
\newcommand{\bseq}{\begin{subequations}}
\newcommand{\eseq}{\end{subequations}}

\definecolor{blue1}{rgb}{0,0.4470,0.7410}
\definecolor{red1}{rgb}{0.8500,0.3250,0.0980}
\definecolor{green1}{rgb}{0.4660,0.6740,0.1880}
\definecolor{carnelian}{rgb}{0.7, 0.11, 0.11}
\definecolor{cadmiumgreen}{rgb}{0.0, 0.42, 0.24}

\newcommand*{\xMin}{0}%
\newcommand*{\xMax}{8}%
\newcommand*{\yMin}{1}%
\newcommand*{\yMax}{9}%

\newcommand*{\xMinone}{0}%
\newcommand*{\xMaxone}{8}%
\newcommand*{\yMinone}{0}%
\newcommand*{\yMaxone}{8}%

\newcommand*\elevation{110}
\newcommand*\anglerot{90}
\pgfmathsetmacro\xc{cos(\anglerot)}   
\pgfmathsetmacro\xs{sin(\anglerot)} 
\pgfmathsetmacro\yc{cos(\elevation)}   
\pgfmathsetmacro\ys{sin(\elevation)}

\tikzmath{\x2 = 2.6; \y2 = 1.6; \x3 = 0.4; \y3 = 2.4; \x1 = \x2 + \x3; \y1 = \y2 + \y3; \u2 = 2.6; \v2 = 0.1; \u3 = -0.6; \v3 = 2.4; \u1 = \u2 + \u3; \v1 = \v2 + \v3; \w2 = 3.0; \p2 = 0.0001; \w3 = 0.0001; \p3 = 2.4; \w1 = \w2 + \w3; \p1 = \p2 + \p3; \s1 = 6.0; \s2 = 12.0; \xx2 = \x2; \yy2 = \y2; \xx3 = \x3; \yy3 = \y3; \xx1 = \xx2 + \xx3; \yy1 = \yy2 + \yy3; \r = 0.3; \t = 0.4; \xx5 = \u2; \yy5 = \v2; \xx6 = \u3; \yy6 = \v3; \xx4 = \xx5 + \xx6; \yy4 = \yy5 + \yy6; \xx8 = \w2; \yy8 = \p2; \xx9 = \w3; \yy9 = \p3; \xx7 = \xx8 + \xx9; \yy7 = \yy8 + \yy9;} 

\usepackage[nameinlink,capitalise,noabbrev]{cleveref}

\makeatother

\begin{document}

\markboth{M.F.P. ten Eikelder et al.}{Constraints for eliminating the Gibbs phenomenon in finite element approximation spaces}

%
\catchline{}{}{}{}{}
%

\title{Constraints for eliminating the Gibbs phenomenon in finite element approximation spaces}
\author{Marco F.P. ten Eikelder\footnote{Corresponding author.}}

\address{Institute for Mechanics, Computational Mechanics Group\\
Technical University of Darmstadt\\
Franziska-Braun-Straße 7, 64287 Darmstadt, Germany\\
marco.eikelder@tu-darmstadt.de}

\author{Stein K.F. Stoter}
\address{Department of Mechanical Engineering\\
Eindhoven University of Technology\\
P.O. Box 513, 5600 MB, Eindhoven, The Netherlands\\
s.k.f.stoter@tue.nl}

\author{Yuri Bazilevs}
\address{School of Engineering\\
Brown University\\
184 Hope Street, Providence, RI 02912, United States\\
yuri\_bazilevs@brown.edu}

\author{Dominik Schillinger}
\address{Institute for Mechanics, Computational Mechanics Group\\
Technical University of Darmstadt\\
Franziska-Braun-Straße 7, 64287 Darmstadt, Germany\\
dominik.schillinger@tu-darmstadt.de}

\maketitle


\begin{abstract}
One of the major challenges in finite element methods is the mitigation of spurious oscillations near sharp layers and discontinuities known as the Gibbs phenomenon. 
In this article, we propose a set of functionals to identify spurious oscillations in best approximation problems in finite element spaces. Subsequently, we adopt these functionals in the formulation of constraints in an effort to eliminate the Gibbs phenomenon. By enforcing these constraints in best approximation problems, we can entirely eliminate over- and undershoot in one dimensional continuous approximations, and significantly suppress them in one- and higher-dimensional discontinuous approximations.
\end{abstract}



\section{Introduction}

\subsection{Historical overview}
The discovery of the Gibbs phenomenon may be traced back to Henry Wilbraham (1848), and the phenomenon was rediscovered by J. Willard Gibbs (1898-1899), in their studies on Fourier series 
\cite{hewitt1979gibbs,gottlieb1997gibbs}. 
It is traditionally described as the inability to recover point values of a discontinuous function by a truncated Fourier expansion. Near the discontinuity, the error does not vanish as the number of terms in the expansion is increased, and the magnitude of the over- and undershoots tends to a fixed limit. The limiting value is known as the Gibbs constant. It is less well known that the Gibbs phenomenon also occurs in truncated expansions of other sets of orthogonal functions \cite{jacob1937sur,kaber2006the}. In fact, the associated Gibbs constants are often identical, as is the case for expansions with Legendre, Hermite, or Laguerre polynomials \cite{fay2006gibbs,davis2022gibbs}. 

Fundamentally, the Gibbs phenomenon has, however, little to do with Fourier series or expansions in orthogonal polynomials. The effect arises from the best approximation in a square integral metric, of which these expansions are examples \cite{foster1991gibbs}. As such, it also occurs in best approximation problems in the $L^2$-metric by piecewise linear polynomials \cite{foster1991gibbs} or splines \cite{richards1991gibbs}. The role of the metric herein is crucial: spurious oscillations that appear in the $L^2$-metric are significantly more  severe than those that occur in the $L^q$-metric when $q$ tends to $1$, for which they are in some cases even completely absent \cite{guermond2004finite}. A detailed study on the possible elimination of the Gibbs phenomenon in $L^q$-best approximation by piecewise linear finite element shape functions is presented in \cite{houston2022gibbs}. In the last few decades, the $L^1(\Omega)$ functional settings has hence been explored as the point of departure for approximating solutions to partial differential equations (PDEs) \cite{lavery1988nonoscillatory,lavery1991solution,lowrie1994numerical,jiang1993non,guermond2004finite}. The main challenge behind these approaches is that they require the minimization of a nondifferentiable functional, which leads to a poorly behaved nonlinear problem. As a consequence, there is a lack of practical algorithms for solving even standard problems in computational mathematics. Additionally, even though approximations in subspaces of $L^1(\Omega)$ reduce spurious oscillations, on some meshes these do not vanish in general \cite{houston2020eliminating,houston2022gibbs}

The more conventional (Bubnov-)Galerkin method produces solutions that are optimal in an inner product induced norm (associated with subspaces of $L^2(\Omega)$). As such, approximations of interior- and boundary layers indeed tend to suffer from spurious oscillations. This issue is well known in the finite element community, and many attempts at resolutions have been proposed. Arguably the most successful 
is the class of residual-based stabilized methods \cite{BroHug82,HugFra86,HugFra89}, which are primarily adopted for fluid mechanics related applications. Residual-based stabilization significantly improves the solution quality in regions free of abrupt changes, but the Gibbs phenomenon still occurs in regions with sharp layers. As a remedy, the finite element formulation is often augmented with a nonlinear stabilization mechanism that locally introduces artificial diffusion \cite{hughes1986new2,BaCaTeHu07,ten2020theoretical}. These methods are referred to as shock- or discontinuity capturing methods.  

In the case of nonlinear (hyperbolic) evolution equations, the above stabilization methods do still not suffice. In order to enhance the quality of numerical approximations for these types of problems, algorithms have been designed to inherit certain stability properties of the underlying PDE. The prevalent example is the entropy stability property possessed by entropy solutions. Weak solutions of nonlinear evolution equations are not unique and the entropy stability property singles out the entropy solution as the physically relevant solution \cite{kruvzkov1970first}.
The entropy stability concept, which reduces for many physical systems to an energy-dissipation property, has frequently been used in the construction of stable finite element methods \cite{liu2013functional,liu2015liquid,EiAk17i,EiAk17ii,ten2021novel,guermond2011entropy}.
Even though the solution quality enhances significantly, numerical solutions that inherit entropy stability do not preclude spurious oscillations. For particular variable sets, the Galerkin method may even exactly satisfy the entropy stability condition but still exhibit spurious oscillations \cite{hughes1986new}.
Evidently, the entropy stability concept is not inextricably linked to the Gibbs phenomenon.
It does, however, seem to be a good indicator for the identification of shock waves, and thus as an indicator of where the Gibbs phenomenon might manifest.

A stability concept that is more directly targeted at removing the Gibbs phenomenon, is the total variation diminishing (TVD) property introduced by Harten \cite{harten1983high,harten1984class}.
Solutions with the TVD property preclude the growth of the total variation of the solution.
The design of numerical schemes with the TVD property still is an active area of research.
The incentive for the design of TVD schemes is the desire to produce numerical approximations that satisfy the maximum principle, as well as certain monotonicity properties.
Despite its success, particularly in the finite difference and finite volume communities, the applicability of TVD schemes limited.
Namely, it is solely suitable for time-dependent and scalar conservation laws and does not provide any information on local solution quality.
Moreover, its introduction in the discrete setting relies on a Cartesian grid and lacks frame-invariance.

The above two observations, namely (i) the occurrence of the Gibbs phenomenon in entropy stable discrete solutions, and (ii) the limitations of TVD schemes, have incentivized the design of a novel stability concept called \textit{variation entropy theory} \cite{ten2019variation}. This theory provides a local continuous generalization of the TVD stability condition for general conservation laws in an entropy framework. Similar to classical entropy solutions, variation entropy solutions satisfy an underlying stability condition. This stability condition serves in the discrete setting as an indicator of the Gibbs phenomenon. It has successfully been employed in the variation multiscale (VMS) paradigm \cite{Hug95,Hug98,BaCaCoHu07} to design a framework for discontinuity capturing methods \cite{ten2020theoretical}.

\subsection{Objective}
Despite the significant attention it has gained, a precise mathematical definition of the Gibbs phenomenon does not exist.
Any attempt at eliminating the Gibbs phenomenon thus first requires an identification strategy. The identifier that we develop is rooted in variation entropy theory. 
We then propose is to eliminate the Gibbs phenomenon via enforcement of constraints. 

This brings us to the main objective of this article: \textit{to identify a set of practical constraints that aim to eliminate the Gibbs phenomenon in the approximation of sharp layers and discontinuities in finite element spaces}.
To facilitate the analysis, we discuss our results in the isogeometric analysis framework, which we think of as a generalization of $\mathcal{C}^0$- and $\mathcal{C}^{-1}$-finite element spaces to higher order continuity.

Some remark are in order. First, it may seem feasible to construct constraints that remove oscillations by explicitly choosing the coefficients of the basis functions such that the numerical approximation does not exceed bounds of the analytical profile. This is however not a strategy that is realizable in practical computations. The challenge is thus to establish a set of constraints that can be adopted in practice. Second, the idea of \textit{a priori} enforcing constraints in numerical methods is not new. A notable contribution in this regard is the work of Evans et al. \cite{evans2009enforcement}, in which a framework for the enforcement of constraints in the VMS framework is presented.

\subsection{Main results}
The main result of this paper is a set of integral constraints that aim to identify and eliminate the Gibbs phenomenon. 
The occurrence of the Gibbs phenomenon in a certain approximation can not solely be inferred from the approximation itself. Rather, it stands in relation to the function being approximated, and depends on the (sub)domain of interest. We propose an indicator of the form:
\begin{align}\label{eq: intro: Gibbs constraint}
    \mathscr{G}_{\phi,\omega}(\phi^*) \leq 0,
\end{align}
where the function $\phi^* \in H^1(\tilde{\Omega})$ is an approximation of the function $\phi\in H^1(\Omega)$ on $\omega \subset \Omega$. Here, $H^1(\tilde{\Omega})$ is a broken Sobolev space and  $\tilde{\Omega}$ is a collection of disjoint subdomains (precise definitions are provided in \cref{sec: Gibbs constraints 1D}). We call this constraint the \textit{Gibbs constraint}. It follows from the \textit{Gibbs functional}, which we define as:
\begin{align}\label{eq intro: large g}
    \mathscr{G}_{\phi,\omega}(\phi^*) := \displaystyle\int_{\omega} g_{\phi}(\phi^*)~{\rm d}x,
\end{align}
with $g_\phi$ as:
\begin{align}\label{eq intro: small g}
    g_{\phi}(\phi^*):= \left\{ \begin{matrix}
       \|\nabla \phi^*\|_2^{-1}\nabla \phi^*\cdot\nabla (\phi^*-\phi) & \text{ for } \nabla\phi^* \neq 0,\\
      0 & \text{ for } \nabla\phi^* = 0.\end{matrix}\right.
\end{align}
We search functions $\phi^*$ as an approximation of $\phi$ for which the Gibbs constraint is satisfied on predetermined sets of subdomains $\omega$. We study the application of the Gibbs constraints in the context of finite element best approximation problems. In particular, we consider the constrained best approximation problem:
\begin{align}
    \phi^h = \underset{\theta^h \in \mathcal{K}_{p,\alpha}}{\rm arginf}~ \|\phi-\theta^h\|_{\mathcal{H}},
\end{align}
where $\|\cdot \|_{\mathcal{H}}$ is a norm induced by a certain inner product and the feasible set is given by:
\begin{align}
  \mathcal{K}_{p,\alpha} := \left\{ \phi^h \in \mathcal{V}^h_{D,p,\alpha} :~\right.&\left. \mathscr{G}_{\phi,\omega_j}(\phi^*) \leq 0,~ j = 1,...,J, ~\omega_j \in \mathcal{T}_\omega \right\}.
\end{align}
Precise definitions of $\mathcal{V}^h_{D,p,\alpha}$ and $\mathcal{T}_\omega$ are provided in later in the paper. We demonstrate for sharp layers that finite element approximations of arbitrary degree and continuity are free of over- and undershoots when they satisfy the Gibbs constraint (for one dimensional continuous approximations), or these oscillations are significantly suppressed (for discontinuous approximation). 

The choice of subdomains $\omega$ depends on the regularity $\alpha$ of the finite element approximation space and the dimension of the domain. In one dimension, the constraints may be applied element-wise ($\omega_j = K_j$) when the finite element space is either discontinuous or $\mathcal{C}^0$-continuous. For higher regularity finite element spaces ($\alpha \geq 1$), the subdomains $\omega$ need to be collections of neighboring elements, and the number of collected elements increases with the regularity. In higher dimensions, the Gibbs constraints are too restrictive for continuous finite element spaces, 
 limiting its applicability 
 to discontinuous finite element spaces.

\subsection{Outline}
The remainder of the paper is structured as follows.
First, in \cref{sec: prelim} we provide preliminaries concerning function spaces and projectors. Then, in \cref{sec: Gibbs review} we present an overview of the Gibbs phenomenon for best approximations in finite element spaces of arbitrary degree and continuity.
In \cref{sec: Gibbs constraints 1D} we present the identification of the constraints in one spatial dimension.
Next, we extend our construction to higher dimensions in \cref{sec: Gibbs higher D}.
Finally, we provide a summary and outlook in \cref{sec: conclusions}.

\section{Preliminaries}\label{sec: prelim}
\subsection{Function spaces}\label{subsec: func spaces}
We adopt the standard functional analysis setting. We denote by $\Omega \subset \mathbb{R}^d$ the bounded, open and connected domain with spatial dimension~$d$, and with boundary $\partial \Omega$. $L^2(\Omega)$ is the Lebesgue space of $2$-integrable functions on $\Omega$. Furthermore, $H^{1}(\Omega) \subset L^2(\Omega)$ is the Sobolev space of $L^2(\Omega)$-functions with their gradient also in $[L^2(\Omega)]^d$. The subspace $H_0^{1}(\Omega) \subset H^{1}(\Omega)$ consists of functions with zero trace on $\partial \Omega$. The associated norms are denoted as $\|\cdot\|_{L^2(\Omega)}=\|\cdot\|_{\Omega}$ and $\|\cdot\|_{H_0^{1}(\Omega)}$. Furthermore, we use standard notation for the $L^2$-inner product on $\Omega$, $(\cdot,\cdot)_{L^2(\Omega)} = (\cdot,\cdot)_\Omega$ and write $\langle \cdot, \cdot \rangle_{D}$ for the duality pairing $H^{-1/2}(D) \times H^{1/2}(D) \rightarrow \mathbb{R}$ on some boundary domain $D$. 

In this article, we consider finite element spaces of arbitrary degree and continuity. As such, we make use of the \textit{isogeometric analysis} framework \cite{HuCoBa04,CoHuBa09,BBCHS06}. We introduce knotvectors, univariate and multivariate B-splines, geometrical mappings and the physical mesh. The ordered knotvector $\Xi$ is defined for degree $p$ and dimensionality $n$ as:
\begin{align}
    \Xi := \left\{-1=\xi_{1},\xi_{2},...,\xi_{n+p+1}=1\right\},
\end{align}
where $\xi_i \in \mathbb{R}$ represents the $i$-th knot with $i = 1,...,n+p+1$. We adopt the convention that $p=0,1,2,\dots$ refers to piecewise constants, linears, quadratics, etc. In this work we restrict ourselves to \textit{open} knotvectors, meaning that the first and last knot appear $p+1$ times. The univariate B-spline basis functions are defined recursively for $p=0,1,2,\dots$. Starting with piecewise constant functions, we have:
\begin{align}
    N_{i,0}(\xi) =  \left\{ 
                      \begin{matrix}
                        1 & \text{if } \xi_i \leq \xi < \xi_{i+1}\\
                        0 & \text{otherwise},
                      \end{matrix}
                    \right.
\end{align}
whereas for $p=1,2,\dots$ the B-spline basis functions are given by:
\begin{align}
    N_{i,p}(\xi) = \dfrac{\xi-\xi_{i}}{\xi_{i+p}-\xi_{i}} N_{i,p-1}(\xi) + \dfrac{\xi_{i+p+1}-\xi}{\xi_{i+p+1}-\xi_{i+1}} N_{i+1,p-1}(\xi).
\end{align}
This definition is augmented with the convention that if a denominator (i.e. $\xi_{i+p}-\xi_{i}$ or $\xi_{i+p+1}-\xi_{i+1}$) is zero, that fraction is taken as zero. B-spline basis functions coincide with standard finite element Lagrange basis functions for $p=0$ and $1$, and differ for $p \geq 2$. The set of B-spline basis functions of degree $p$ consists of non-negative piecewise $p$th-order polynomial functions with local support, that form a partition of unity. Linear combinations of B-spline basis functions are referred to as B-splines. We introduce the vector $\boldsymbol{\zeta} = \left\{\zeta_1,...,\zeta_m\right\}$ consisting of all knots without repetitions. 
The open knot vector implies that the basis functions are interpolatory at the ends of the interval. Inside a knot interval B-spline basis functions are smooth, whereas the repetition of a knot reduces the continuity of the B-spline basis function at that knot. More precisely, a B-spline basis function of degree $p$ at a knot $\xi_i$ with multiplicity $k_i$ has $\alpha_i:=p-k_i$ continuous derivatives at $\xi_i$ (note that $\alpha_1 = \alpha_m = -1$). We denote the space of B-splines of polynomial degree $p$ and regularity $\boldsymbol{\alpha}=\left\{\alpha_1,\dots,\alpha_m\right\}$ as:
\begin{align}
    S_{\boldsymbol{\alpha}}^p:= {\rm span}\left\{N_{i,p}\right\}_{i=1}^n.
\end{align}
B-spline basis functions of degree $p$ with uniform internal multiplicity $p$ are interpolatory and span the same space as standard $\mathcal{C}^0$-Lagrange basis functions. Similarly, B-spline basis functions of degree $p$ with uniform internal multiplicity $p+1$ are discontinuous and span the same space as discontinuous Lagrange basis functions. 

The construction of multivariate B-splines follows from taking a tensor-product of the univariate B-splines. We introduce the open knot vectors:
\begin{align}
    \Xi_l := \left\{\xi_{1,l},\xi_{2,l},...,\xi_{n_l+p_l+1,l}\right\},
\end{align}
for polynomial degrees $p_l$ and dimensionality integers $n_l$ for $l = 1,\dots,d$. 
We define for each knot vector $\Xi_l$ univariate B-spline basis functions $N_{i_l,p_l,l}$ of polynomial degree $p_l$ for $i_l = 1,...,n_l$.
Again we introduce the vector of knots with repetition $\boldsymbol{\zeta}_l = \left\{\zeta_{1,l},...,\zeta_{m_l,l}\right\}$ and regularity vector $\boldsymbol{\alpha}_l=\left\{\alpha_{1,l},\dots,\alpha_{m_l,l}\right\}$. The Cartesian mesh on the parametric domain $\hat{\Omega} = (-1,1)^d\subset \mathbb{R}^d$ is now given by:
\begin{align}
    \hat{\mathcal{T}} = \left\{ Q = \displaystyle\otimes_{l=1,\dots,d} (\zeta_{i_l,l},\zeta_{i_l+1,l}), 1 \leq i_l \leq m_l - 1 \right\}.
\end{align}
The boundary of an open element $Q \in \hat{\mathcal{T}}$ is denoted $\partial Q$. The multivariate tensor-product B-spline basis functions are defined on the parametric mesh $\hat{\mathcal{T}}$ as 
\begin{align}
    N_{i_1,\dots,i_d,p_1,\dots,p_d} = N_{i_1,p_1,1} \otimes \dots \otimes N_{i_d,p_d,d}.
\end{align}
The associated tensor-product B-spline function space on $\hat{\mathcal{T}}$ is given by:
\begin{align}
    S_{\boldsymbol{\alpha}_1,\dots,\boldsymbol{\alpha}_d}^{p_1,\dots,p_d}:= {\rm span}\left\{N_{i_1,\dots,i_d,p_1,\dots,p_d}\right\}_{i_1=1,\dots,i_d=1}^{n_1,\dots,n_d}.
\end{align}
Throughout this paper we restrict ourselves to a uniform regularity vector $\boldsymbol{\alpha}=\boldsymbol{\alpha}_l$ in the interior, i.e. $\alpha_{2,l} = \dots = \alpha_{m-1,l} = \alpha$, and use equal polynomial degrees $p_1= \dots= p_d =p$.
We assume that the physical domain can be exactly described by the continuously differentiable geometrical map (with continuously differentiable inverse) $\mathbf{F}:\boldsymbol{\xi} \in \hat{\Omega} \rightarrow \bx \in \Omega$. The physical mesh on $\Omega$ follows by applying the geometrical map $\mathbf{F}$ on elements of the parametric mesh:
\begin{align}
    \mathcal{T} = \left\{ K: K = \mathbf{F}(Q), Q \in \hat{K} \right\}.
\end{align}
As usual, we demand the element $K \in \mathcal{T}$ to be shape-regular. The boundary of an element $K \in \mathcal{T}$ is denoted as $\partial K$. 
We define the Jacobian of the mapping $\mathbf{F}$ as $\mathbf{J} = \partial \bx/\partial \boldsymbol{\xi}$. Lastly, we introduce the finite element approximation space $\mathcal{V}^h_{p,\alpha}=\left\{\mathbf{F}\left(\mathcal{S}^{p}_{\alpha}\right)\right\}:=\left\{\mathbf{F}\left(\mathcal{S}^{p,\dots, p}_{\boldsymbol{\alpha},\dots, \boldsymbol{\alpha}}\right)\right\}$, and its subspaces $\mathcal{V}^h_{0,p,\alpha} \subset \mathcal{V}^h_{p,\alpha}$ and $\mathcal{V}^h_{D,p,\alpha} \subset \mathcal{V}^h_{p,\alpha}$ consisting of those functions that satisfy homogeneous and inhomogeneous boundary conditions on $\partial \Omega$, respectively. 

\subsection{Projection operators}\label{subsec: proj}
In this subsection we introduce some orthogonal projection operators $\mathscr{P}:\mathcal{V} \rightarrow \mathcal{V}^h_{D,p,\alpha}$. Consider first the (constrained) $\mathcal{H}$-best approximation problem:
\begin{align}\label{eq: best approx H}
    \phi^h = \underset{\theta^h \in \mathcal{V}^h_{p,\alpha}}{\rm arginf}~ \|\phi-\theta^h\|_{\mathcal{H}},
\end{align}
subject to the trace equality:
\begin{align}\label{eq: trace equalities}
    \phi^h|_{\partial \Omega} = \phi|_{\partial \Omega},
\end{align}
where $\|\cdot\|_{\mathcal{H}}$ is a norm induced by the inner product $(\cdot,\cdot)_{\mathcal{H}}$. The constraint \eqref{eq: trace equalities} may be homogenized via the adoption of a lift argument, as is standard in finite element methods, whereby  the approximation space becomes $\mathcal{V}^h_{0;p,\alpha} \subset \mathcal{V}^h_{p,\alpha}$. 
The $\mathcal{H}$-best approximation may be determined by solving the first-order optimality conditions obtained from taking the Gateaux derivative in \eqref{eq: best approx H}:\\

\textit{find $\phi^h \in \mathcal{V}_{D;p,\alpha}^h$ such that for all $w^h \in \mathcal{V}_{0;p,\alpha}^h$:}
\begin{align}
   (\mathscr{P}_{\mathcal{H}}\phi - \phi,w^h)_{\mathcal{H}} = (\phi^h - \phi,w^h)_{\mathcal{H}} = 0.
\end{align}

For approximation spaces $\mathcal{V}^h_{p,\alpha}$ consisting of continuous functions (i.e. $\alpha  \geq 0$), we introduce the $L^2$- and $H_0^1$-orthogonal projectors respectively as:
\begin{subequations}
  \begin{align}
   (\mathscr{P}_{L^2}\phi - \phi,w^h)_{\Omega} = (\phi^h - \phi,w^h)_{\Omega} =&~ 0,\\
    (\mathscr{P}_{H_0^1}\phi - \phi,w^h)_{H_0^1(\Omega)} =(\phi^h - \phi,w^h)_{H_0^1(\Omega)} =&~ 0.
  \end{align}
\end{subequations}

Next, we consider approximation spaces $\mathcal{V}^h_{p,\alpha}$ consisting of discontinuous functions (i.e. $\alpha =-1$). In this context, the standard $H_0^1$-norm is not a suitable best approximation problem. 
In order to introduce a suitable alternative to the $H^1_0$-best approximation problem, we first introduce some additional notation. We define the union of ($n_{\rm el}$) open element domains and the associated interface skeleton as:
\begin{subequations}
\begin{alignat}{2}
  \tilde{\Omega} =&~ \displaystyle\bigcup_{i=1}^{n_{\rm el}} K_i,\label{eq: def tilde Omega fem}\\
  \Gamma =&~ \displaystyle\bigcup_{i=1}^{n_{\rm el}} \partial K_i,
\end{alignat}
\end{subequations}
and introduce $\Gamma^0 = \Gamma \backslash \partial \Omega$ as the interior part of the interface skeleton. Next, we introduce some trace operators that are convenient in the context of discontinuous basis functions. For an interior edge $e$, shared by elements $K^+$ and $K^-$, we define the outward pointing unit normal vectors on $e$ as $\mathbf{n}^+$ and $\mathbf{n}^-$, respectively. Denoting $\phi^+= \phi |_{\partial K^+}$ and $\phi^- = \phi|_{\partial K^-}$ of a scalar quantity $\phi$, we define the average $\avg{\phi}$ and jump $\jump{\phi}$ on $\Gamma^0$ as:
\begin{subequations}
  \begin{alignat}{2}
    \avg{\phi}  =&~ \frac{1}{2}\left(\phi^+ + \phi^-\right), \\
    \jump{\phi} =&~ \phi^+ \mathbf{n}^+ + \phi^-\mathbf{n}^-.
  \end{alignat}
\end{subequations}
For a vector-valued quantity $\boldsymbol{\psi}$ on $e$ we define $\boldsymbol{\psi}^+$ and $\boldsymbol{\psi}^-$ analogously and introduce the average $\avg{\boldsymbol{\psi}}$ on $\Gamma^0$ as:
\begin{align}
    \avg{\boldsymbol{\psi}} = \frac{1}{2}\left(\boldsymbol{\psi}^+ + \boldsymbol{\psi}^-\right).
  \end{align}
We do not require the jump of a vector quantity and leave it undefined. 

We now provide an alternative to the $H^1_0$-best approximation problem that is well-posed. This alternative is closely related the well-known interior penalty method. This method requires the evaluation of boundary flux $\partial_n \phi$ on $\Gamma^0$, which is an unbounded operator and yields a double-valued function for $\phi \in \mathcal{V} = H^1(\Omega)$. To circumvent this issue we first introduce the broken space:
\begin{align}\label{eq: broken space fem}
  H^1(\tilde{\Omega}) = \left\{ \phi \in L^2(\Omega): \phi|_K \in H^1(K) ~\text{ for all }K \in \mathcal{T} \right\},
\end{align}
and consider solution functions $\phi \in \tilde{\mathcal{V}} = H^1(\tilde{\Omega})$. To mitigate the issue of the unbounded boundary flux operator, we introduce a suitable additional function space. Recalling that the boundary flux is double-valued, we introduce the function space of the boundary fluxes as the product space $\mathcal{Q} \times \mathcal{Q}$, where $\mathcal{Q} = H^{-1/2}(\Gamma^0)$. Consider then the following operator:
\begin{align}\label{eq: proj IP}
  \mathscr{P}_{\rm IP}: \tilde{\mathcal{V}} \times \mathcal{Q} \times \mathcal{Q} &\longrightarrow {\rm ran} \mathscr{P}_{\rm IP}, \nonumber \\
  (\phi,\mu^+,\mu^-)  &\longrightarrow \left(\phi^h,\partial_n \phi^{+,h},\partial_n \phi^{-,h}\right),
\end{align}
with
\begin{align}\label{eq: best approx IP}
  \phi^h = \underset{\theta^h \in \mathcal{V}^h_{D,p,\alpha}}{\rm arginf}& \frac{1}{2} \| \phi- \theta^h \|^2_{H_0^1(\tilde{\Omega})} - \left\langle \avg{\mu \mathbf{n}-\nabla\theta^h}, \jump{\phi-\theta^h}\right\rangle_{\Gamma^0} \nn\\
  &~+ \frac{1}{2}\left\langle  \eta \jump{\phi-\theta^h},\jump{\phi-\theta^h}\right\rangle_{\Gamma^0}.
\end{align}
The range of $\mathscr{P}_{\rm IP}$ is given by:
\begin{align}
    {\rm ran} \mathscr{P}_{\rm IP} =  \left\{ (w^h,\partial_n w^{+,h},\partial_n w^{-,h}): w^h \in \mathcal{V}^h_{D,p,\alpha} \right\},
\end{align}
with dimension ${\rm dim}\left({\rm ran} \mathscr{P}_{\rm IP}\right) = {\rm dim} \mathcal{V}^h_{D,p,\alpha}$. Additionally, the mapping $\mathscr{P}_{\rm IP}$ is idempotent, and is a linear and bounded operator on the space $\tilde{\mathcal{V}} \times \mathcal{Q}\times \mathcal{Q}$. As a consequence, $\mathscr{P}_{\rm IP}$ is a projector, and we refer to it as the \textit{interior penalty projector} \cite{stoter2021nitsche,stoter2022discontinuous}. The penalty parameter $\eta$ penalizes mismatches of interface jumps. Well-posedness is ensured when the penalty parameter satisfies a certain lower bound. In this work we base the value of the penalty parameter $\eta$ on the work of Shahbazi \cite{shahbazi2005explicit}.

\begin{remark}[Interpretation boundary flux]\label{rmk: interpretation boundary flux}
  The double-valued quantity $\mu$ acts as a surrogate to $\partial_n \phi$ in \eqref{eq: best approx IP}, and is introduced to make the projector a bounded operator. In practice, we simply use $\nabla \phi$ in place of $\mu \mathbf{n}$ in \eqref{eq: best approx IP}. To clarify the consistency of this replacement, we expand a part of the second integrand in \eqref{eq: best approx IP}:
  \begin{align}
      \avg{\mu \mathbf{n}}\cdot \jump{\phi-\phi^h} =&~ \frac{1}{2}\left( \mu^+ (\phi-\phi^h)^+ - \mu^+(\phi-\phi^h)^- \right.\nn\\
      &~\left.- \mu^-(\phi-\phi^h)^+ + \mu^- (\phi-\phi^h)^- \right).
  \end{align}
  By then replacing $\mu$ by $\partial_n\phi$, we arrive at:
  \begin{align}
      \avg{\partial_n \phi \mathbf{n}}\cdot \jump{\phi-\phi^h} =&~ \frac{1}{2}\left( \mathbf{n}^+\cdot \nabla \phi^+ (\phi-\phi^h)^+ - \mathbf{n}^+\cdot \nabla \phi^+(\phi-\phi^h)^- \right. \nonumber\\
      &~\left.- \mathbf{n}^-\cdot \nabla \phi^-(\phi-\phi^h)^+ + \mathbf{n}^-\cdot \nabla \phi^- (\phi-\phi^h)^- \right) \nonumber\\
      =&~\avg{\nabla \phi}\cdot \jump{\phi-\phi^h},
  \end{align}
  which reveals the connection between $\mu \mathbf{n}$ and $\nabla \phi$. 
\end{remark}

To employ the interior penalty projector for the best approximation of sufficiently smooth functions, we may replace $\mu \mathbf{n}$ by $\nabla \phi$ (see \cref{rmk: interpretation boundary flux}). 
By taking the Gateaux derivative of \eqref{eq: best approx IP}, and using this substitution, we obtain the following first-order optimality condition:\\

\textit{find $\phi^h \in \mathcal{V}_{D;p,0}^h$ such that for all $w^h \in \mathcal{V}_{0;p,0}^h$:}
\begin{align}
    &~( \phi^h - \phi, w^h)_{H_0^1(\tilde{\Omega})} - \langle \jump{\phi^h-\phi} , \avg{\nabla w^h} \rangle_{\Gamma^0} - \langle \avg{\nabla \phi^h-\nabla \phi} , \jump{w^h} \rangle_{\Gamma^0} \nn\\
    &~+ \langle \eta \jump{\phi^h-\phi}, \jump{w^h} \rangle_{\Gamma^0} = 0.
\end{align}
With this substitution it is easy to see that the interior penalty projector is indeed associated with a best approximation problem:
\begin{align}
    \phi^h = \underset{\theta^h \in \mathcal{V}^h_{D,p,-1}}{\rm arginf}~ \|\phi-\theta^h\|_{{\rm IP}(\Omega)},
\end{align}
where the norm is defined as:
\begin{align}
  \|v\|^2_{\rm IP(\Omega)}:= \|v\|_{H_0^1(\tilde{\Omega})}^2 - 2 \left\langle \avg{\nabla v}, \jump{v}\right\rangle_{\Gamma^0} + \left\langle  \eta \jump{v},\jump{v}\right\rangle_{\Gamma^0}.
\end{align}

\section{An exposition of the Gibbs phenomenon for best approximations in finite element spaces}\label{sec: Gibbs review}

In this section we demonstrate the occurrence of the Gibbs phenomenon in best approximation problems that involve finite element approximation spaces of arbitrary continuity and degree. For simplicity, we work with B-splines with equal knot spacing. 
We consider the one-dimensional case in \cref{sec: review Gibbs 1D} and the two-dimensional case in \cref{sec: review Gibbs 2D}.

\subsection{The Gibbs phenomenon in one dimension}\label{sec: review Gibbs 1D}
We consider best approximations of a step function $\phi \in L^2(\Omega)$ defined as:
\begin{align}\label{eq: exact step}
    \phi=\phi_a(x) = \left\{
      \begin{matrix}
        1 & x > a\\
        -1 & x < a
      \end{matrix}
    \right.\,,
\end{align}
where $a$ denotes the location of the jump discontinuity. As some best approximation statements involve weak derivatives, we wish to work with solution functions in $H^1(\Omega)$. Therefore, we introduce the following smooth (differentiable) approximation $\phi \in H^1(\Omega)=:\mathcal{V}$ of the step function:
\begin{align}\label{eq: approx step}
    \phi=\phi_a^{\epsilon}(x) = \tanh\left(\dfrac{x-a}{\epsilon}\right),
\end{align}
where $\epsilon \ll 1$ is a smoothing parameter. 

\begin{figure}[!ht]
\begin{subfigure}{0.49\textwidth}
\centering
\includegraphics[width=0.95\textwidth]{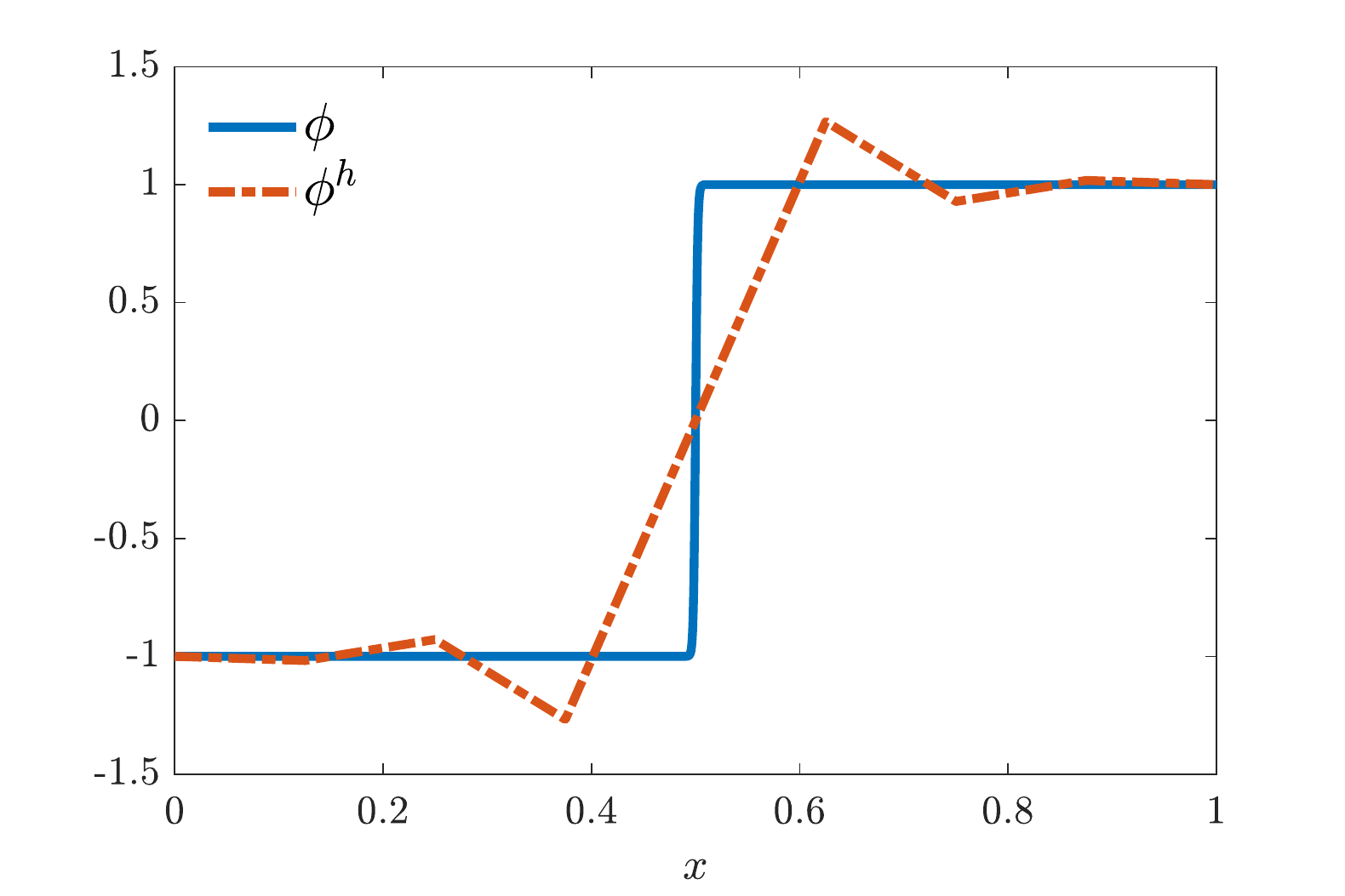}
\caption{$n_{\rm el} =8$.}
\end{subfigure}
\begin{subfigure}{0.49\textwidth}
\centering
\includegraphics[width=0.95\textwidth]{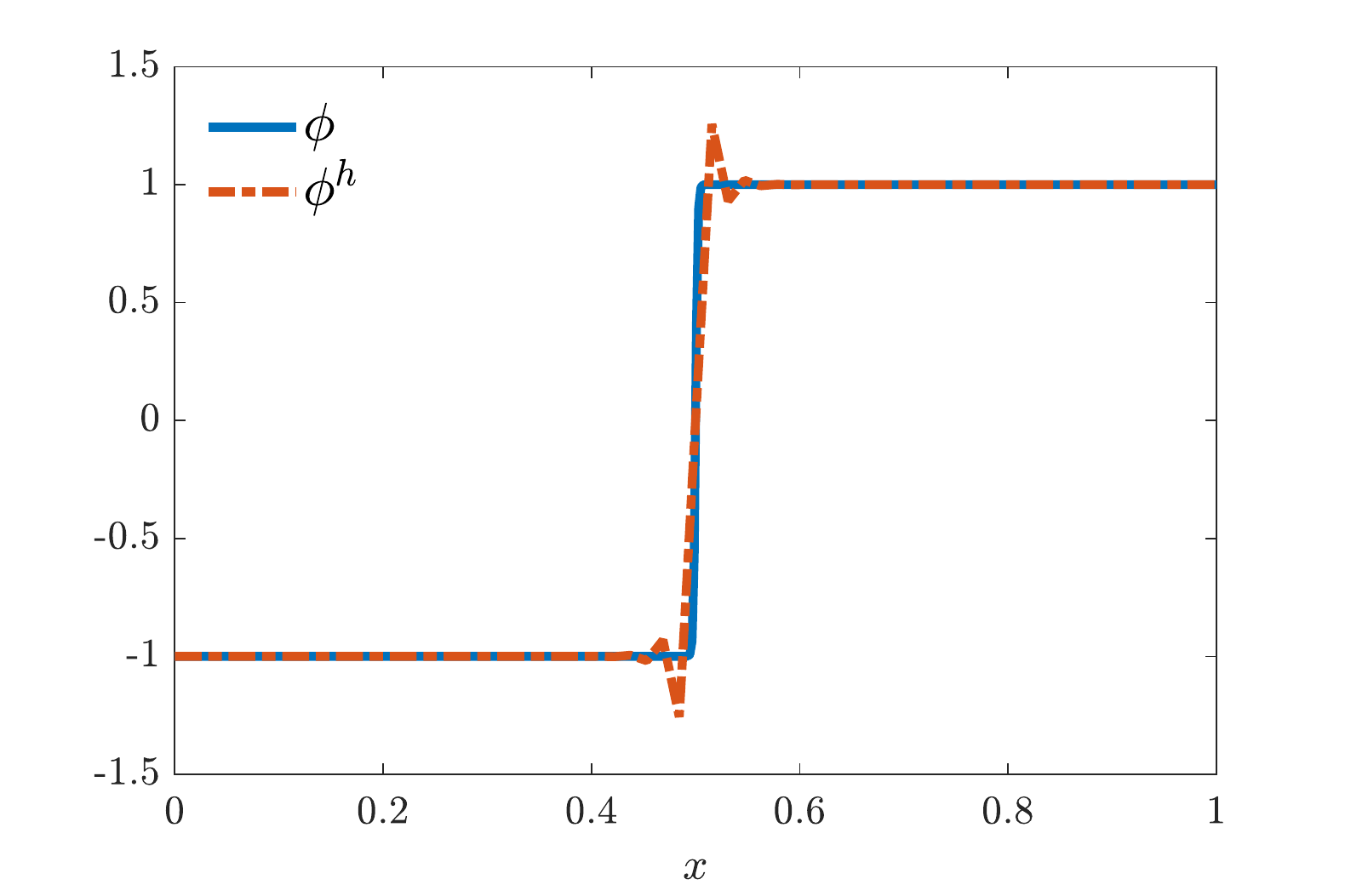}
\caption{$n_{\rm el} =64$.}
\end{subfigure}
\caption{The $L^2$-best approximation $\phi^h \in \mathcal{V}_{D;1,0}^h$ of the smooth step function $\phi = \phi_{0.5}$.}
\label{fig: L2. Continuous Piecewise linears}
\end{figure}

We start off with the case in which the approximation space $\mathcal{V}^h_{p,\alpha}$ consists of continuous functions, i.e. $\alpha \geq 0$. The $L^2$-best approximation $\phi^h \in \mathcal{V}_{D;1,0}^h$ (the space spanned by continuous piecewise linear basis functions) of the smooth step function is illustrated in \cref{fig: L2. Continuous Piecewise linears}.
We observe that the numerical approximation $\phi^h$ contains over- and undershoots near the sharp layer. These oscillations do not vanish when the number of elements is increased. In fact, the over- and undershoots on each side of the discontinuity converge to the value $1-\sqrt{3}/2\approx 0.13$ as the number of elements is increased 
(assuming that the layer is `sufficiently sharp') \cite{foster1991gibbs}. 
We note that the Gibbs phenomenon is often mistakenly interpreted as related to approximation with higher-order basis functions. This example illustrates that this is \textit{not} the case.

\begin{figure}[!ht]
\begin{subfigure}{0.49\textwidth}
\centering
\includegraphics[width=0.95\textwidth]{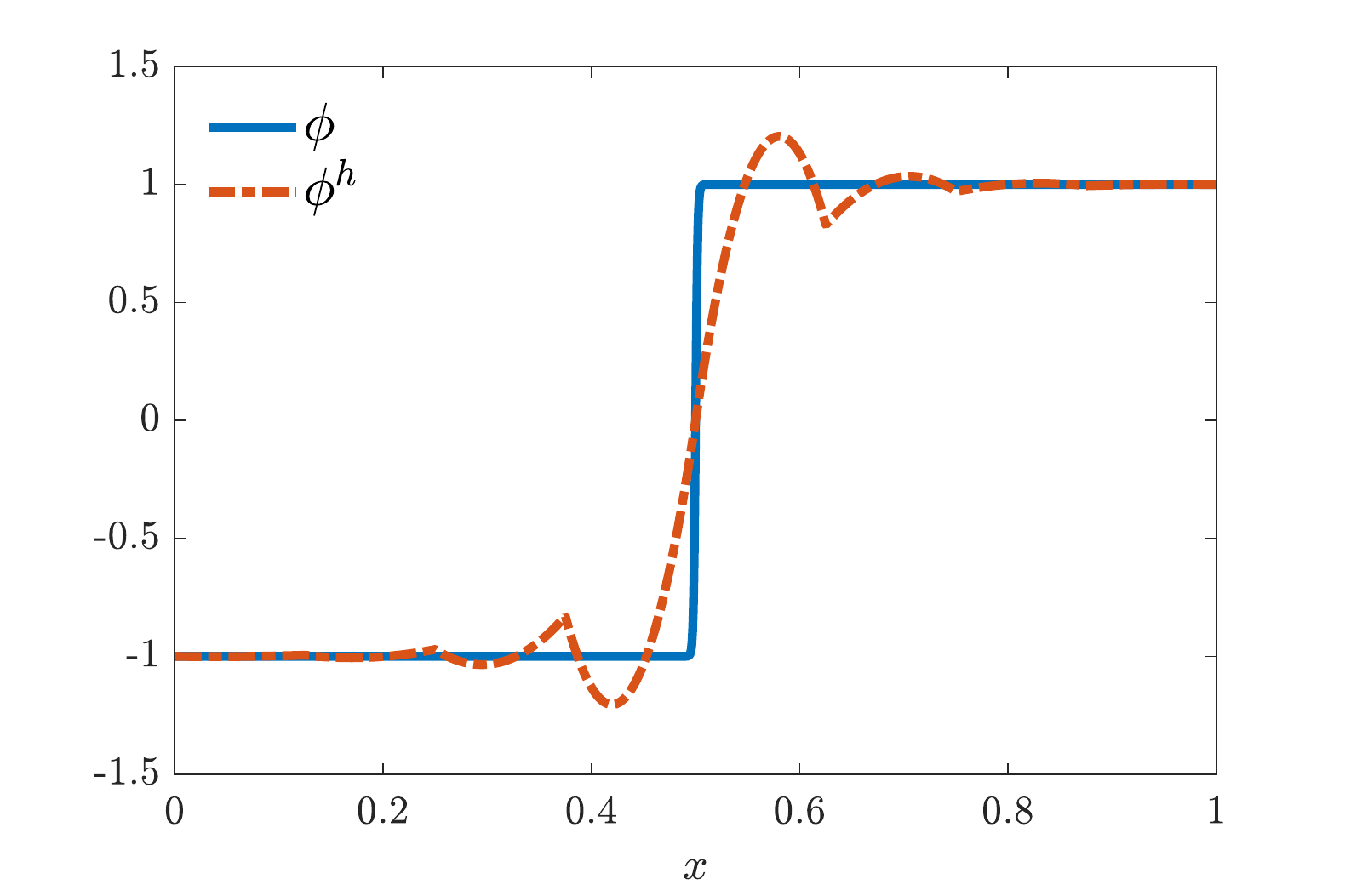}
\caption{$p=2$ and $\alpha = 0$.}
\end{subfigure}
\begin{subfigure}{0.49\textwidth}
\centering
\includegraphics[width=0.95\textwidth]{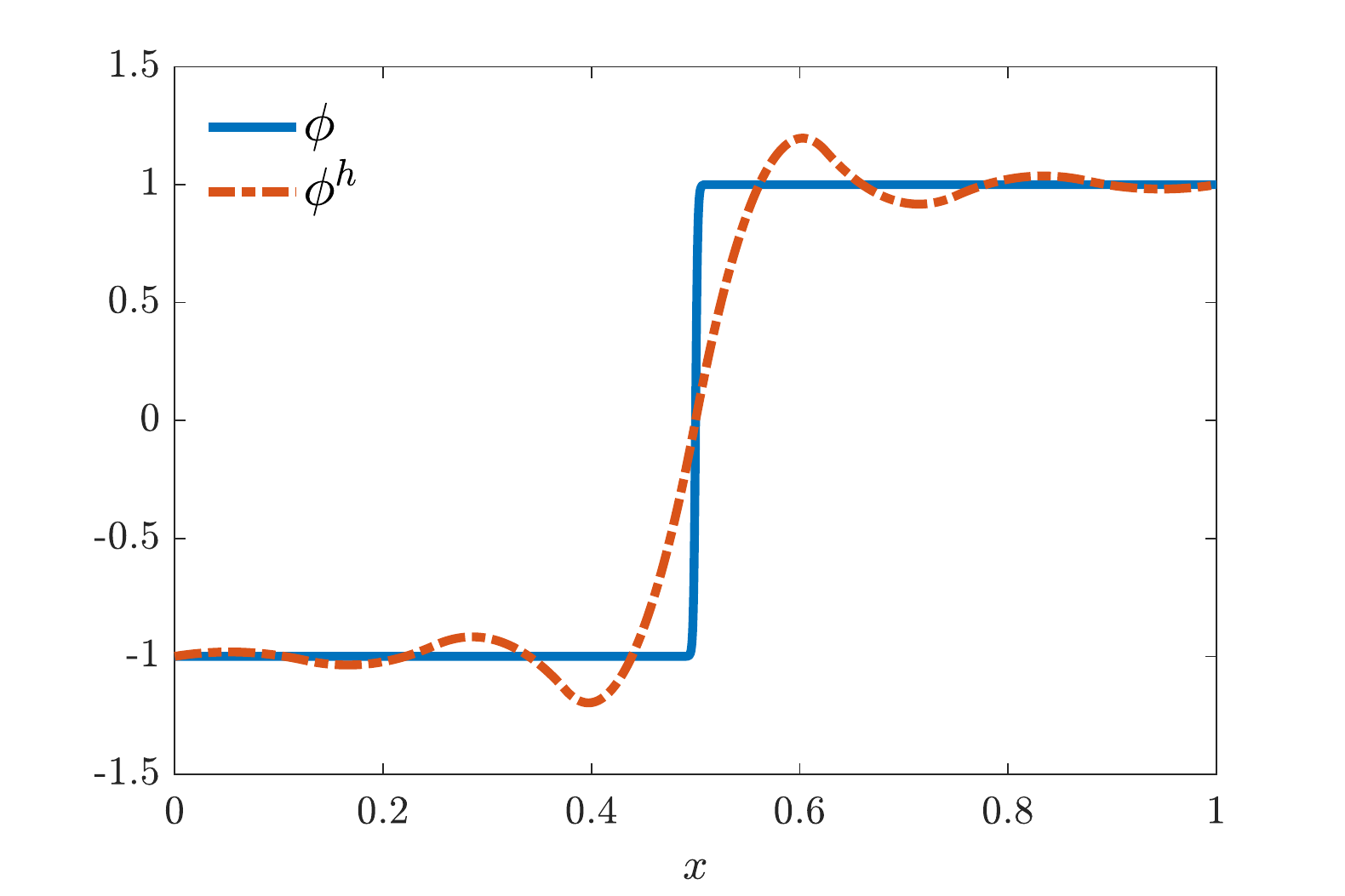}
\caption{$p=2$ and $\alpha = 1$.}
\end{subfigure}
\caption{The $L^2$-best approximations $\phi^h \in \mathcal{V}_{D;2,0}^h$ (left) and $\phi^h \in \mathcal{V}_{D;2,1}^h$ (right), both with $n_{\rm el} =8$, of the smooth step function $\phi = \phi_{0.5}$.}
\label{fig: L2-best approximation. Quadratic Lagrange}
\end{figure}

\cref{fig: L2-best approximation. Quadratic Lagrange} shows the approximations in $\mathcal{V}_{D;2,0}^h$ and $\mathcal{V}_{D;2,1}^h$, the spaces of continuous quadratic finite elements.
The figure shows over- and undershoots of roughly the same magnitude as those for the linear approximation. For the quadratic B-splines the over- and undershoots on each side of the discontinuity converge with the number of elements to a value of approximately $0.10$ \cite{richards1991gibbs}. The Gibbs phenomenon persists when the polynomial order $p$ of the maximum regularity B-spline basis functions is increased. 
Moreover, the magnitude of the over- and undershoots converges to the same value as that of a truncated Fourier series. This value is  approximately $0.09$ and is known as the \textit{Gibbs constant}.

\begin{remark}[Different degrees of freedom]\label{rmk: dofs}
It is important to realize that  the number of degrees-of-freedom (dofs) is significantly different for results with the same number of elements and polynomial degree $p$ but with different regularity $\alpha$. In this situation we have $21$ dofs for $\phi^h \in \mathcal{V}^h_{D,2,0}$ and only $10$ dofs for $\phi^h \in \mathcal{V}^h_{D,2,1}$ (the count includes boundary dofs).
\end{remark}

\begin{figure}[!ht]
\begin{subfigure}{0.49\textwidth}
\centering
\includegraphics[width=0.95\textwidth]{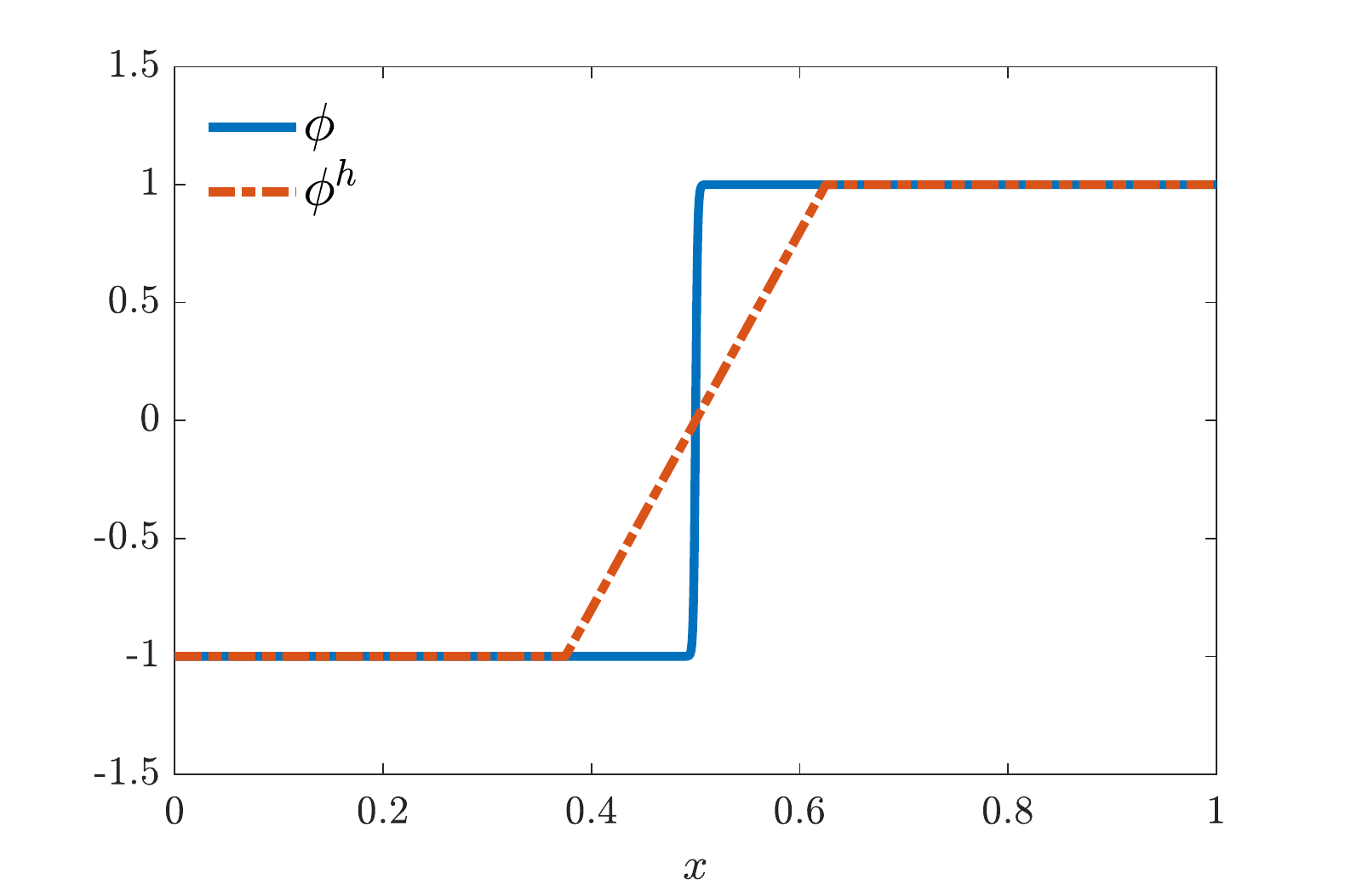}
\caption{$a = 0.5$.}
\end{subfigure}
\begin{subfigure}{0.49\textwidth}
\centering
\includegraphics[width=0.95\textwidth]{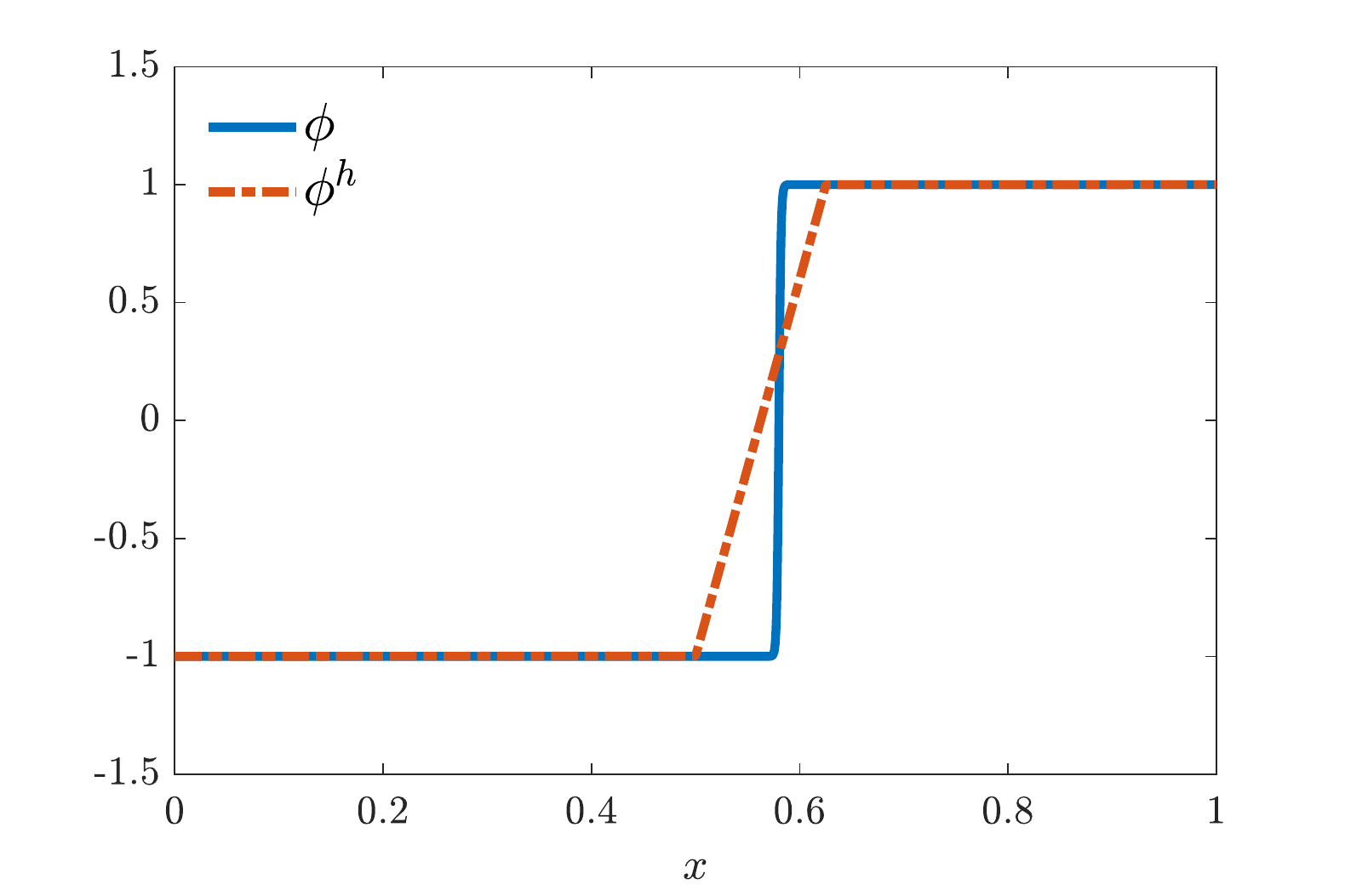}
\caption{$a = 0.58$.}
\end{subfigure}
\caption{The $H_0^1$-best approximations $\phi^h \in \mathcal{V}_{D;1,0}^h$, with $n_{\rm el} =8$, of the smooth step function $\phi = \phi_a$ for different $a$.}
\label{fig: H01-best approximation. Piecewise linears. 8 elements}
\end{figure}

In \cref{fig: H01-best approximation. Piecewise linears. 8 elements}, we visualize the $H_0^1$-best approximation $\phi^h \in \mathcal{V}_{D;1,0}^h$ (the space spanned by continuous piecewise linear basis functions) for $a=0.5$ and $a=0.58$. We observe nodally exact numerical approximations for both cases. The combination of the linear basis functions with the nodal exactness implies that the numerical approximations are free of over- and undershoots, i.e. the Gibbs phenomenon is not present.

\begin{figure}[!ht]
\begin{subfigure}{0.49\textwidth}
\centering
\includegraphics[width=0.95\textwidth]{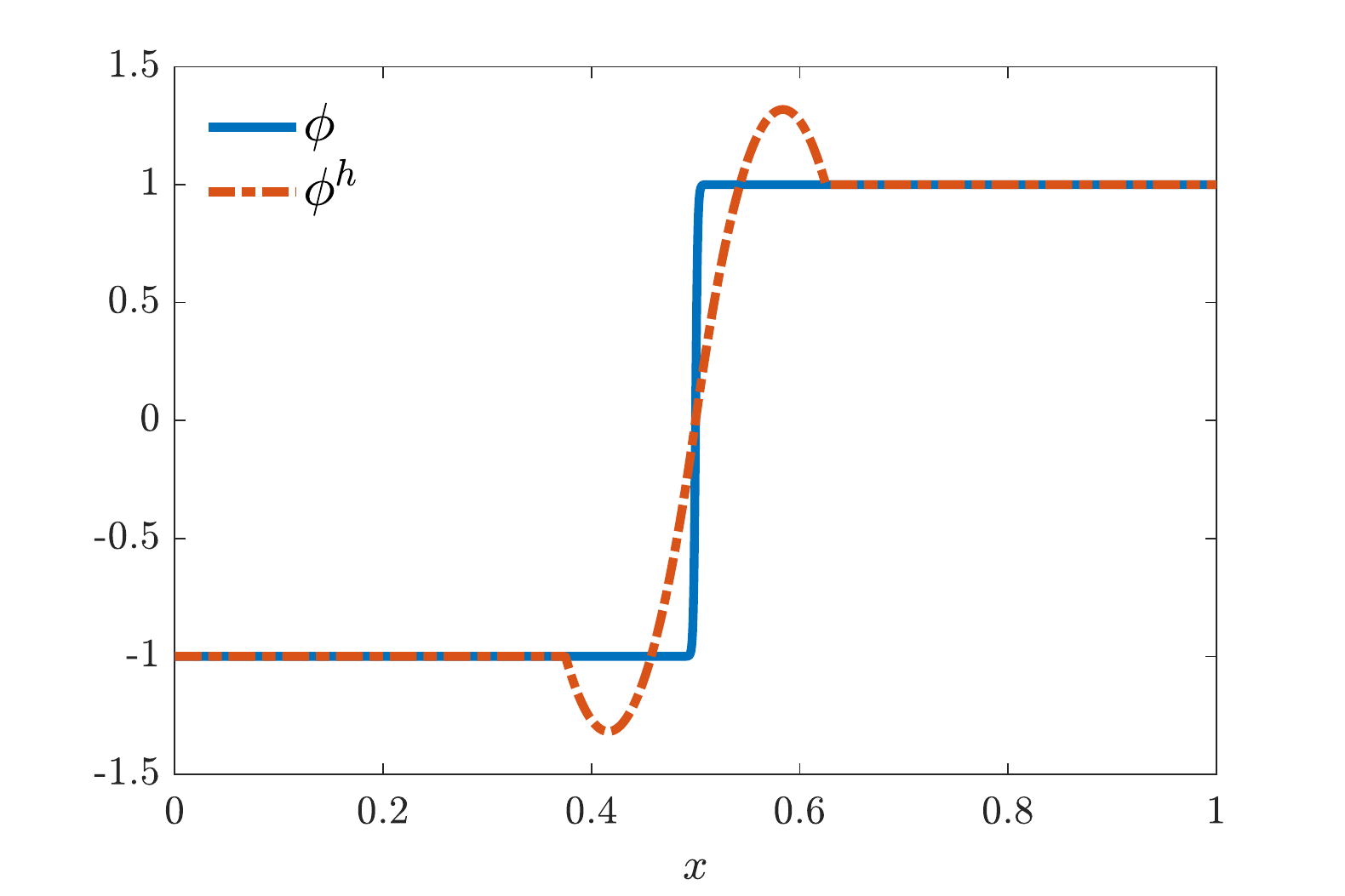}
\caption{$p=2$ and $\alpha = 0$.}
\end{subfigure}
\begin{subfigure}{0.49\textwidth}
\centering
\includegraphics[width=0.95\textwidth]{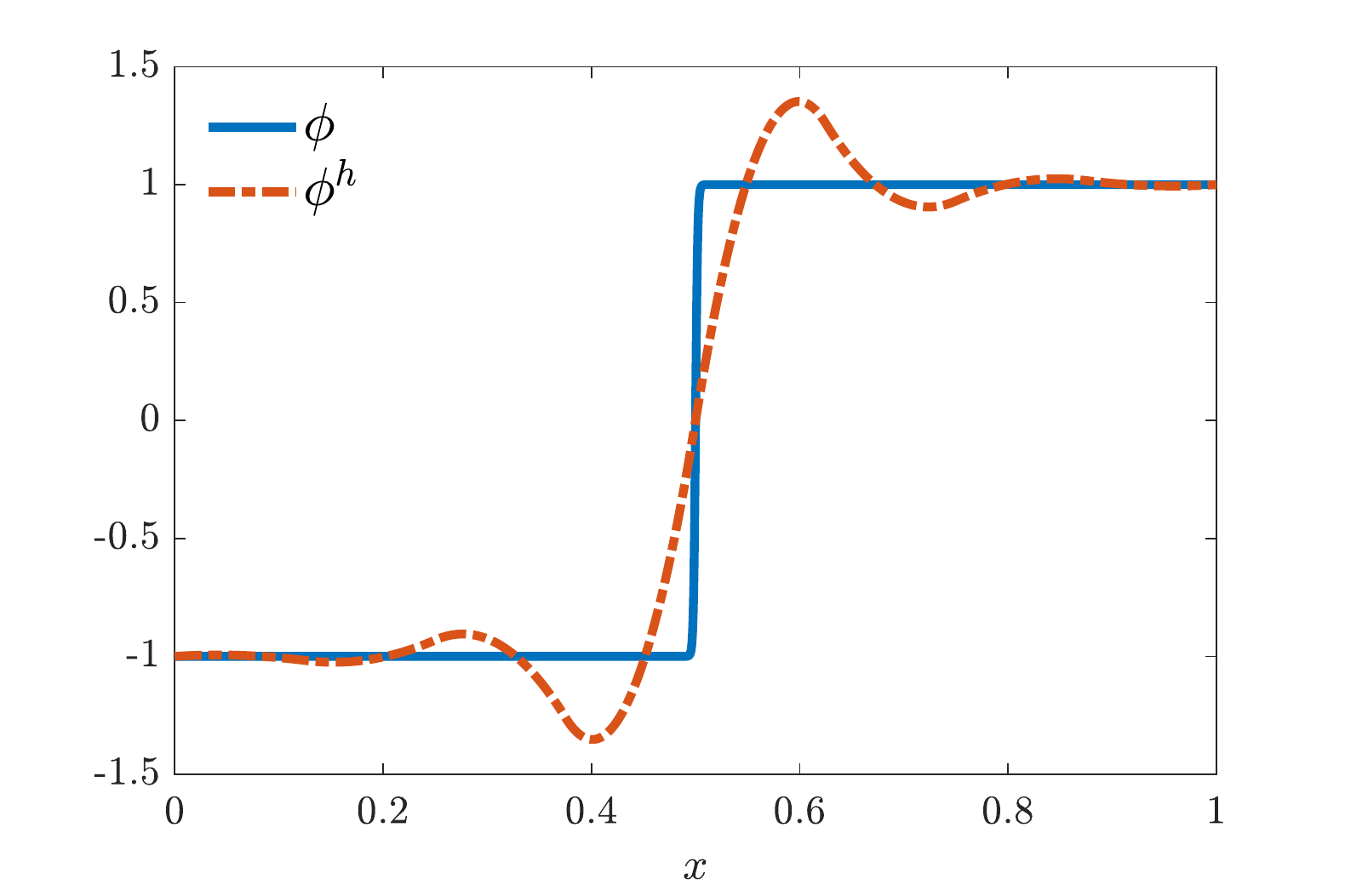}
\caption{$p=2$ and $\alpha = 1$.}
\end{subfigure}
\caption{The $H_0^1$-best approximation $\phi^h \in \mathcal{V}_{D;2,0}^h$ (left) and $\phi^h \in \mathcal{V}_{D;2,1}^h$ (right), with $n_{\rm el} =8$, of the smooth step function $\phi = \phi_{0.5}$.}
\label{fig: H01-best approximation. Quadratics}
\end{figure}

In \cref{fig: H01-best approximation. Quadratics}, we illustrate the $H_0^1$-best approximations for the approximation spaces of quadratic finite elements $\mathcal{V}_{D;2,0}^h$ and $\mathcal{V}_{D;2,1}^h$. We observe over- and undershoots for both the quadratic Lagrange polynomials and the B-spline functions. Again, these oscillations persist with mesh refinement. For the case of the Lagrange basis functions, we have nodal exactness at element boundary nodes, but the monotonicity property is lost in the element interiors \cite{HugSan06}.

\begin{lemma}[Nodal interpolant]\label{lem: monotonicty H01}
The $H_0^1$-best approximation in the space $\mathcal{V}_{D;p,0}^h$ in one dimension is nodally interpolatory at the element boundary nodes. For linear elements ($p=1$) this best approximation is monotonicity preserving, while for higher-order basis functions ($p>1$) monotonicity inside the elements is in general lost.
\end{lemma}

\begin{figure}[!ht]
\begin{subfigure}{0.49\textwidth}
\centering
\includegraphics[width=0.95\textwidth]{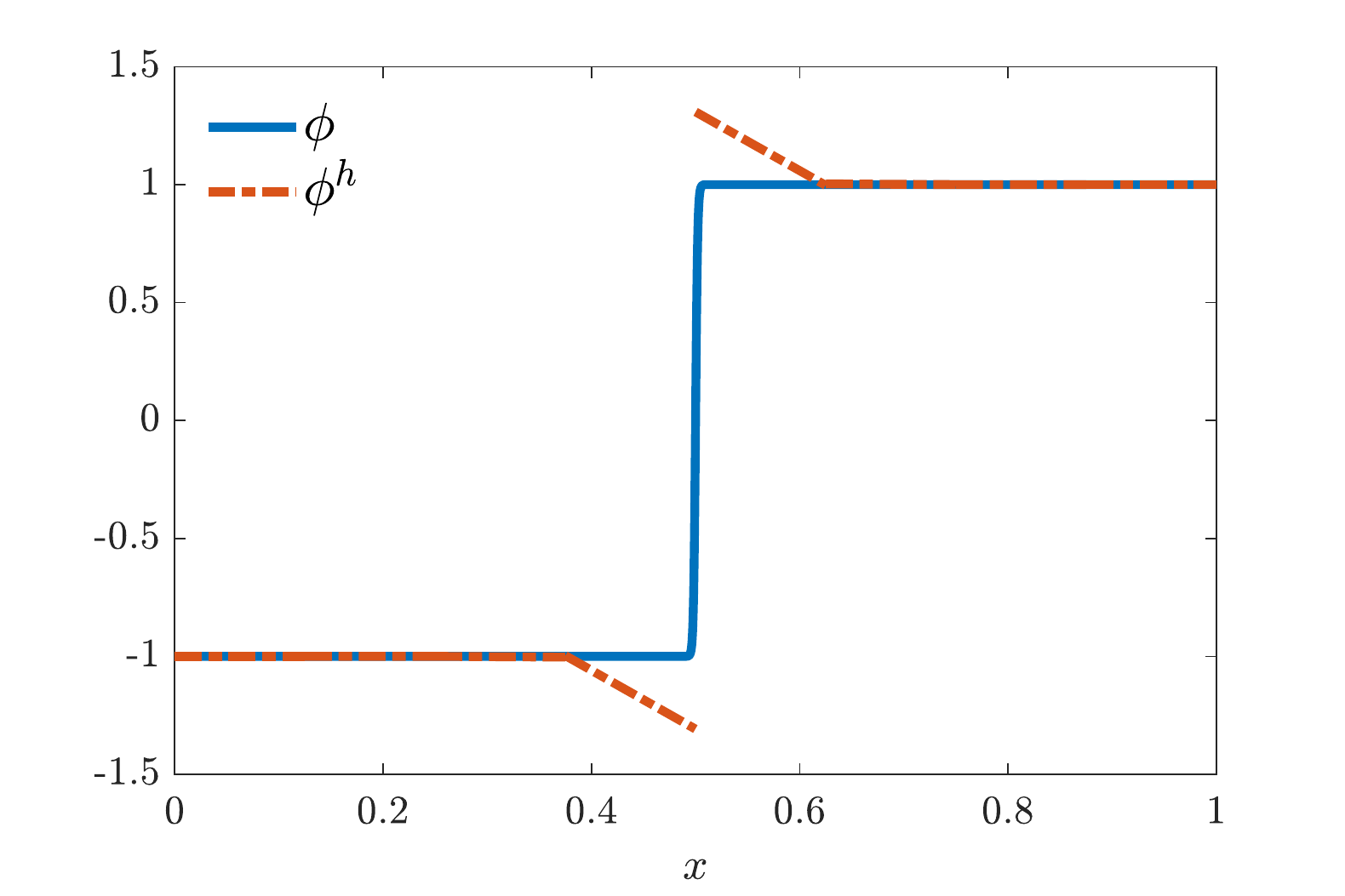}
\caption{$p=1$.}
\end{subfigure}
\begin{subfigure}{0.49\textwidth}
\centering
\includegraphics[width=0.95\textwidth]{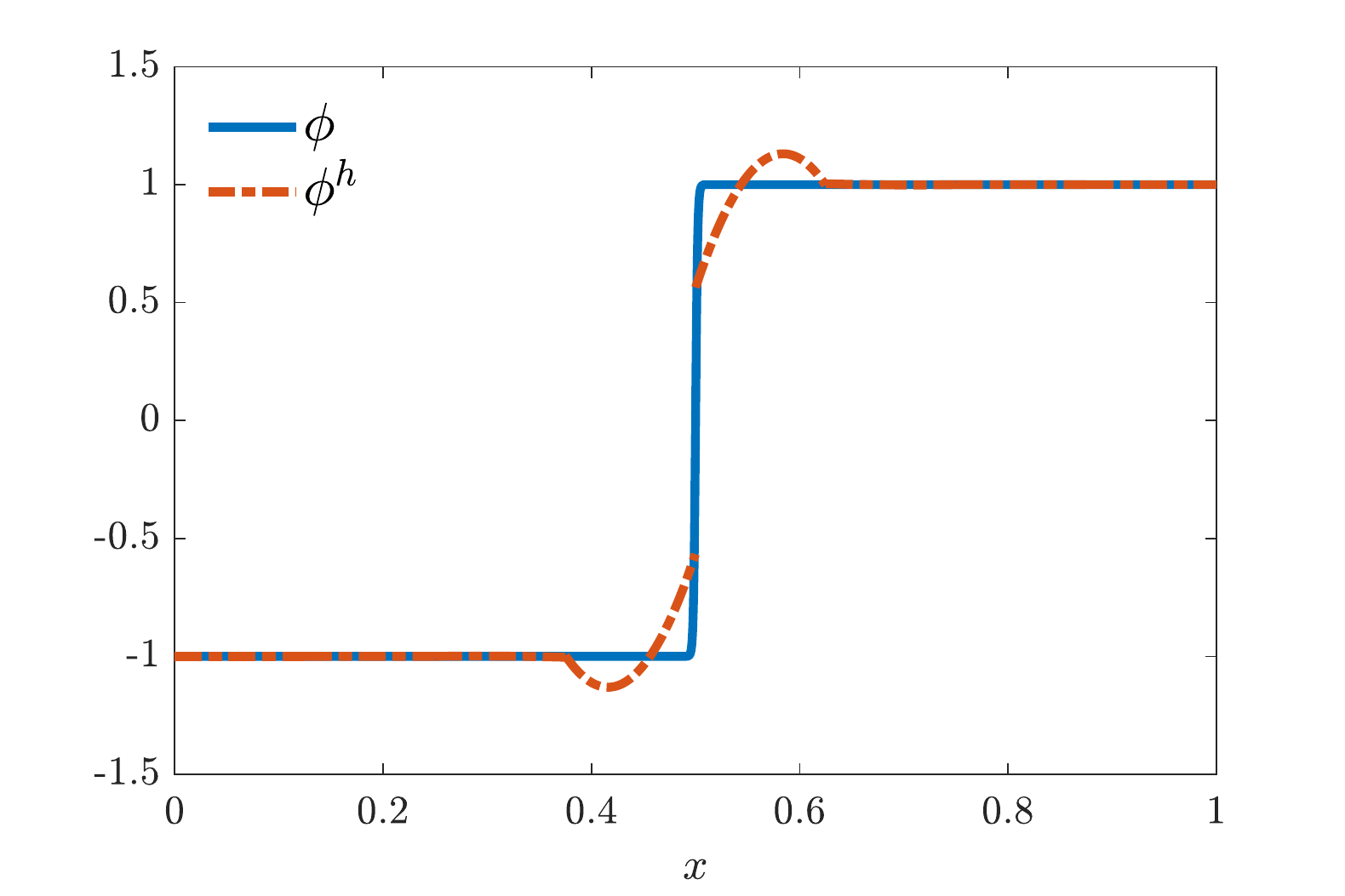}
\caption{$p=2$.}
\end{subfigure}
\caption{The interior penalty-best approximation $\phi^h \in \mathcal{V}_{D;p,-1}^h$, with $n_{\rm el} =8$, of the smooth step function $\phi = \phi_{0.5}$, for different $p$.}
\label{fig: IP. Linears an Quadratics}
\end{figure}

Next, we turn our attention to discontinuous approximations, i.e. $\alpha = -1$. We select as penalty parameter $\eta = 6(p+1)^2/h$. The interior penalty-best approximation of the smooth step function is illustrated in \cref{fig: IP. Linears an Quadratics} for the approximation spaces $\mathcal{V}^h_{D,1,-1}$ (discontinuous piecewise linears) and $\mathcal{V}^h_{D,2,-1}$ (discontinuous piecewise quadratics). For this best approximation problem, we see that the numerical approximations $\phi^h$ contain over- and undershoots near the sharp layer; the nodal exactness of the $H_0^1$-best approximation problem is not inherited. Furthermore, we observe that the average of the approximation $\phi^h$ at the element boundaries coincides with the value of $\phi$. This is a property of the interior penalty projector \cite{stoter2021nitsche}.

\begin{proposition}[Vanishing average error on element boundaries]\label{lem: IP prop}
  The interior penalty best approximation $\phi^h \in \mathcal{V}^h_{D,p,-1}$ of $\phi$ in one-dimension satisfies the property:
  \begin{align}
      \avg{\phi^h - \phi}\Big|_{\Gamma^0} = 0,
  \end{align}
  for all polynomial orders $p$.
\end{proposition}

\subsection{The Gibbs phenomenon in two dimensions}\label{sec: review Gibbs 2D}

We consider best approximations of a two-dimensional step function $\phi \in L^2(\Omega)$ on the square domain $\Omega = (-1,1)^2$:
\begin{align*}
    \phi=\phi(x,y) = \left\{
      \begin{matrix}
        1 & x-y > 0\\
        -1 & x-y < 0
      \end{matrix}
    \right. \,.
\end{align*}
Again, we work with a smooth approximation $\phi \in H^1(\Omega)=:\mathcal{V}$ of the step function:
\begin{align*}
    \phi=\phi^{\epsilon}(x,y) = \tanh\left(\dfrac{y-x}{\epsilon}\right),
\end{align*}
where $\epsilon \ll 1$ is a smoothing parameter.

\begin{figure}
\begin{subfigure}{0.49\textwidth}
\includegraphics[width=0.95\textwidth]{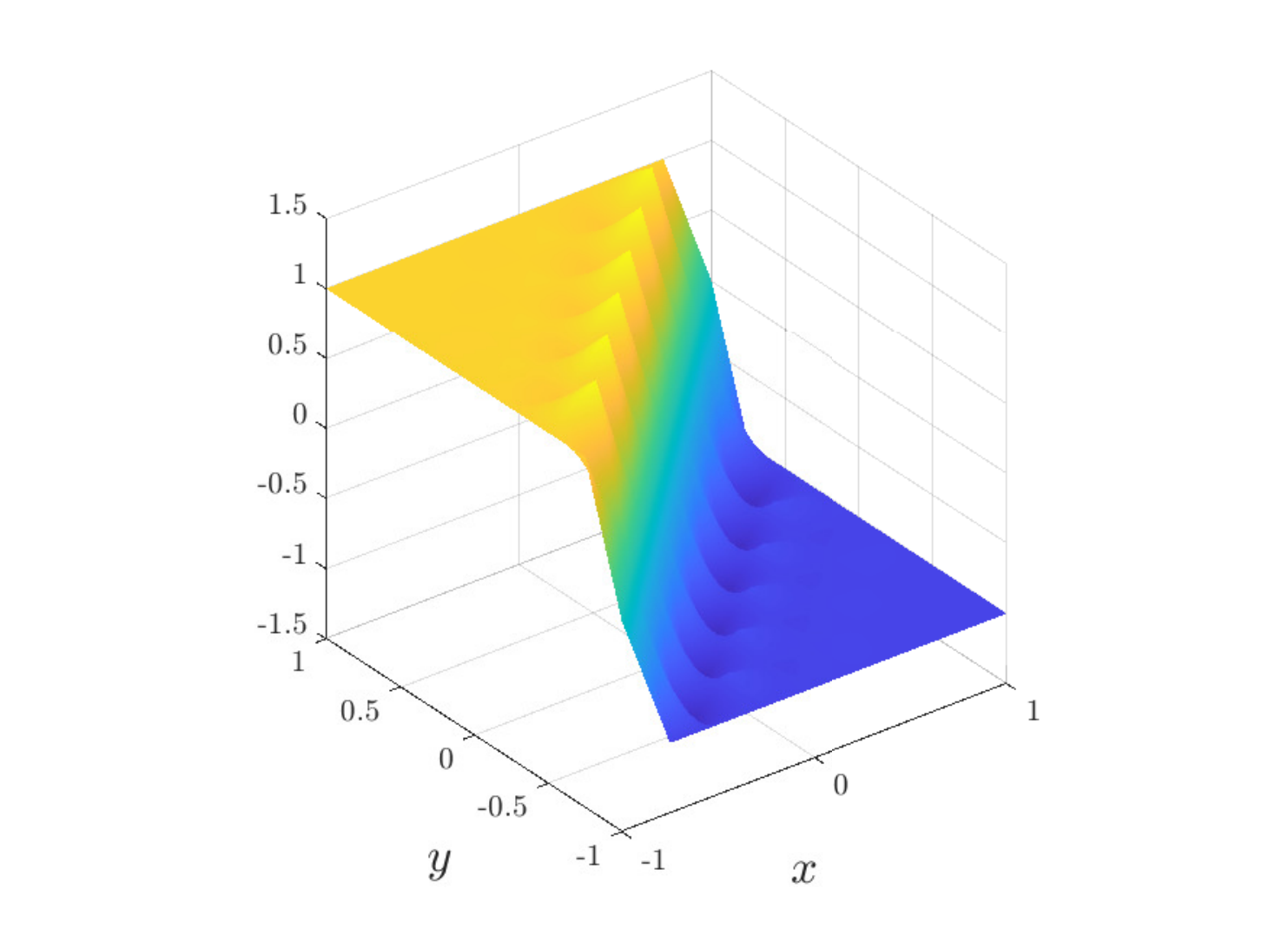}
\caption{$p=1$ and $\alpha = 0$.}
\end{subfigure}
\begin{subfigure}{0.49\textwidth}
\includegraphics[width=0.95\textwidth]{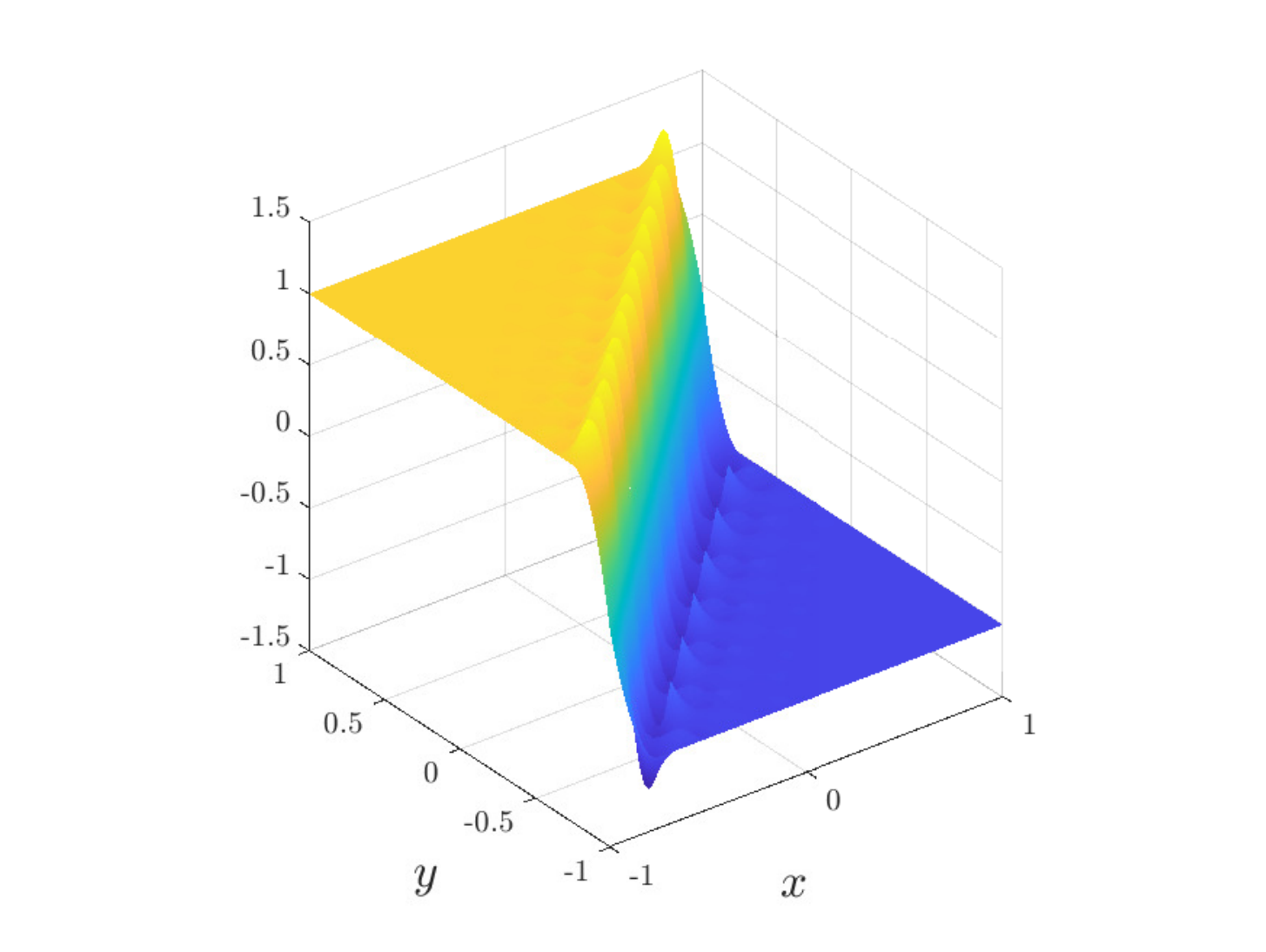}
\caption{$p=2$ and $\alpha = 0$.}
\end{subfigure}
\begin{center}
\begin{subfigure}{0.49\textwidth}
\includegraphics[width=0.95\textwidth]{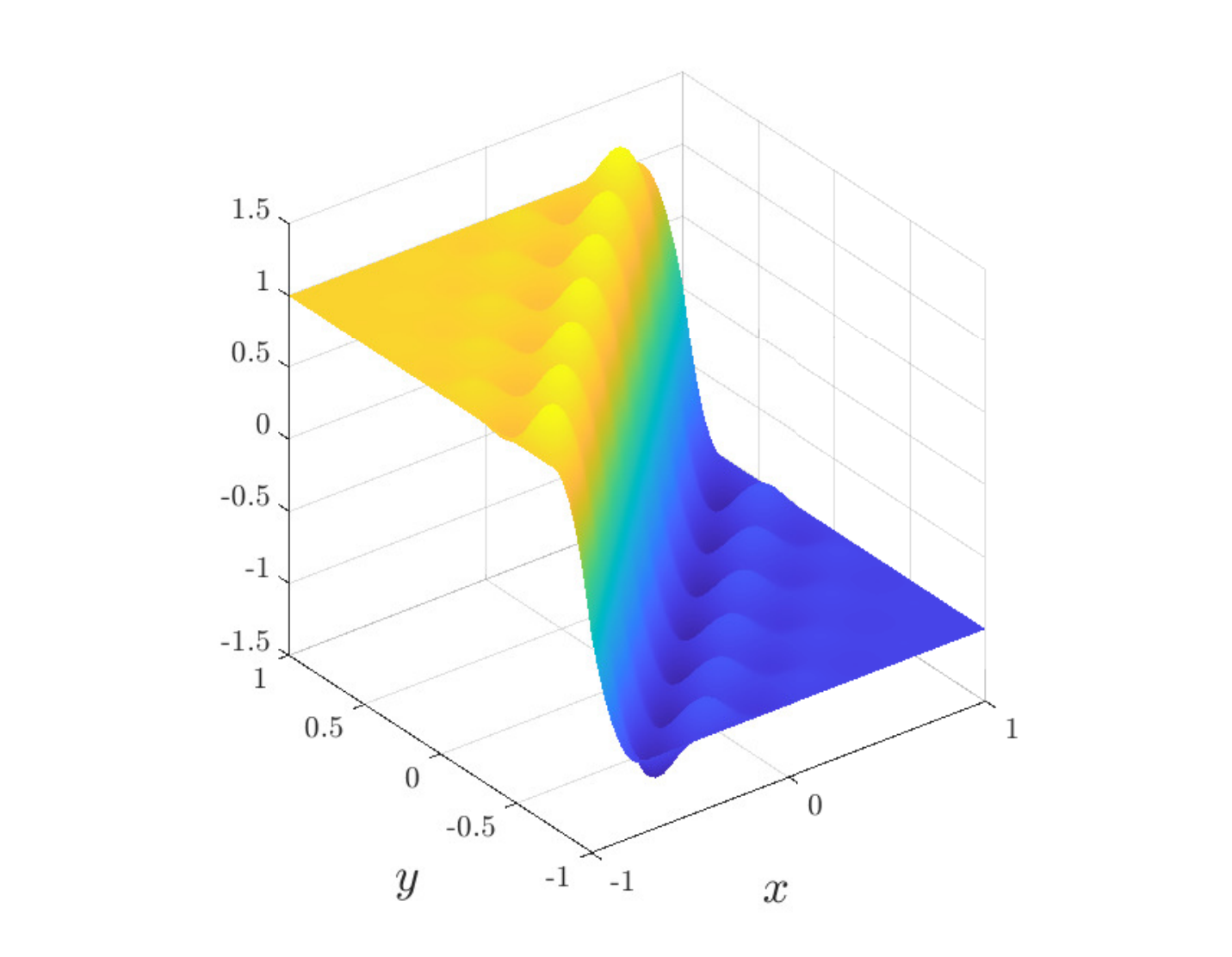}
\caption{$p=2$ and $\alpha = 1$.}
\end{subfigure}
\end{center}
\caption{The $H_0^1$-best approximation $\phi^h \in \mathcal{V}_{D;p,\alpha}^h$, with $n_{{\rm el}}=8\times8$), of the smooth two-dimensional step function, for different $p$ and $\alpha$.}
\label{fig: 2D H01-best approximations}
\end{figure}

Analogous to the one-dimensional case, we begin with approximation spaces $\mathcal{V}^h_{D,p,\alpha}$ consisting of continuous functions ($\alpha \geq 0$). Recall from the one-dimensional case 
that the $L^2$-best approximation contains over- and undershoots. This is also the case in higher dimensions, and we omit the visualization. We display the $H_0^1$-best approximation for the continuous finite element approximation spaces $\mathcal{V}^h_{D,1,0}, \mathcal{V}^h_{D,2,0}$ and $\mathcal{V}^h_{D,2,1}$ in \cref{fig: 2D H01-best approximations}. We observe the occurrence of over- and undershoots for each of the approximations.

\begin{figure}[!ht]
\begin{center}
\begin{subfigure}{0.49\textwidth}
\includegraphics[width=0.77\textwidth]{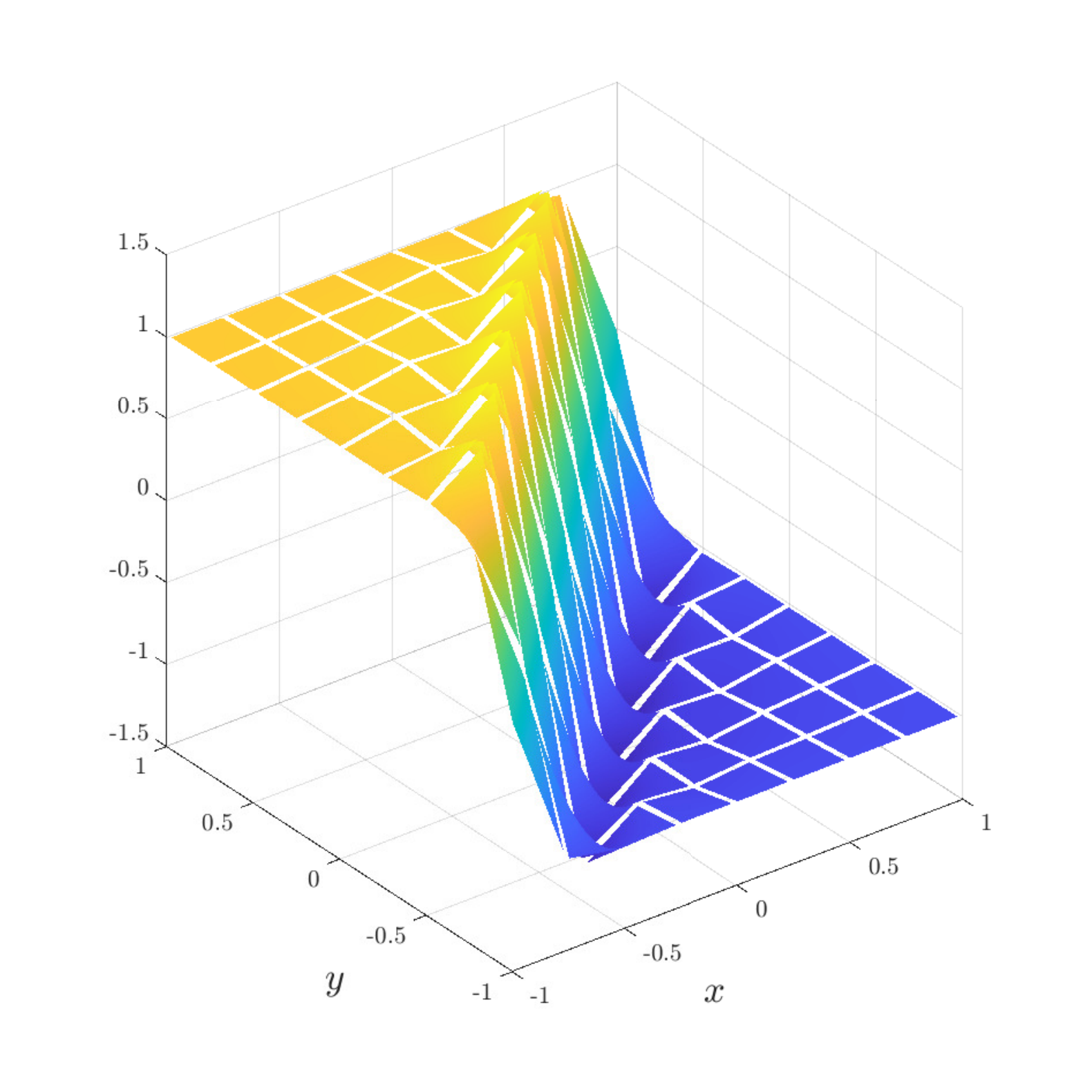}
\caption{$p = 1$}
\end{subfigure}
\begin{subfigure}{0.49\textwidth}
\includegraphics[width=0.77\textwidth]{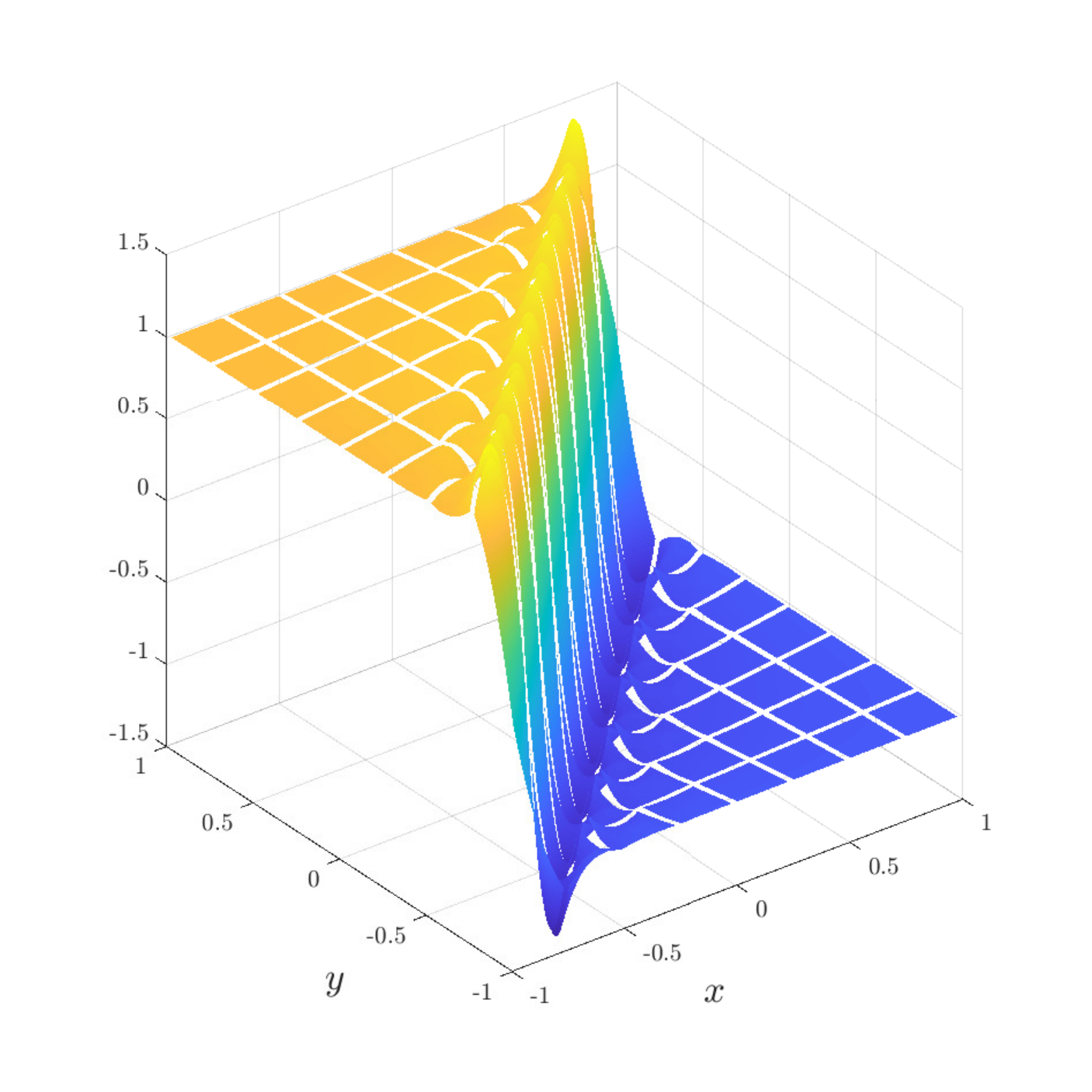}
\caption{$p = 2$}
\end{subfigure}
\caption{The interior penalty-best approximation $\phi^h \in \mathcal{V}_{D;p,-1}^h$, with $n_{{\rm el}}=8\times8$, of the smooth two-dimensional step function, for different $p$.}
\label{fig: 2D IP-best approximations}
\end{center}
\end{figure}

Next, we consider discontinuous approximations ($\alpha = -1$). Analogous to the one-dimensional case, we solely consider the interior penalty-best approximation. We select as a penalty parameter $\eta = 2(2 p +1)(2p+2)/h$. In \cref{fig: 2D IP-best approximations}, we visualize the interior penalty best approximations for approximation spaces $\mathcal{V}^h_{1,-1}$ and $\mathcal{V}^h_{2, -1}$, i.e. for linear and quadratic discontinuous basis functions. Again, we observe over- and undershoots of the finite element approximation in both cases.

\section{Eliminating the Gibbs phenomenon in one dimension}\label{sec: Gibbs constraints 1D}

In this section we present the constraints for the elimination of the Gibbs phenomenon in finite element spaces in one dimension. We call these constraints the \textit{Gibbs constraints}. To this purpose, we first describe the construction of our proposed Gibbs constraints in general approximation spaces in \cref{subsec: Gibbs 1D: general}, and present the properties of the constraints in \cref{subsec: Properties Gibbs}. Next, we advance the discussion to best approximation problems in \cref{subsec: Gibbs 1D: opt}. Finally, in \cref{subsec: Gibbs 1D fem} we apply the Gibbs constraints to finite element spaces of arbitrary continuity, and perform numerical experiments. 

\subsection{Gibbs constraints}\label{subsec: Gibbs 1D: general}
One of the challenges of dealing with the Gibbs phenomenon is the uncertainty in the level of locality that is required to identify the phenomenon. Pointwise evaluations of functions, or their derivatives, carry insufficient information to be able to infer the occurrence of the Gibbs phenomenon. On the other hand, the information that may be deduced from global evaluations, such as global integrals, is too coarse-grained to establish the existence of spurious oscillations on a local scale. In this subsection, we construct a constraint for the elimination of the Gibbs phenomenon on a \textit{given subdomain}. The selection of the subdomains then remains an important matter, which we discuss extensively in \cref{subsec: Gibbs 1D fem} in the context of finite element approximations.

Consider the one-dimensional simply connected domain $\Omega \subset \mathbb{R}$, and let $\phi^*: \Omega \rightarrow \mathbb{R}$ denote an approximation of $\phi: \Omega \rightarrow \mathbb{R}$. The main ingredient in the elimination of the Gibbs phenomenon relies on the \textit{fundamental theorem of Lebesgue integral calculus}, which we recall here.

\begin{theorem}[Fundamental theorem of Lebesgue integral calculus]\label{thm: fundamental calculus}
Let $\theta: \Omega \rightarrow \mathbb{R}$ be an absolutely continuous function, then $\theta$ is differentiable almost everywhere and for each $\omega = [x_L,x_R] \subset \Omega$ we have
\begin{align}\label{eq: fundamental theorem}
    \displaystyle\int_{x_L}^{x_R} {\rm D}\theta ~{\rm d}x = \theta_R-\theta_L,
\end{align}
with trace equalities $\theta(x_L) = \theta_L$ and $\theta(x_R) = \theta_R$, and ${\rm D}\theta \in L^1(\Omega)$.
\end{theorem}
The fundamental theorem communicates that the trace values (the right-hand side in \eqref{eq: fundamental theorem}) are controlled by the integral. Still, the trace values do not provide any information on the oscillatory behavior of the function $\theta$ inside $\omega$. In contrast, the \textit{total variation} is a concept that does incorporate this.

\begin{definition}[Total variation]\label{def: total variation}
Let $\theta: \Omega \rightarrow \mathbb{R}$ be a given function, and let\\ $P=\left\{x_L = x_0, x_1, \dots, x_{N-1}, x_N = x_R\right\}$ denote a partition of $\omega = [x_L,x_R] \subset \Omega$. The variation of $\theta$ with respect to partition $P$ is defined as:
\begin{align}
    V_{\omega,P}(\theta):= \sum_{i=0}^{N} |\theta(x_{i+1})-\theta(x_i)|.
\end{align}
Denote by $\mathfrak{P}$ the set of all possible partitions of $\omega$. The total variation is given by:
\begin{align}
    V_{\omega}(\theta):= \sup_{P \in \mathfrak{P}} V_{\omega,P}(\theta).
\end{align}
\end{definition}
The total variation $V_{\omega}(\theta)$ represents a measure of the fluctuations of $\theta$ on $\omega$. A function $\theta$ with the property $V_\omega(\theta)<\infty$ is said to have \textit{bounded variation} and we write $\theta \in {\rm BV}(\omega)$. An absolutely continuous function has bounded variation. 
In case~$\theta$ is a continuous differentiable function, the total variation $V_\omega(\theta)$ may be evaluated as follows.

\begin{lemma}[Total variation continuous differentiable function]\label{lem: TV of C1 function}
Let $\theta \in \mathcal{C}^1(\Omega)$. The total variation of $\theta$ on $\omega \subset \Omega$ is given by:
\begin{align}
    V_{\omega}(\theta) = \displaystyle\int_{\omega} \left|{\rm D}\theta\right| ~{\rm d}x.
\end{align}
\end{lemma}

If the function $\theta$ is only piecewise continuously differentiable, the expression in \cref{lem: TV of C1 function} needs to be augmented with jump terms, as is expressed in \cref{lem: TV of piecewise C1 function}.

\begin{lemma}[Total variation piecewise-continuous function]\label{lem: TV of piecewise C1 function}
Let $\theta \in {\rm BV}(\Omega)$ be a function that has a continuous derivative on each $(a_i,a_{i+1}) \subset \Omega, i=0, \dots, M$ and jump discontinuities at $a_i, i=1, \dots, M$. Denote with $a_{i}^{-}, a_{i}^+$ the left and right limits of the discontinuity at $a_i$. The total variation of $\theta$ on $\omega \subset \Omega$ is given by:
\begin{align}
    V_{\omega}(\theta) = \displaystyle\sum_{i=0}^{M} \displaystyle\int_{a_i}^{a_{i+1}} \left|{\rm D}\theta\right| ~{\rm d}x + \displaystyle\sum_{i=1}^M \left|\theta(a_i)-\theta(a_i^-)\right| +  \left|\theta(a_i^+)-\theta(a_i)\right|.
\end{align}
\end{lemma}

Next, we note that $V_\omega$ is convex and satisfies a homogeneity property. 
\begin{proposition}[Convexity and homogeneity of the total variation]\label{prop: homogeneity}
  Given the same assumptions on $\theta$ as in \cref{lem: TV of piecewise C1 function}, the functional $V_\omega$ is convex and satisfies the homogeneity property:
  \begin{align}
    {\rm d}V_\omega(\theta)(\theta) = V_\omega(\theta).  
  \end{align}
\end{proposition}

\begin{remark}
  A convex functional that satisfies the homogeneity property is\\ termed a \textit{variation entropy} \cite{ten2019variation}. Variation entropy functionals form the basis of an entropy stability theory for (hyperbolic) conservation laws called {variation entropy theory}. 
\end{remark}

The total variation of a monotonic function is related to its trace values in the following way.

\begin{lemma}[Total variation monotonic function]\label{lem: total variation bounds}
 Let $\omega = [x_L,x_R] \subset \Omega$ and let $\mathcal{C}^0_D(\Omega)$ denote the space of continuous functions $\theta$ that satisfy trace equalities $\theta(x_L) = \theta_L$ and $\theta(x_R) = \theta_R$. For $\theta \in \mathcal{C}^0_D(\Omega)$ we have 
 \begin{align}
    V_\omega(\theta) \geq  \left|\theta_R-\theta_L\right|,
\end{align}
where equality only holds when $\theta$ is monotonic.
\end{lemma}
This is an important ingredient in the design of the Gibbs constraints. It communicates that the jump of the trace values are controlled by the total variation. 

Let us introduce the quantity
\begin{align}
\mathscr{V}_{\phi,\omega}(\phi^*):=V_{\omega}(\phi^*) - V_\omega(\phi).
\end{align}
With the aim of eliminating the Gibbs phenomenon, one could introduce the constraint $\mathscr{V}_{\phi,\omega}(\phi^*) \leq 0$. The disadvantage is that total variation does not take into account the sign of its argument:
\begin{align}
    V_{\omega}(-\theta) = V_{\omega}(\theta).
\end{align}
As a consequence, we have:
\begin{align}
\mathscr{V}_{\phi,\omega}(-\phi)=0.
\end{align}

We now generalize the element-based definition of the subdivision of domain $\Omega$, \eqref{eq: def tilde Omega fem}, to the union of ($J$) disjoint general subdomains, i.e.
\begin{align}
  \tilde{\Omega} := \displaystyle\bigcup_{j=1}^J \omega_j,
\end{align}
with $\omega_{j}\cap \omega_k = \emptyset$ for $j\neq k$. Furthermore, we redefine the broken space \eqref{eq: broken space fem} for this subdivision as:
\begin{align}
  H^1(\tilde{\Omega}) := \left\{v \in L^2(\Omega): v|_{\omega} \in H^1(\omega) \text{ for all } \omega \in \mathcal{T}_\omega \right\},
\end{align}
with $\mathcal{T}_\omega$ the collection of $\omega_j, j = 1, \dots J$.

Motivated by the above observation, we wish to find approximations that, besides bounding the size of the jump of the approximation, also carry information about the direction of the analytical solution $\phi$. To this purpose we now propose the \textit{Gibbs functional} and the associated \textit{Gibbs constraint}.

\begin{definition}[Gibbs functional]\label{def: Gibbs 1D}
The Gibbs functional of the function $\phi^* \in H^1(\tilde{\Omega})$, on $\omega=\omega_j$ (for some $j$) and with respect to the given function $\phi \in H^1(\Omega)$, is defined as:
\begin{align}\label{eq: large g}
    \mathscr{G}_{\phi,\omega}(\phi^*) := \displaystyle\int_{\omega} g_{\phi}(\phi^*)~{\rm d}x,
\end{align}
where the functional $g_{\phi}$ is defined as:
\begin{align}\label{eq: small g}
    g_{\phi}(\phi^*):= \left| {\rm D}\phi^* \right| - {\rm sgn}\left({\rm D}\phi^*\right){\rm D}\phi= - {\rm sgn}\left({\rm D}\phi^*\right){\rm D}\phi'.
\end{align}
Here ${\rm sgn}$ is the sign function (i.e. ${\rm sgn}(t)=t/|t|$ for $t\neq 0$ and ${\rm sgn}(0)=0$), and $\phi':\Omega \rightarrow \mathbb{R}$ defined as $\phi':= \phi - \phi^*$ denotes the error function.
\end{definition}

\begin{definition}[Gibbs constraint]\label{def: Gibbs constraints 1D}
The Gibbs constraint of the function $\phi^*\in H^1(\tilde{\Omega})$,  on $\omega$ and with respect to the given function $\phi\in H^1(\Omega)$, is defined as:
\begin{align}\label{eq: gibbs constraint}
    \mathscr{G}_{\phi,\omega}(\phi^*) \leq 0.
\end{align}
\end{definition}

Note that the incorporation of the correct sign may be recognized in the Gibbs constraint via the equivalence:   \begin{align}\label{eq: localized constraint}
     \left. \begin{matrix}  \left| {\rm D}\phi^* \right|- \left| {\rm D}\phi\right| \leq 0\\[10pt]
     {\rm sgn}\left( {\rm D}\phi^*\right)- {\rm sgn}\left( {\rm D}\phi\right) =0 \end{matrix} \right\}\Leftrightarrow g_{\phi}(\phi^*)\leq 0,
   \end{align}
the integration of which provides the Gibbs constraint.

The objective is now to search for functions $\phi^*$ as approximations of $\phi$ that satisfy the Gibbs constraint on certain $\omega$.
\subsection{Properties of the Gibbs constraint}\label{subsec: Properties Gibbs}

It is the purpose of this subsection to discuss the properties of the Gibbs constraint and to establish its connection with the well-known concepts of monotonic solutions and the maximum principle.

We have the simple but important property that $\phi$ as an approximation of itself satisfies the Gibbs constraint.
\begin{proposition}[Perfect approximation]\label{prop: perfect approx}
  The Gibbs functional vanishes for a perfect approximation ($\phi^* = \phi$):
\begin{align}
    \mathscr{G}_{\phi,\omega}(\phi) = 0.
\end{align}
\end{proposition}
Furthermore, we have the following lower bound of the Gibbs functional.
\begin{lemma}[Gibbs functional bound]\label{prop: lower bound G}
The Gibbs functional satisfies the lower bound:
\begin{align}
    \mathscr{G}_{\phi,\omega}(\phi^*)\geq \mathscr{V}_{\phi,\omega}(\phi^*).
\end{align}
\end{lemma}

We now proceed with establishing connections between the Gibbs constraint and certain properties of the function approximation. We first provide a characterization of the Gibbs functional.

\begin{lemma}[Characterization Gibbs functional]\label{lem: Gibbs Cont}
Let $\phi^* \in H^1(\tilde{\Omega})$ and $\phi \in H^1(\Omega)$, and let $\omega = [x_L,x_R]\subset \Omega$ be given. Denote the locations of sign changes of ${\rm D}\phi^*$ by $x_i, i=1, \dots, N$ with $x_i<x_{i+1}$. The form of Gibbs functional $\mathscr{G}_{\phi,\omega}$ depends on the sign of ${\rm D}\phi^*$ on $[x_L, x_1]$ and the number of sign changes $N$:
\begin{enumerate}
    \item for $N$ odd and ${\rm D}\phi^*\geq 0$ on $[x_L,x_1]$:
    \begin{align}\label{eq: form G N odd Dphih geq 0}
      \mathscr{G}_{\phi,\omega}(\phi^*)=\phi'(x_L) + 2 \sum_{i=1}^N (-1)^i \phi'(x_i) + \phi'(x_R) ,
    \end{align}
    \item for $N$ odd and ${\rm D}\phi^*\leq 0$ on $[x_L,x_1]$:
    \begin{align}\label{eq: form G N odd Dphih leq 0}
      \mathscr{G}_{\phi,\omega}(\phi^*)=-\phi'(x_L) - 2 \sum_{i=1}^N (-1)^i \phi'(x_i) - \phi'(x_R),
    \end{align}
    \item for $N$ even and ${\rm D}\phi^*\geq 0$ on $[x_L,x_1]$:
    \begin{align}\label{eq: form G N even Dphih geq 0}
      \mathscr{G}_{\phi,\omega}(\phi^*)=\phi'(x_L) + 2 \sum_{i=1}^N (-1)^i \phi'(x_i) - \phi'(x_R),
    \end{align}
    \item for $N$ even and ${\rm D}\phi^* \leq 0$ on $[x_L,x_1]$:
    \begin{align}\label{eq: form G N even Dphih leq 0}
      \mathscr{G}_{\phi,\omega}(\phi^*)=-\phi'(x_L) - 2 \sum_{i=1}^N (-1)^i \phi'(x_i) + \phi'(x_R) ,
    \end{align}
\end{enumerate}
where we recall $\phi'= \phi - \phi^*$.
\end{lemma}

A direct consequence of this characterization is the following lemma.
\begin{lemma}[Interpolatory monotonic approximation]\label{lem: interp mono}
An interpolatory monotonic approximation $\phi^* \in H^1(\tilde{\Omega})$ of $\phi \in H^1(\Omega)$ on $\omega \subset \Omega$ satisfies the Gibbs constraint:
\begin{align}
    \mathscr{G}_{\phi,\omega}(\phi^*) \leq 0.
\end{align}
\end{lemma}

Additionally, an approximation of a monotonic analytical profile that is free of the Gibbs phenomenon satisfies a bound on the trace values.
\begin{lemma}[Approximation of monotonic function]\label{lem: Gibbs ID local}
Suppose that $\phi \in H^1(\Omega)$ is monotonically increasing (decreasing) on $\omega = [x_L,x_R]\subset \Omega$ and the approximation $\phi^*\in H^1(\tilde{\Omega})$ satisfies the Gibbs constraint:
\begin{align}\label{eq: G leq 0}
    \mathscr{G}_{\phi,\omega}(\phi^*) \leq 0,
\end{align}
then the increase (decrease) of $\phi^*$ is bounded by the increase (decrease) of $\phi$, i.e. 
\begin{subequations}\label{eq: G leq 0 charact}\begin{align}
    \phi^*(x_R)-\phi^*(x_L) \leq \phi(x_R)-\phi(x_L) &\quad \quad (\text{if }\phi \text{ is monotonically increasing}),\label{eq: monot incr}\\
    \phi^*(x_L)-\phi^*(x_R) \leq \phi(x_L)-\phi(x_R) &\quad \quad (\text{if }\phi \text{ is monotonically decreasing}).\label{eq: monot decr}
\end{align}
\end{subequations}
The equality in \eqref{eq: G leq 0 charact} holds when \eqref{eq: G leq 0} holds with equality.
\end{lemma}
\begin{proof}
We omit the proof of the general case and consider instead two simple cases. Without loss of generality, assume that $\phi$ is monotonically increasing.  Suppose first that ${\rm D}\phi^*$ has no sign changes, i.e. $\phi^*$ is monotonically increasing on $\omega=[x_L,x_R]$. Then, \eqref{eq: form G N even Dphih geq 0} reduces to:
  \begin{align}\label{eq: monot proof}
      \mathscr{G}_{\phi,\omega}(\phi^*) = (\phi^*(x_R)-\phi^*(x_L)) - (\phi(x_R)-\phi(x_L)),
  \end{align}
  which is negative if and only if \eqref{eq: monot incr} holds. Now suppose that ${\rm D}\phi^*>0$ has a single change of sign, say at $x_1$, and that ${\rm D}\phi^*>0$ for $x<x_1$. Additionally, since we may shift $\phi^*$ by $\phi'(x_L)$, we take $\phi'(x_L)=0$. It is easy to verify that $\phi^*(x_1)\leq \phi(x_R)$, as otherwise $\mathscr{G}_{\phi,\omega}(\phi^*)\geq 0$. Since $\phi^*$ is decreasing on $[x_1,x_R]$ and $\phi$ is increasing, we immediately have $\phi'(x_R) \geq 0$. The case ${\rm D}\phi^*<0$ for $x<x_1$ follows from a similar argument.
\end{proof}

Next, we introduce the classical definition of the maximum principle.
\begin{definition}[Maximum principle]
An approximation $\phi^* : \Omega \rightarrow \mathbb{R}$ of $\phi  : \Omega \rightarrow \mathbb{R}$ satisfies the maximum principle on $\omega \subset \Omega$ if and only if it does not exceed the bounds of $\phi$ on $\omega$:
\begin{subequations}
\begin{alignat}{2}
  \inf_\omega \phi^* &~ \geq  \inf_\omega \phi, \label{eq: monot approx inf}\\
  \sup_\omega \phi^* &~ \leq  \sup_\omega \phi. \label{eq: monot approx sup}
\end{alignat}
\end{subequations}
\end{definition}

\begin{theorem}[Interpolatory approximation of monotonic function]
  Suppose that $\phi \in H^1(\Omega)$ is monotonic on $\omega = [x_L,x_R]\subset \Omega$ and $\phi^* \in H^1(\tilde{\Omega})$ is an interpolatory approximation, i.e. $\phi^*(x_L)=\phi(x_L)$ and $\phi^*(x_R)=\phi(x_R)$. We have the following results:
  \begin{enumerate}
      \item $\mathscr{G}_{\phi,\omega}(\phi^*) \geq 0$,
  \item $\phi^*$ is a monotonic function if and only if 
  $\mathscr{G}_{\phi,\omega}(\phi^*) = 0$,
  \item if $\mathscr{G}_{\phi,\omega}(\phi^*) = 0$
then $\phi^*$ satisfies the maximum principle on $\omega$.
  \end{enumerate}
\end{theorem}
\begin{proof}
1. Without loss of generality, assume $\phi$ is increasing. In case ${\rm D}\phi^*$ does not change sign, we invoke \eqref{eq: monot proof} of \cref{lem: Gibbs ID local} and obtain $\mathscr{G}_{\phi,\omega}(\phi^*) = 0$. In the other case, denote the locations of sign changes of ${\rm D}\phi^*$ as $x_i \in \omega, i = 1, \dots, N$ with $x_i\leq x_{i+1}$. 
\cref{lem: Gibbs Cont} provides:
\begin{itemize}
    \item for ${\rm D}\phi^*\geq 0$ on $[x_L,x_1]$:
    \begin{align}\label{eq: form G N odd Dphih geq 0 2}
      \mathscr{G}_{\phi,\omega}(\phi^*)= 2 \sum_{i=1}^N (-1)^i \phi'(x_i),
    \end{align}
    \item for ${\rm D}\phi^*\leq 0$ on $[x_L,x_1]$:
    \begin{align}\label{eq: form G N odd Dphih leq 0 2}
      \mathscr{G}_{\phi,\omega}(\phi^*)=- 2 \sum_{i=1}^N (-1)^i \phi'(x_i) ,
    \end{align}
\end{itemize}
where we recall $\phi'= \phi - \phi^*$.
To show non-negativity of \eqref{eq: form G N odd Dphih geq 0 2} and \eqref{eq: form G N odd Dphih leq 0 2}, one has to consider various cases. For example, for $N=1$ we have $\mathscr{G} = 2 \left|\phi'(x_1)\right|\geq 0$. We omit a detailed proof of the general case.\\

2. `$\Rightarrow$': Suppose $\phi^*$ is a monotonic function. Without loss of generality, assume $\phi$ is increasing. Since $\phi^*$ is a monotonic function, the sums in \eqref{eq: form G N odd Dphih leq 0 2} and \eqref{eq: form G N odd Dphih geq 0 2} are empty and the expressions vanish. 

`$\Leftarrow$': Suppose that $\mathscr{G}_{\phi,\omega}(\phi^*)=0$. \cref{prop: lower bound G} implies $\mathscr{V}_{\phi,\omega}(\phi^*)\leq 0$. Without loss of generality, suppose that $\phi$ is monotonically increasing. We arrive at:
\begin{align}
    V_\omega(\phi^*) \leq \phi(x_R)-\phi(x_L). 
\end{align}
Since $\phi^*$ is interpolatory, \cref{lem: total variation bounds} implies that $\phi^*$ is monotonically increasing.

3. If $\mathscr{G}_{\phi,\omega}(\phi^*)=0$ then $\phi^*$ is monotonic via claim 2, and \eqref{eq: monot incr} provides the maximum principle via:
  \begin{subequations}
  \begin{alignat}{2}
      &\sup_\omega \phi^* = \phi^*(x_R) = \phi(x_R) = \sup_\omega \phi,\\
      &\inf_\omega \phi^* = \phi^*(x_L) = \phi(x_L) = \inf_\omega \phi.
  \end{alignat}
  \end{subequations}

\end{proof}

We have the following connection between the Gibbs constraints and the maximum principle.
\begin{lemma}[Gibbs constraints and maximum principle]
Let the analytical profile $\phi^* \in H^1(\Omega)$ be monotonic on each $\omega_j, j=1, \dots, J$, and let $\phi^* \in H^1(\tilde{\Omega})$ be an approximation that satisfies $\mathscr{G}_{\phi,\omega_j}(\phi^*)\leq 0, j = 1, \dots, J$, such that $\phi^*$ is interpolatory on $\partial\Omega$. Then $\phi^*$ satisfies the maximum principle on $\Omega$.
\end{lemma}

In order to preclude the Gibbs phenomenon on the entire domain, the strategy is to require $\mathscr{G}_{\phi,\omega_j}(\phi^*) \leq 0$, for $j=1, \dots, J$. We note that the practical applicability of this strategy relies on an appropriately chosen subdivision of $\Omega$. In the subsequent subsection, we discuss the Gibbs constraints in the context of best approximation problems. We return to the domain subdivision problem in the context of finite elements in \cref{subsec: Gibbs 1D fem}.

\subsection{Best approximation problems under Gibbs constraints}\label{subsec: Gibbs 1D: opt}

Let now $\phi \in H^1(\Omega)$ be given and consider the best approximation problem:
\begin{subequations}\label{eq: optimization problem}
\begin{align}
    \phi^* =&~ \underset{\theta^* \in \mathcal{K}}{\rm arginf}~ \|\phi-\theta^*\|_{\HH},
\end{align}
where the feasible set is defined as:
\begin{align}\label{eq: constraint set relax}
  \mathcal{K} := \left\{ \phi^* \in H^1(\tilde{\Omega}) :~\right.&\left. \mathscr{G}_{\phi,\omega_j}(\phi^*) \leq 0,~j = 1,...,J\right\}.
\end{align}
\end{subequations}
We note that, in general $\mathcal{K}$, has no strictly feasible point. For example, in the case $\phi \equiv 0$ in $\Omega$, we have $\mathscr{G}_{\phi,\omega_j}(\phi^*) = 0$ for $j=1,\dots,J$. In the following we exclude this trivial case.

The standard techniques to study best approximation problems are the gradient-based methods. We remark, however, that the Gateaux derivative of $\phi^* \rightarrow \mathscr{G}_{\phi,\omega_j}(\phi^*)$ does not exist due to the occurrence of the sign function. To permit the adoption of standard gradient methods, we regularize the nondifferentiable constraint function as follows. We introduce the parameter $\varepsilon \in \mathbb{R}$ and define the differentiable regularizations $\left|\cdot\right|_\varepsilon: \mathbb{R}\rightarrow \mathbb{R}_+$ and ${\rm sgn}_\varepsilon: \mathbb{R}\rightarrow \mathbb{R}$ as:
\begin{subequations}
  \begin{alignat}{2}
    \left|r\right|_\varepsilon :=& \left|r^2 + \varepsilon^2\right|^{1/2},\\
    {\rm sgn}_\varepsilon(r) :=& r \left|r\right|_\varepsilon^{-1},
  \end{alignat}
\end{subequations}
which satisfy the relation:
\begin{align}
    {\rm D}\left|r\right|_\varepsilon = {\rm sgn}_\varepsilon(r).
\end{align}
Next, we introduce the regularized functional
\begin{align}\label{eq: large g reg}
    \mathscr{G}^\varepsilon_{\phi,\omega}(\phi^*) := \displaystyle\int_{\omega} g^\varepsilon_{\phi}(\phi^*)~{\rm d}x,
\end{align}
where the functional $g^\varepsilon_{\phi}$ is defined as:
\begin{align}\label{eq: small g reg}
    g^\varepsilon_{\phi}(\phi^*) := \left| {\rm D}\phi^* \right|_\varepsilon- {\rm sgn}_\varepsilon\left({\rm D}\phi^*\right){\rm D}\phi= - {\rm sgn}_\varepsilon\left({\rm D}\phi^*\right){\rm D}\phi' + \varepsilon^2 \left|{\rm D}\phi^*\right|_\varepsilon^{-1}.
\end{align}
It is now clear that the regularized functional $\mathscr{G}^\varepsilon_{\phi,\omega}(\phi^*)$ is differentiable in the Gateaux sense. The Gateaux derivative of $\mathscr{G}_{\phi,\omega}^\varepsilon(\phi^*)$, denoted ${\rm d}\mathscr{G}_{\phi,\omega}^\varepsilon(\phi^*)(w)$, is given by:
\begin{align}\label{eq: Gat der eps}
    {\rm d}\mathscr{G}^\varepsilon_{\phi,\omega}(\phi^*)(w)= \displaystyle\int_{\omega_j} & \left|{\rm D}\phi^*\right|_\varepsilon^{-1}{\rm D}\phi^* {\rm D}w - \varepsilon^2 \left|{\rm D}\phi^*\right|_\varepsilon^{-3}{\rm D}\phi{\rm D}w ~{\rm d}x.
\end{align}

We have the following property regarding the convexity of the regularized Gibbs functional.

\begin{proposition}[Convexity of Gibbs functional]\label{prop: conv Gibbs func} The functional $g^\varepsilon_{\phi}$ is quasi-convex:
  \begin{align}
      g^\varepsilon_{\phi,\omega}(\zeta\phi_1^* + (1-\zeta)\phi_2^*) \leq \max\left\{ g^\varepsilon_{\phi,\omega}(\phi_1^*),g^\varepsilon(\phi_2^*) \right\},
  \end{align}
  for all $\phi_1^*, \phi_2^* \in H^1(\tilde{\Omega})$, $\zeta \in [0,1]$. The functional $ \mathscr{V}_{\phi,\omega}$ is convex, but $\mathscr{G}^\varepsilon_{\phi,\omega}$ is in general not (quasi-)convex.
\end{proposition}

We now consider the best approximation problem:
\begin{subequations}\label{eq: optimization problem reg}
\begin{align}
    \phi^* =&~ \underset{\theta^* \in \mathcal{K}^\varepsilon}{\rm arginf}~ \|\phi-\theta^*\|_{\HH},
\end{align}
where the regularized feasible set is defined as:
\begin{align}\label{eq: constraint set relax reg}
  \mathcal{K}^\varepsilon := \left\{ \phi^* \in H^1(\tilde{\Omega}) :~\right.&\left. \mathscr{G}^\varepsilon_{\phi,\omega_j}(\phi^*) \leq 0,~j = 1,...,J \right\}.
\end{align}
\end{subequations}

A consequence of \cref{prop: conv Gibbs func} is that the feasible set $\mathcal{K}^\varepsilon$ is not convex. This is a result of the occurrence of the (regularized) sign function in $g_{\phi}^\varepsilon$.

\begin{remark}
  The feasible set determined by (a regularized) functional\\ $\mathscr{V}_{\phi,\omega}(\phi^*) \leq 0$, i.e.
  \begin{align}\label{eq: constraint set V}
  \left\{ \phi^* \in H^1(\tilde{\Omega}) :~\right.&\left. \mathscr{V}^\varepsilon_{\phi,\omega_j}(\phi^*) \leq 0, j = 1,...,J \right\},
\end{align}
with
\begin{align}\label{eq: constraint set def V eps}
  \mathscr{V}^\varepsilon_{\phi,\omega_j}(\phi^*)=\displaystyle\int_\omega |{\rm D}\phi^*|_\varepsilon - |{\rm D} \phi|_\varepsilon ~{\rm d}x,
\end{align}
is convex. Furthermore, note that quasi-convexity of the functional is sufficient for the associated feasible set to be convex.
\end{remark}

We now introduce the first-order optimality Karush-Kuhn-Tucker (KKT) conditions of the constrained best approximation problem \eqref{eq: optimization problem reg}. We note that the solution of the KKT conditions is only guaranteed to be a local optimum due to the lack of convexity of the feasible set $\mathcal{K}^\varepsilon$. Replacing $\mathcal{K}^\varepsilon$ by the feasible set from \eqref{eq: constraint set V} would imply global optimality.

\begin{theorem}[Karush-Kuhn-Tucker conditions]\label{thm: KKT 1} The function $\phi^* \in \mathcal{K}$ is a local optimum of the problem \eqref{eq: optimization problem reg} if and only if there exist Lagrange multipliers $\lambda_j \in \mathbb{R}$ ($j=1,...,J$) such that the following Karush-Kuhn-Tucker (KKT) conditions hold:\\\\
\begin{subequations}
\indent \textit{Stationarity:}\\\\
\indent \textit{Find }$\phi^* \in H^1(\tilde{\Omega}), \lambda_j \in \mathbb{R}$ \text{for } $j=1,...,J$ \textit{such that}
  \begin{align}\label{eq: thm: stationarity}
     \left(\phi^*-\phi,w\right)_{\HH} &+ \sum_{j=1}^J \lambda_j {\rm d}\mathscr{G}^\varepsilon_{\phi,\omega_j}(\phi^*)(w) = 0 \quad \text{for all }w \in H^1(\tilde{\Omega}),
  \end{align}
\indent \textit{Primal feasibility:}
\begin{align}
   \mathscr{G}^\varepsilon_{\phi,\omega_j}(\phi^*)  \leq 0  \quad \text{ for } j=1,...,J,
\end{align}
\indent \textit{Dual feasibility:}
\begin{align}\label{eq: thm: dual feasibility}
    \lambda_j \geq 0 \quad \text{ for } j=1,...,J,
\end{align}
\indent \textit{Complementary slackness:}
\begin{align}
   \lambda_j \mathscr{G}^\varepsilon_{\phi,\omega_j}(\phi^*) = 0 \quad \text{ for } j=1,...,J.
\end{align}
\end{subequations}
\end{theorem}

\begin{proposition}[Homogeneity Gibbs functional]\label{prop: homogeneity G}
  The functional $\mathscr{G}^{\varepsilon}_{\phi,\omega}$ satisfies the property:
  \begin{align}
    \lim_{\varepsilon \rightarrow 0} {\rm d}\mathscr{G}^{\varepsilon}_{\phi,\omega}(\phi^*)(\phi^*) =     V_\omega(\phi^*) \geq 0.    
  \end{align}
\end{proposition}
\begin{proof}
Substitution of $w=\phi^* \in H^1(\tilde{\Omega})$ into \eqref{eq: Gat der eps} provides
\begin{subequations}\label{eq: homogen Gibss}
\begin{align}
    {\rm d}\mathscr{G}^\varepsilon_{\phi,\omega}(\phi^*)(\phi^*)=&~ \displaystyle\int_{\omega} {\rm d}g^\varepsilon_{\phi}(\phi^*)(\phi^*) ~{\rm d}x,\\
    {\rm d}g^\varepsilon_{\phi}(\phi^*)(\phi^*) =&~\left|{\rm D}\phi^*\right|_\varepsilon - \varepsilon^2 \left|{\rm D}\phi^*\right|_\varepsilon^{-3}\left(\left|{\rm D}\phi^*\right|_\varepsilon^2+{\rm D}\phi{\rm D}\phi^*\right).
\end{align}
\end{subequations}
The integrand ${\rm d}g^\varepsilon_{\phi}(\phi^*)(\phi^*)$ vanishes for ${\rm D}\phi^*=0$. On the other hand, for ${\rm D}\phi^* \neq 0$ we have
\begin{align}\label{eq: positivity Gateaux w=phi^h}
   \lim_{\varepsilon \rightarrow 0} {\rm d}g^\varepsilon_{\phi}(\phi^*)(\phi^*) =  |{\rm D}\phi^*| \geq 0.
\end{align}
\end{proof}

\begin{lemma}[Positivity]\label{lem: positive inner product}
The solution of the optimization problem \eqref{eq: optimization problem} has the following inner product form that is positive:
\begin{align}
    \left(\phi',\phi^*\right)_{\HH} \geq 0.
\end{align}
\end{lemma}
\begin{proof}
Substitution of $w = \phi^*$ into the stationarity condition \eqref{eq: thm: stationarity} provides:
  \begin{align}
     -\left(\phi',\phi^*\right)_{\HH} &+ \sum_{j=1}^J \lambda_j {\rm d}\mathscr{G}^\varepsilon_{\phi,\omega}(\phi^*)(\phi^*) = 0,
\end{align}
where we have employed the substitution $\phi'=\phi-\phi^*$. The result now follows from taking the limit $\varepsilon \rightarrow 0$, utilizing \cref{prop: homogeneity G} and invoking the dual feasibility property \eqref{eq: thm: dual feasibility}.
\end{proof}

\subsection{Best approximations in finite element spaces}\label{subsec: Gibbs 1D fem}

In this subsection, we seek for finite element approximations that satisfy the Gibbs constraints. 
Consider the best approximation problem:
\begin{subequations}\label{eq: optimization problem fem all}
\begin{align}\label{eq: optimization problem fem}
    \phi^h =&~ \underset{\theta^h \in \mathcal{K}_{p,\alpha}}{\rm arginf}~ \|\phi-\theta^h\|_{\HH},
\end{align}
where the feasible set $\mathcal{K}_{p,\alpha}$ for a finite element approximation space of polynomial degree $p$ and regularity $\alpha$ is defined as:
\begin{align}\label{eq: constraint set fem}
  \mathcal{K}_{p,\alpha} := \left\{ \phi^h \in \mathcal{V}^h_{D;p,\alpha} :~\right.&\left. \mathscr{G}_{\phi,\omega_j}(\phi^h) \leq 0, j = 1,...,J,\right\}.
\end{align}
\end{subequations}

To proceed, we regularize the non-differentiable constraint according to \eqref{eq: large g reg}-\eqref{eq: small g reg} and introduce the KKT conditions of the regularized problem.\\

\begin{subequations}\label{eq: KKT fem}
\indent \textit{Stationarity:}\\\\
\indent \textit{Find }$\phi^h \in \mathcal{V}^h_{D,p,\alpha}, \lambda_j \in \mathbb{R}$ \text{for } $j=1,...,J$ \textit{such that}
  \begin{align}\label{eq: thm: stationarity fem}
     \left(\phi^h-\phi,w^h\right)_{\HH} &+ \sum_{j=1}^J \lambda_j {\rm d}\mathscr{G}^\varepsilon_{\phi,\omega_j}(\phi^h)(w^h) = 0 \quad \text{for all }w^h \in \mathcal{V}^h_{0,p,\alpha},
  \end{align}
\indent \textit{Primal feasibility:}
\begin{align}
   \mathscr{G}^\varepsilon_{\phi,\omega_j}(\phi^h)  \leq 0  \quad \text{ for } j=1,...,J,
\end{align}
\indent \textit{Dual feasibility:}
\begin{align}\label{eq: thm: dual feasibility fem}
    \lambda_j \geq 0 \quad \text{ for } j=1,...,J,
\end{align}
\indent \textit{Complementary slackness:}
\begin{align}\label{eq: thm: complementary slackness fem}
   \lambda_j \mathscr{G}^\varepsilon_{\phi,\omega_j}(\phi^h) = 0 \quad \text{ for } j=1,...,J.
\end{align}
\end{subequations}

The problem \eqref{eq: KKT fem} takes the following algebraic form:\\
  \textit{Find }$\boldsymbol{\phi}^h, \boldsymbol{\lambda}$ such that: 
  \begin{subequations}\label{eq: prob 1}
  \begin{align}
      \mathbf{M} \boldsymbol{\phi}^h =&~ \mathbf{M} \boldsymbol{\phi} - \boldsymbol{\lambda}^T\mathbf{G}(\boldsymbol{\phi}^h),\\
      \mathbf{g}(\boldsymbol{\phi}^h) \leq &~ 0,\\
      \boldsymbol{\lambda} \geq&~ 0,\\
      \boldsymbol{\lambda}^T \mathbf{g}(\boldsymbol{\phi}^h) =&~0,
  \end{align}
\end{subequations}
  where the matrices $\mathbf{M} = [M_{AB}] \in \mathbb{R}^{n_{{\rm dof}}\times n_{{\rm dof}}}$ and $\mathbf{G}(\boldsymbol{\phi}^h) = [G_{jB}] \in \mathbb{R}^{n_{{\rm el}}\times n_{{\rm dof}}}$ are given by:
  \begin{subequations}
  \begin{align}
 M_{AB} =&~ (N_A, N_B)_{\mathcal{H}} ,\\
G_{jB} = &~ \displaystyle\int_{\omega_j} \left|{\rm D}\phi^h\right|_\varepsilon^{-1}{\rm D}\phi^h {\rm D}N_B - \varepsilon^2 \left|{\rm D}\phi^h\right|_\varepsilon^{-3}{\rm D}\phi{\rm D}N_B ~{\rm d}x,
  \end{align}
\end{subequations}
and the vectors are $\mathbf{g}(\boldsymbol{\phi}^h) = [\mathscr{G}_{\phi,\omega_j}(\phi^h)] \in \mathbb{R}^{n_{{\rm el}}}$ and $\boldsymbol{\lambda} = [\lambda_j]$. Here, $n_{{\rm dof}}$ denotes the number of degrees of freedom.

\begin{remark}[Computation constrained solutions]
In this article, we use carefully selected examples that allow the construction of solutions of the constrained best approximation problem \eqref{eq: optimization problem fem all}. The computation of constrained solutions with standard gradient-based methods is in general very difficult. There are a number of challenges. First, in some situations the number of degrees of freedom in the feasible solution set may be very small. It can even occur that the feasible set consists of a single function (see \cref{rmk: unique feasible sol}). This is extremely difficult to find with gradient-based methods. Second, the problem is in general non-convex. As a consequence, gradient-based methods may get stuck in local optima. This excludes a large class of powerful methodologies from convex optimization (that often rely on the KTT conditions). Third, the problem is highly nonlinear. This means that standard Newton-Raphson linearization methods often do not converge, and one has to work with less efficient approaches such as quasi-Newton type methods. Finally, we note that similar issues occur in the computation of $L^q$-best approximations when taking $q\rightarrow 1$, and especially when $q=1$.

\end{remark}

We start with discontinuous approximation spaces ($\alpha = -1$). In order to apply the Gibbs constraints in practice, the subdomains $\omega_j$ need to be selected. We select the subsets $\omega_j \subset \Omega, j=1,\dots,J$ as the finite elements: $\omega_j = K_j$. Note that the set $\mathcal{K}_{p,-1}$ is not empty (it contains at least all piecewise constants). The feasible set is convex for discontinuous piecewise linear basis functions (recall that it is general not convex). 

\begin{theorem}[Convexity of $\mathcal{K}_{1,-1}$]\label{thm: convexity Km11}
  The feasible set $\mathcal{K}_{1,-1}$ is convex.
\end{theorem}
\begin{proof}
  Let $\phi^h_1, \phi^h_2 \in \mathcal{K}_{1,-1}$ be given. Fix an element $K_i=(x_{i},x_{i+1})$; on this element the functions $\phi^h_1 $ and $\phi^h_2$ have a representation $\phi^h_1 = a_1 x + b_1$ and $\phi^h_2 = a_2 x + b_2$ for some scalars $a_1, a_2, b_1, b_2 \in \mathbb{R}$. Since $\phi^h_1, \phi^h_2 \in \mathcal{K}_{1,-1}$, we now have:
  \begin{subequations}\label{eq: GK phi1 phi2 d=1 p=1}
    \begin{align}
        \begin{split}\mathscr{G}_{\phi,K_i}(\phi^h_1)&~= \displaystyle\int_{x_{i}}^{x_{i+1}} \left|{\rm D} \phi_1^h\right|- {\rm sgn}\left({\rm D} \phi_1^h\right){\rm D}\phi ~{\rm d}x \\
        &~= (x_{i+1}-x_{i})|a_1| - {\rm sgn}\left(a_1\right)(\phi(x_{i+1})-\phi(x_i))\leq 0,\end{split}\\
        \begin{split}\mathscr{G}_{\phi,K_i}(\phi^h)&~= \displaystyle\int_{x_{i}}^{x_{i+1}} \left|{\rm D} \phi_2^h\right|- {\rm sgn}\left({\rm D} \phi_2^h\right){\rm D}\phi {\rm d}x \\
        &~= (x_{i+1}-x_{i})|a_2| - {\rm sgn}\left(a_2\right)(\phi(x_{i+1})-\phi(x_i)) \leq 0.\end{split}
    \end{align}
  \end{subequations}
Without loss of generality, we assume $\phi(x_{i+1})-\phi(x_i)>0$. From \eqref{eq: GK phi1 phi2 d=1 p=1}, we find $a_1>0$ and $a_2>0$, and we can write:
  \begin{subequations}\label{eq: GK phi1 phi2 d=1 p=1 2}
    \begin{align}
        (x_{i+1}-x_{i})a_1 - (\phi(x_{i+1})-\phi(x_i))&\leq 0,\\
        (x_{i+1}-x_{i})a_2 - (\phi(x_{i+1})-\phi(x_i)) &\leq 0.
    \end{align}
  \end{subequations}
Now, let $\zeta \in [0,1]$ and define $\phi^h_{\zeta} = \zeta \phi^h_1 + (1-\zeta) \phi_2^h$. Convexity is then a direct consequence of \eqref{eq: GK phi1 phi2 d=1 p=1 2}:
\begin{align}
        \mathscr{G}_{\phi,K_i}(\phi^h_\zeta)=&~ \displaystyle\int_{x_{i}}^{x_{i+1}} \left|{\rm D} \phi_\zeta^h\right|- {\rm sgn}\left({\rm D} \phi_\zeta^h\right){\rm D}\phi {\rm d}x\nn\\
        =&~ (x_{i+1}-x_{i})(\zeta a_1 + (1-\zeta) a_2) - (\phi(x_{i+1})-\phi(x_i))\nn\\
        \leq&~ 0.
    \end{align}
\end{proof}

\begin{figure}[!ht]
\begin{subfigure}{0.49\textwidth}
\centering
\includegraphics[width=0.95\textwidth]{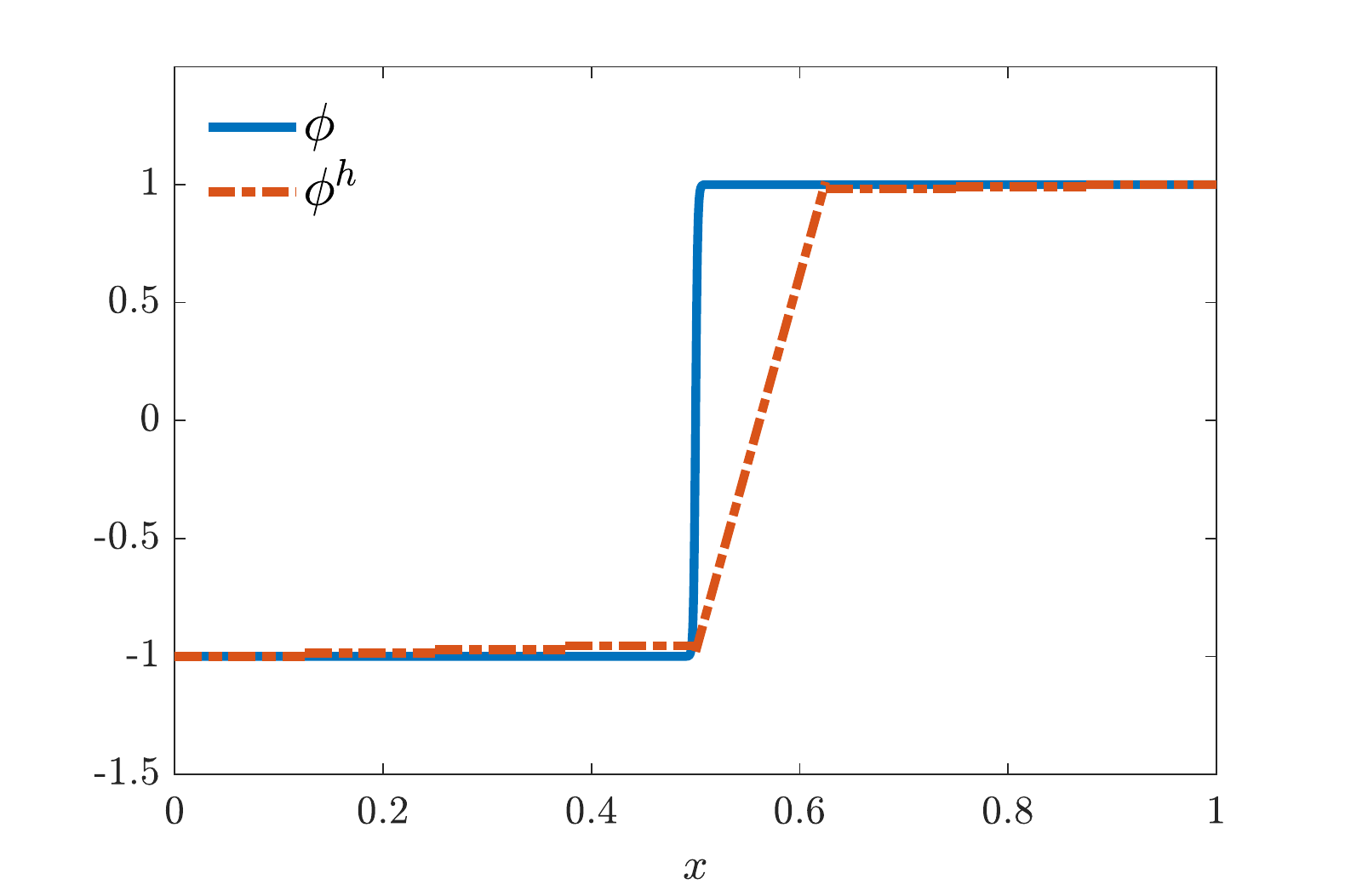}
\caption{$p=1$.}
\end{subfigure}
\begin{subfigure}{0.49\textwidth}
\centering
\includegraphics[width=0.95\textwidth]{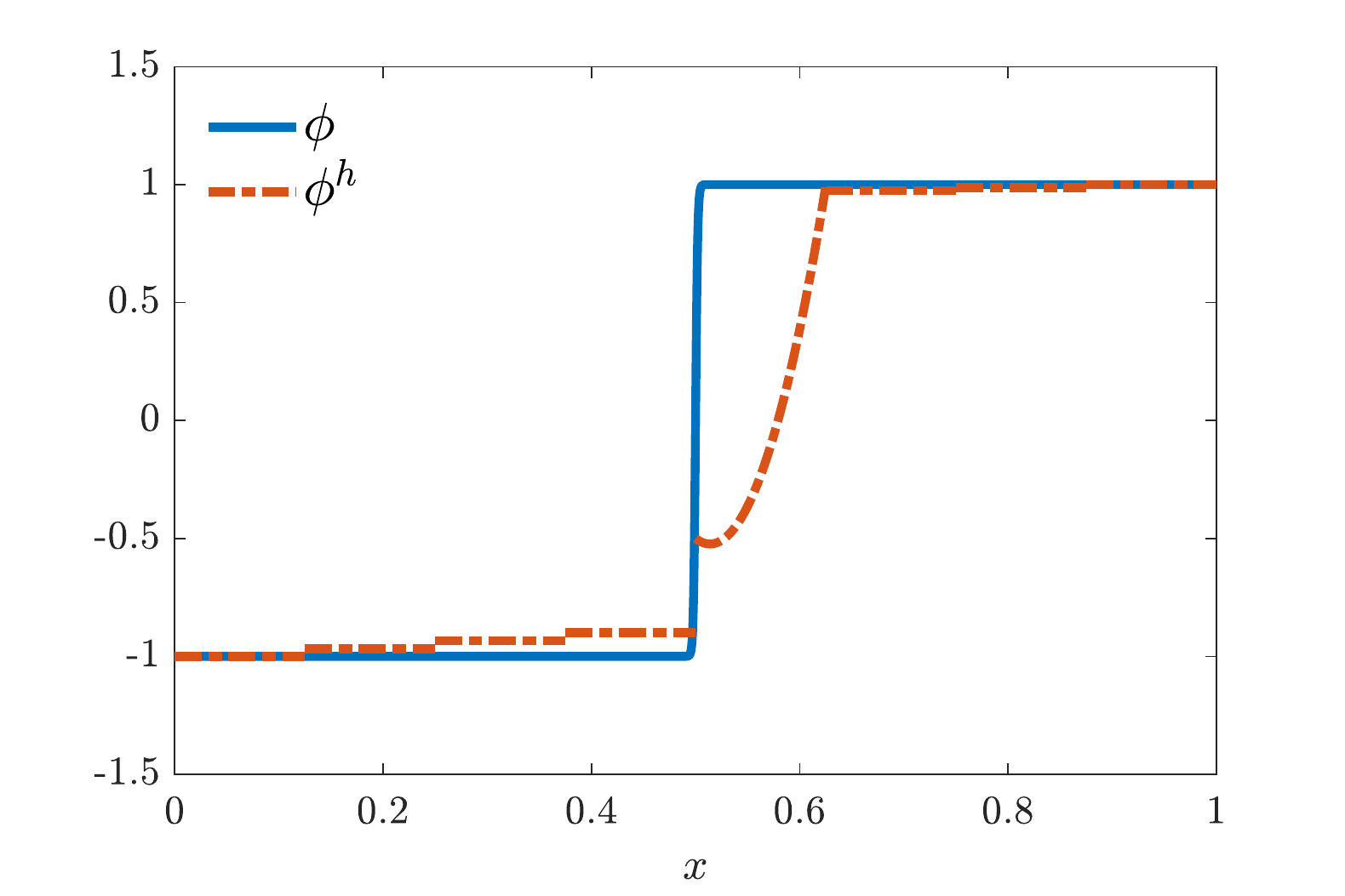}
\caption{$p=2$.}
\end{subfigure}
\caption{The interior penalty-best approximation $\phi^h \in \mathcal{V}_{D;p,-1}^h$, with $n_{\rm el} =8$, of the smooth step function $\phi = \phi_{0.5}$, subject to element-wise Gibbs constraints, for different for $p$.}
\label{fig: IP Gibbs free 1D}
\end{figure}

Consider now the best approximation problem \eqref{eq: optimization problem fem}-\eqref{eq: constraint set fem} with the interior penalty optimality \eqref{eq: proj IP}-\eqref{eq: best approx IP}, subject to the element-wise Gibbs constraints. Note that the Gateaux derivatives ${\rm d}\mathscr{G}^{\varepsilon}_{\phi,\omega_j}$ are linearly independent. Therefore, the Lagrange multiplier $\lambda_j$ solely depends on quantities defined on $K_j$. We visualize the interior penalty-best approximation $\phi^h \in \mathcal{V}_{D;p,-1}^h$ of the smooth step function, subject to the element-wise Gibbs constraints, in \cref{fig: IP Gibbs free 1D} for $p=1, 2$.

We observe that the constrained interior penalty-best approximations do not show over- or undershoots. Note that this property is not valid in general. Not all feasible solutions are monotonic, since the addition of a piecewise constant does not alter the Gibbs functional. Furthermore, both approximations deviate by a small piecewise constant away from the sharp layer. This behavior diminishes as the penalty parameter $\eta$ is increased.

Next, we focus on the case $\alpha = 0$, i.e. continuous approximation spaces. Again, we take $\omega_j=K_j$. It is easy to see that $\mathcal{K}_{p,0}$ is not empty. Namely, $\mathcal{K}_{1,0}$ contains the piecewise linear interpolant, and $\mathcal{K}_{1,0} \subset \mathcal{K}_{p,0}$ for $p \geq 1$. Considering the best approximation problem \eqref{eq: optimization problem fem}-\eqref{eq: constraint set fem} with $\mathcal{H}=H_0^1(\Omega)$, we have the following lemma. 

\begin{lemma}[$H_0^1$-orthogonality continuous linears]\label{lem: vanishing element IP for 1D P1}
  Let $\phi^h \in \mathcal{V}^h_{D,1,0}$ be the $H_0^1$-best approximation of $\phi \in \mathcal{V}$. We have the property:
  \begin{align}
   \int_{K_i} {\rm D}\phi^h {\rm D}\phi' ~{\rm d}x=0,
\end{align}
for all elements $K_i = (x_{L,i},x_{R,i}), i = 1,\dots, n_{{\rm el}}$.
\end{lemma}
\begin{proof}
Applying the first Green's identity we find:
\begin{align}
   \int_{K_i} {\rm D}\phi^h {\rm D}\phi'~{\rm d}x =\int_{\partial {K_i}}  n {\rm D} \phi^h \phi'~{\rm d}a - \int_{K_i} {\rm D}^2 \phi^h \phi'~{\rm d}x,
\end{align}
with $n$ the outward unit normal. Noting that for piecewise linear polynomials the second term vanishes and we are left with:
\begin{align}\label{eq: first-order 1D requirement}
  \int_{K_i} {\rm D}\phi^h {\rm D}\phi'~{\rm d}x = \phi'(x_{R,i}){\rm D} \phi^h(x_{R,i})-\phi'(x_{L,i}){\rm D} \phi^h(x_{L,i}).
\end{align}
Recalling now from \cref{lem: monotonicty H01} that the $H_0^1$-projector provides nodally exact solutions (i.e. $\phi'(x_{R,i})=\phi'(x_{L,i})=0$) completes the proof.
\end{proof}
\begin{lemma}[Vanishing Lagrange multipliers]\label{lem: lag KKT}
  Suppose that the approximation satisfies homogeneous boundary conditions, $\phi^h \in \mathcal{V}^h_{0,1,0}$. The Lagrange multiplier $\lambda_j$ in the KKT conditions \eqref{eq: KKT fem} vanishes if $\phi^h$ is not constant on $\omega_j$, for $\varepsilon \rightarrow 0$.
\end{lemma}
\begin{proof}
  Noting that the $H_0^1$-best approximation is monotonic and interpolatory, \cref{lem: interp mono} ensures satisfaction of the Gibbs constraints $\mathscr{G}_{\phi,K_i}(\phi^h) \leq 0$ (for $\varepsilon \rightarrow 0$). As a consequence, the first term in \eqref{eq: thm: stationarity fem} vanishes (consistent with \cref{lem: vanishing element IP for 1D P1}), and we are left with:
  \begin{align}
    \sum_{j=1}^J \lambda_j {\rm d}\mathscr{G}^\varepsilon_{\phi,K_j}(\phi^h)(w^h) = 0,
  \end{align}
  for all $w^h \in \mathcal{V}^h_{0,1,0}$. Substituting $w^h=\phi^h$, taking the limit $\varepsilon \rightarrow 0$ and invoking the homogeneity property of \cref{prop: homogeneity G} yields:  
  \begin{align}
    \sum_{j=1}^J \lambda_j V_{\omega_j}(\phi^h) = 0.
  \end{align}
  Noting that $\lambda_j \geq 0$ and $V_\omega(\phi^h)\geq 0$, we get $\lambda_j V_{\omega_j}(\phi^h)= 0$, for $j=1, \dots J$. If $\phi^h$ is not constant on $\omega_j$, we have $V_{\omega_j}(\phi^h)>0$ and thus $\lambda_j=0$.
\end{proof}

\begin{figure}[!ht]
\begin{subfigure}{0.49\textwidth}
\centering
\includegraphics[width=0.95\textwidth]{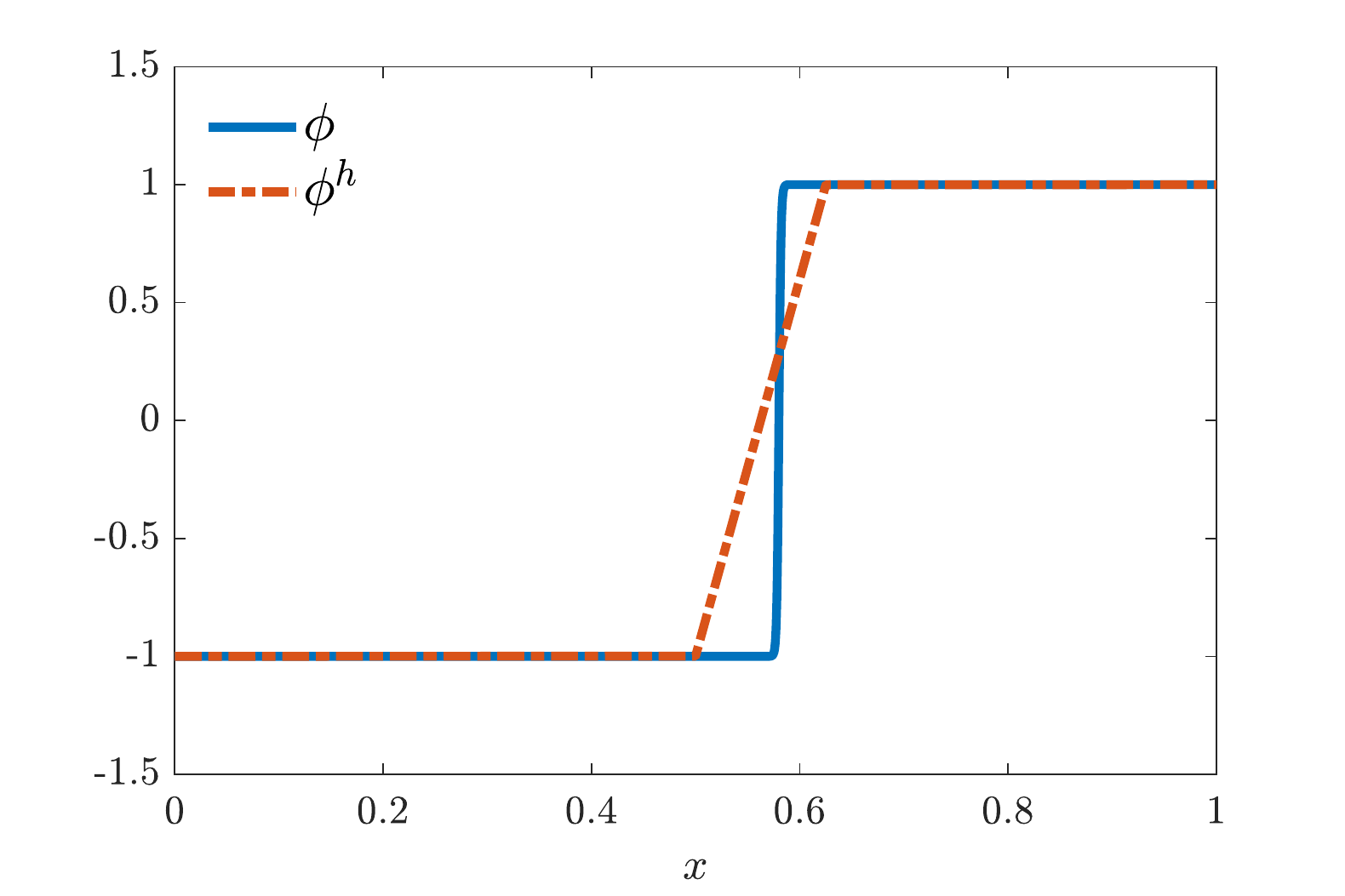}
\caption{$p=1$.}
\end{subfigure}
\begin{subfigure}{0.49\textwidth}
\centering
\includegraphics[width=0.95\textwidth]{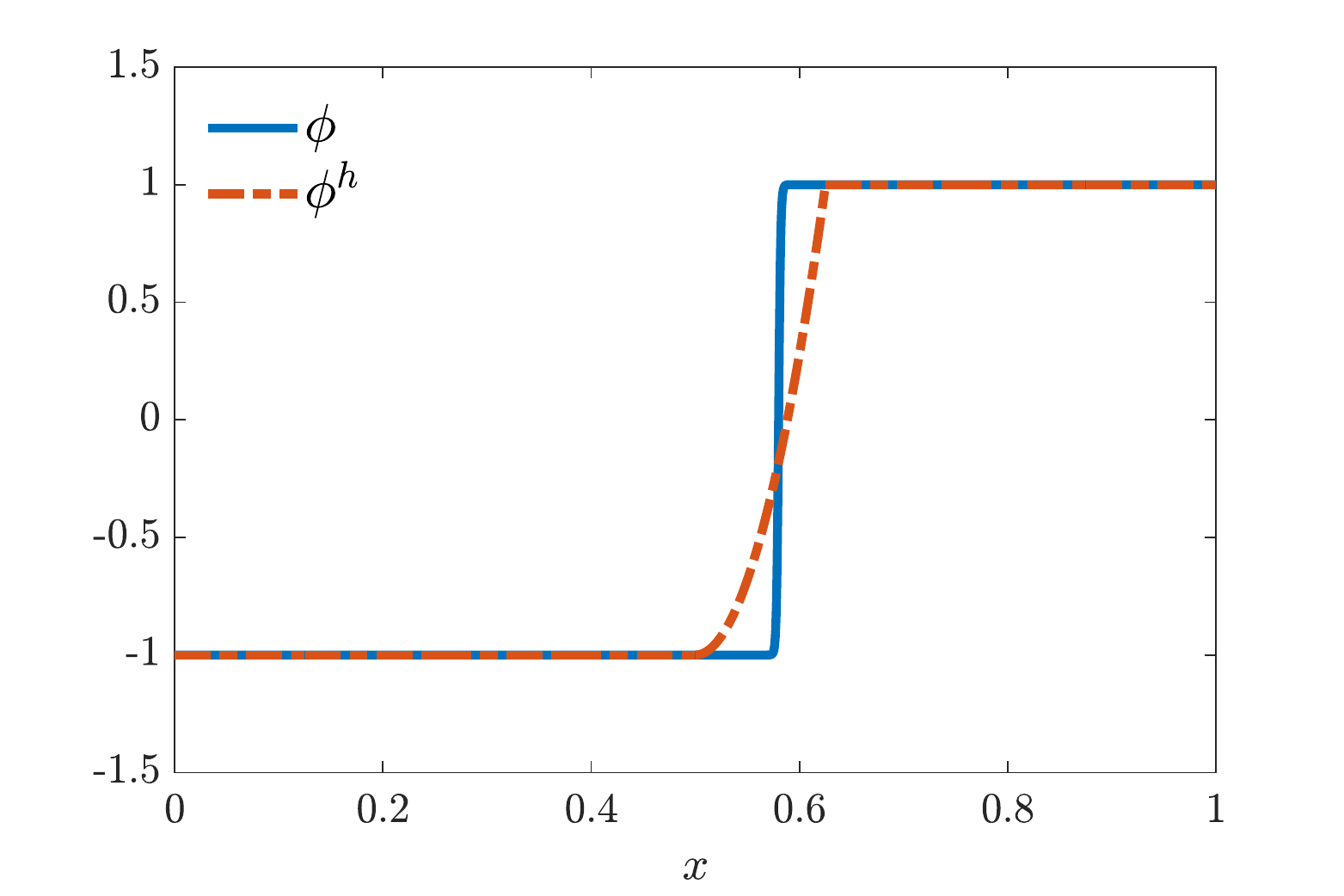}
\caption{$p=2$.}
\end{subfigure}
\caption{The $H_0^1$-best approximation $\phi^h \in \mathcal{V}_{D;p,0}^h$, with $n_{\rm el} =8$,  of the smooth step function $\phi = \phi_{0.58}$, subject to element-wise Gibbs constraints, for different $p$.}
\label{fig: H01 Gibbs free 1D}
\end{figure}

Note that, as a consequence of the overlapping support of the basis functions, the Gateaux derivatives ${\rm d}\mathscr{G}^\varepsilon_{\phi,K_j}$ are linearly dependent. Therefore, the Lagrange multipliers $\lambda_j, j=1,...,n_{el}$ are \textit{non-local}, in the sense that $\lambda_j$ does not solely depend on quantities defined on $\omega_j$. The $H_0^1$-best approximation $\phi^h \in \mathcal{V}_{D;p,0}^h$, $p=1, 2$, of the smooth step function (with $a=0.58$) subject to the element-wise Gibbs constraints is illustrated in \cref{fig: H01 Gibbs free 1D}.
We observe in both figures an interpolatory monotonic approximation. Moreover, we have $\mathscr{G}_{\phi,K_j}(\phi^h)=0$ for all element numbers $j=1,\dots, J$. In general, the following theorem holds.

\begin{theorem}[Monotonic interpolant for regularity $\alpha = 0$]\label{thm: monotonic interpolant}
  The constrained best approximation $\phi^h \in \mathcal{V}_{D,p,0}^h$ defined by the problem \eqref{eq: optimization problem fem}-\eqref{eq: constraint set fem} of a monotonic profile $\phi$ is a monotonic interpolant.
\end{theorem}

\begin{remark}[Uniqueness]\label{rmk: unique feasible sol}
  By \cref{thm: monotonic interpolant}, the best approximation for polynomial degree $p=1$ is the sole feasible solution in $\mathcal{K}_{1,0}$. This function is thus independent of the optimality condition $\mathcal{H}$. 
\end{remark}

We now turn our attention to higher-order smooth ($\alpha \geq 1$) approximation spaces $\mathcal{V}^h_{D,p,\alpha}$. The following proposition precludes existence of a solution for $\omega_j=K_j$ in general.

\begin{proposition}[Infeasible elementwise constraints]\label{prop: infeasible elementwise constraints}
  The constrained best approximation problem \eqref{eq: optimization problem fem}-\eqref{eq: constraint set fem} with $\omega_j=K_j$ and $a=0.58$ has in general no solution for $\alpha \geq 1$.
\end{proposition}

\begin{figure}
\centering
    \includegraphics[width=0.4655\textwidth]{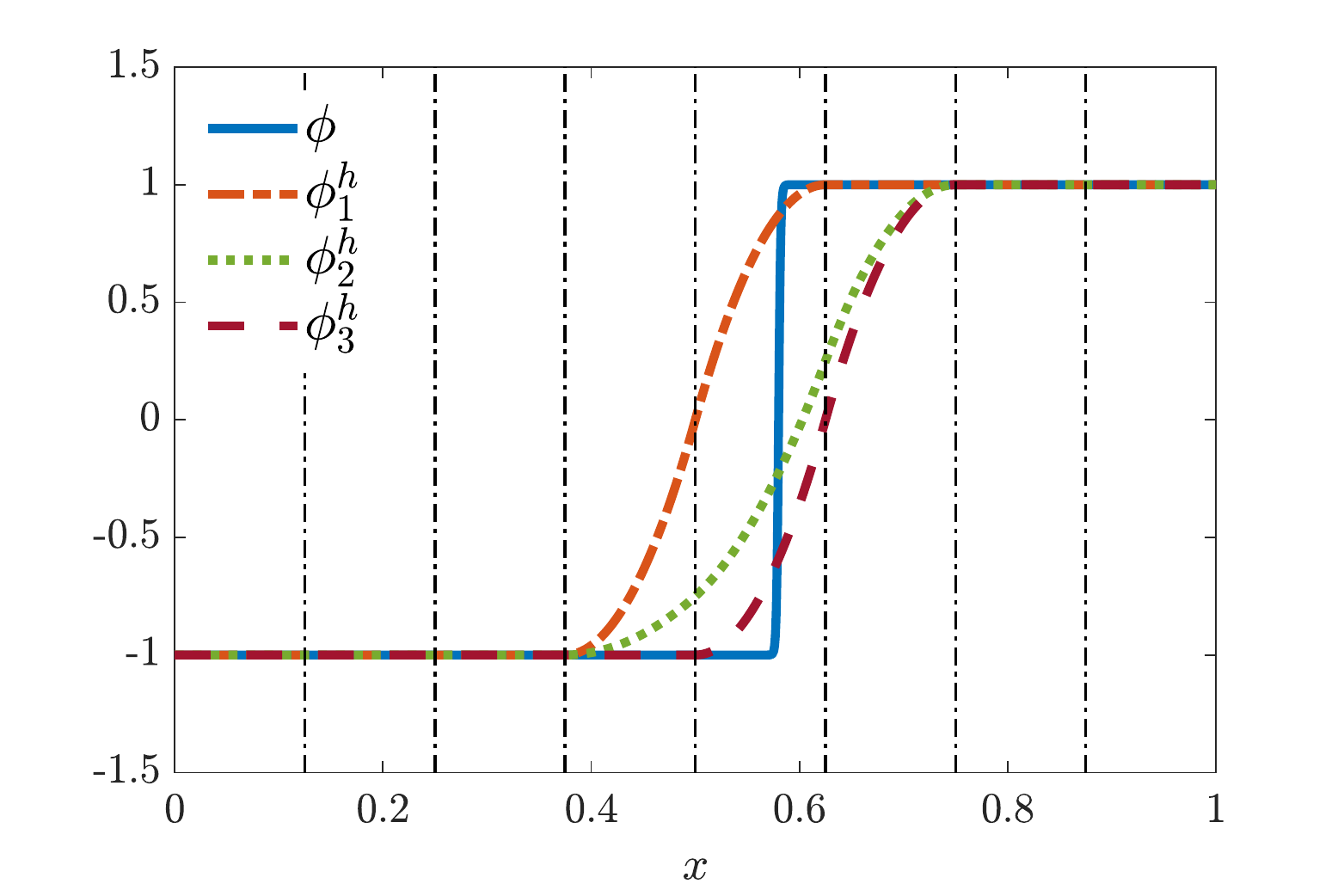}
    \caption{Sharp finite element approximations $\phi_1^h,\phi_2^h,\phi_3^h \in \mathcal{V}^h_{D,2,1}$, with $n_{\rm el} =8$, of the smooth one-dimensional step function $\phi = \phi_{0.58}$. The vertical dashed lines indicate the element boundaries.}
    \label{fig: sharp 1D}
\end{figure}

This is a consequence of the observation that interpolatory solutions exist for $\alpha = -1$ and $\alpha = 0$, but not for $\alpha \geq 1$. We illustrate this for the smooth step function \eqref{eq: approx step} with parameter $a=0.58$ in \cref{fig: sharp 1D}. Here, we display several sharp finite element approximations with quadratic basis functions ($p=2$) and regularity $\alpha = 1$. 
The finite element approximation requires at least $2$ elements to capture the sharp layer. There thus necessarily exists an element $K_i$ for which $\mathscr{G}_{\phi,K_i}(\phi^h)>0$, and the constrained best approximation problem \eqref{eq: optimization problem fem}-\eqref{eq: constraint set fem} therefore has no solution. 

\begin{figure}[!ht]
\begin{subfigure}{0.49\textwidth}
\centering
\includegraphics[width=0.95\textwidth]{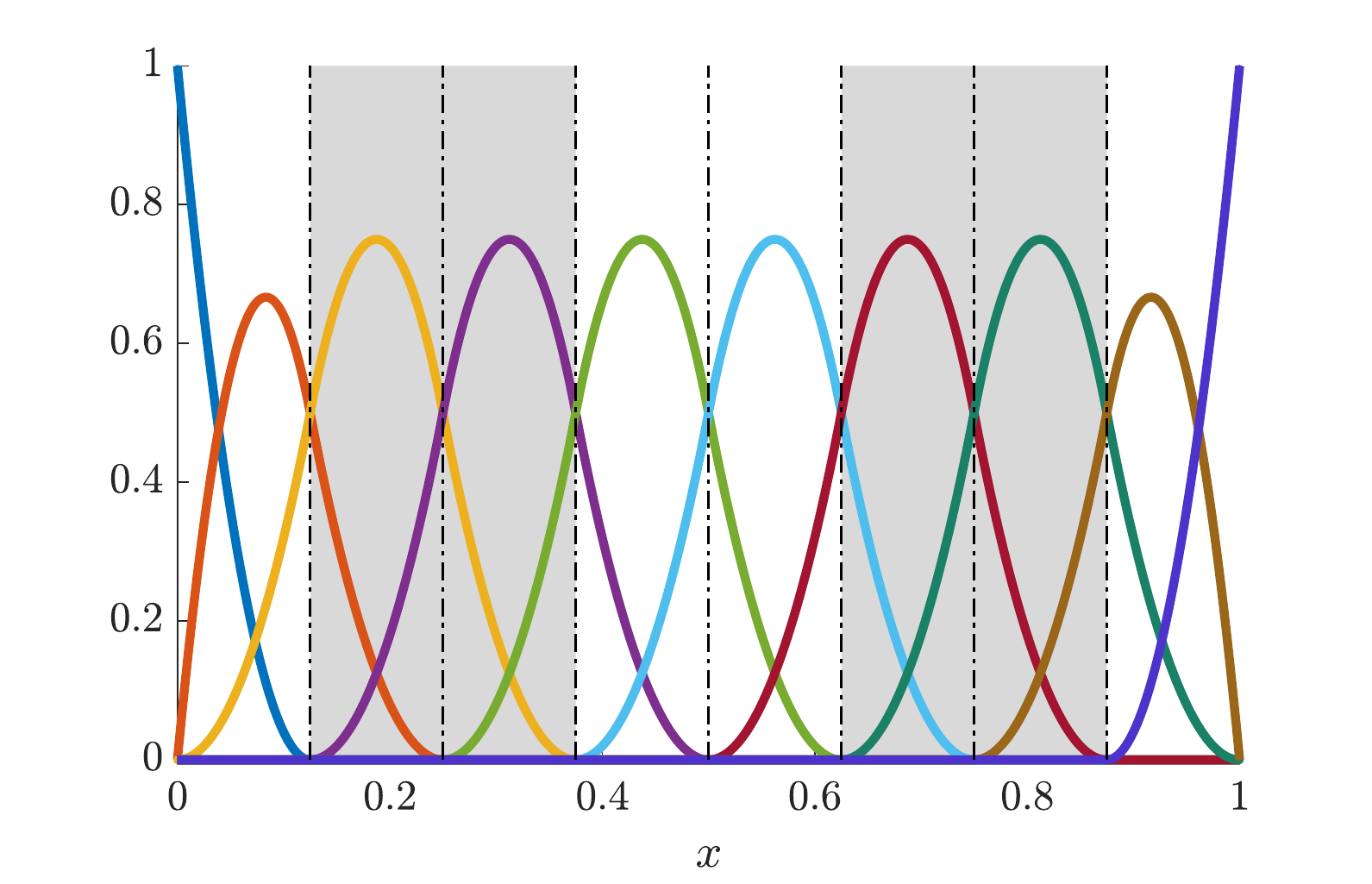}
\caption{Subdivision 1.}
\end{subfigure}
\begin{subfigure}{0.49\textwidth}
\centering
\includegraphics[width=0.95\textwidth]{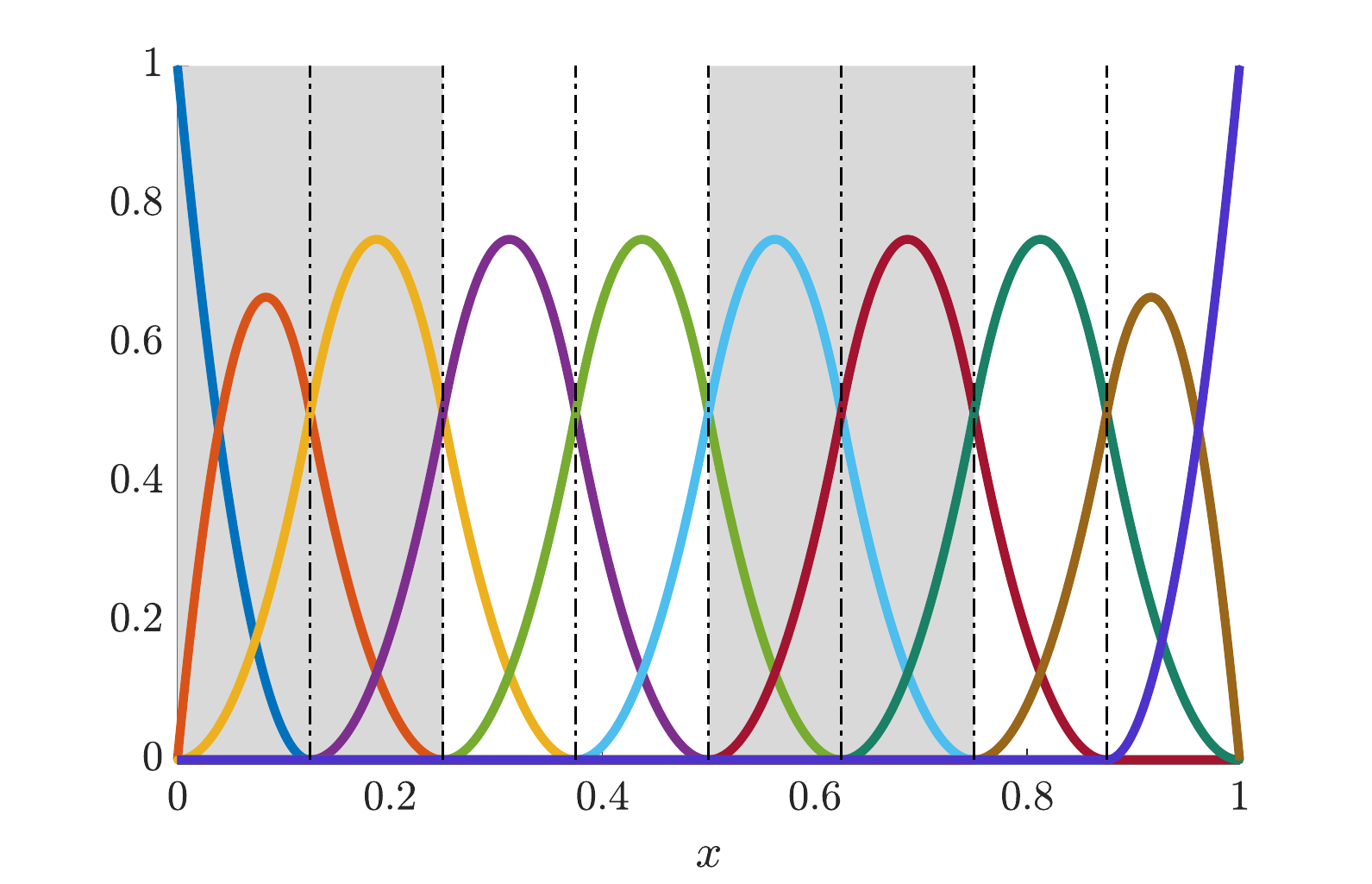}
\caption{Subdivision 2.}
\end{subfigure}
\caption{The two possible subdivisions for $\mathcal{V}^h_{2,1}$ when $n_{{\rm el}}=8$. The vertical dashed lines represent the element boundaries.}
\label{fig: subdivisions 1D}
\end{figure}

In general, B-spline basis functions of degree $p$ have the following local support:
\begin{align}
    {\rm supp}(N_{i,p}) = (\xi_i,\xi_{i+p+1}),
\end{align}
which depends on the regularity $\alpha$. The B-spline basis function $N_{i,p}$ has support in at most $p-k+2 = \alpha + 2$ elements. Neighboring B-spline basis functions share the support:
\begin{align}
    {\rm supp}(N_{i,p})\cap{\rm supp}(N_{i+1,p})  = (\xi_{i+1},\xi_{i+p+1}),
\end{align}
which consists of at most $\alpha + 1$ elements. With the aim of finding a subdivision of $\Omega$ that permits the sharpest approximations, we select $\omega_j$ as the union of $\alpha+1$ neighboring elements. We subdivide domain $\Omega$ into disjoint subdomains $\omega_j = \cup_{i \in \mathcal{I}_j} K_i$, where $\mathcal{I}_j$ denotes an index set. 
Note that shifting groups of elements yields a different subdivision (for $\alpha > 0$). The construction is therefore not unique and we consider $\alpha +1$ possibilities. For the first index set $\mathcal{I}_1$ we have the options $\left\{1\right\} , \dots, \left\{1, \dots, \alpha+1\right\}$. The consecutive index sets contain the next $\alpha +1$ consecutive numbers, where the last set terminates with the last element number. We visualize two possible subdivisions in \cref{fig: subdivisions 1D}.

\begin{figure}[!ht]
\begin{subfigure}{0.49\textwidth}
\centering
\includegraphics[width=0.95\textwidth]{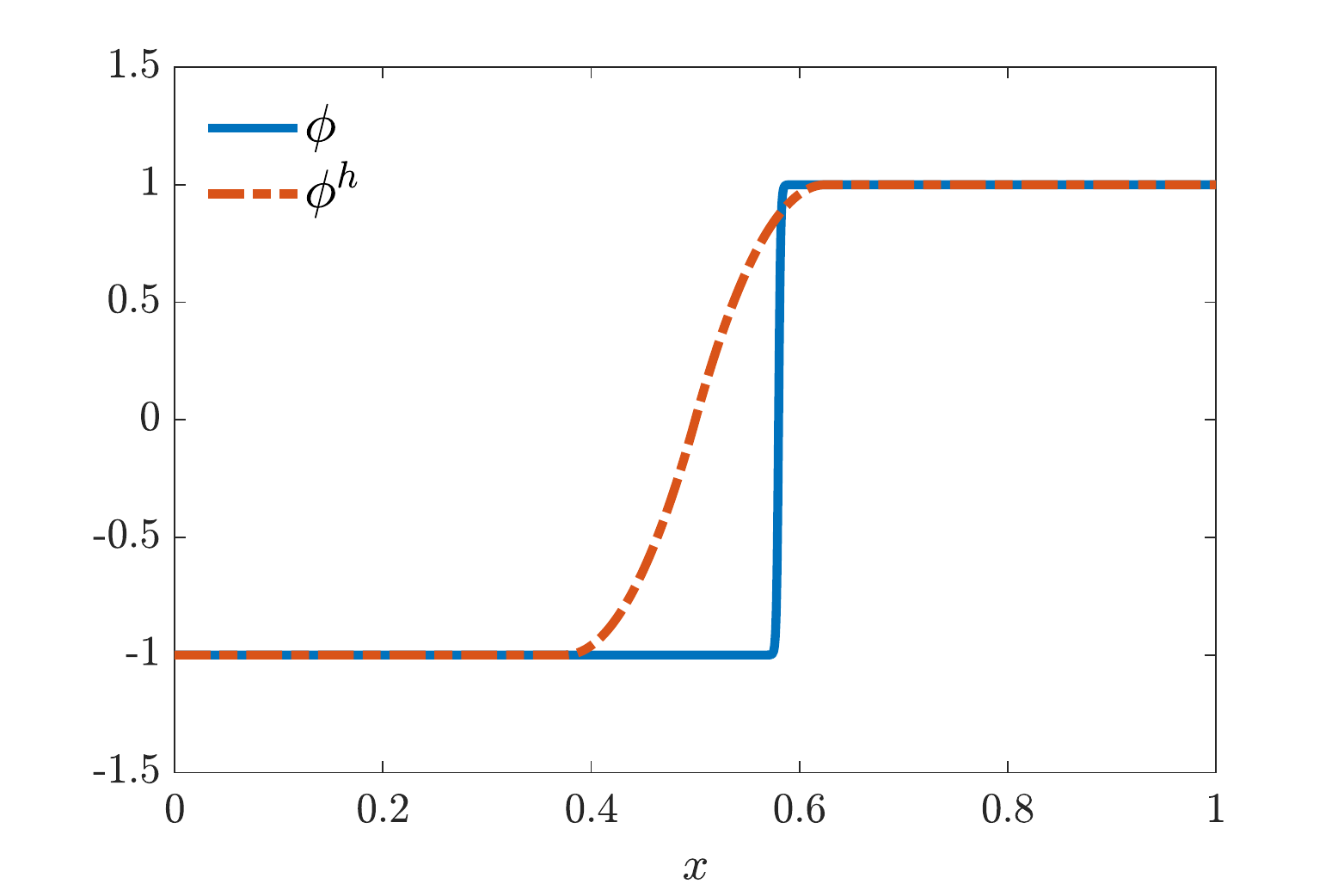}
\caption{Subdivision 1.}
\end{subfigure}
\begin{subfigure}{0.49\textwidth}
\centering
\includegraphics[width=0.95\textwidth]{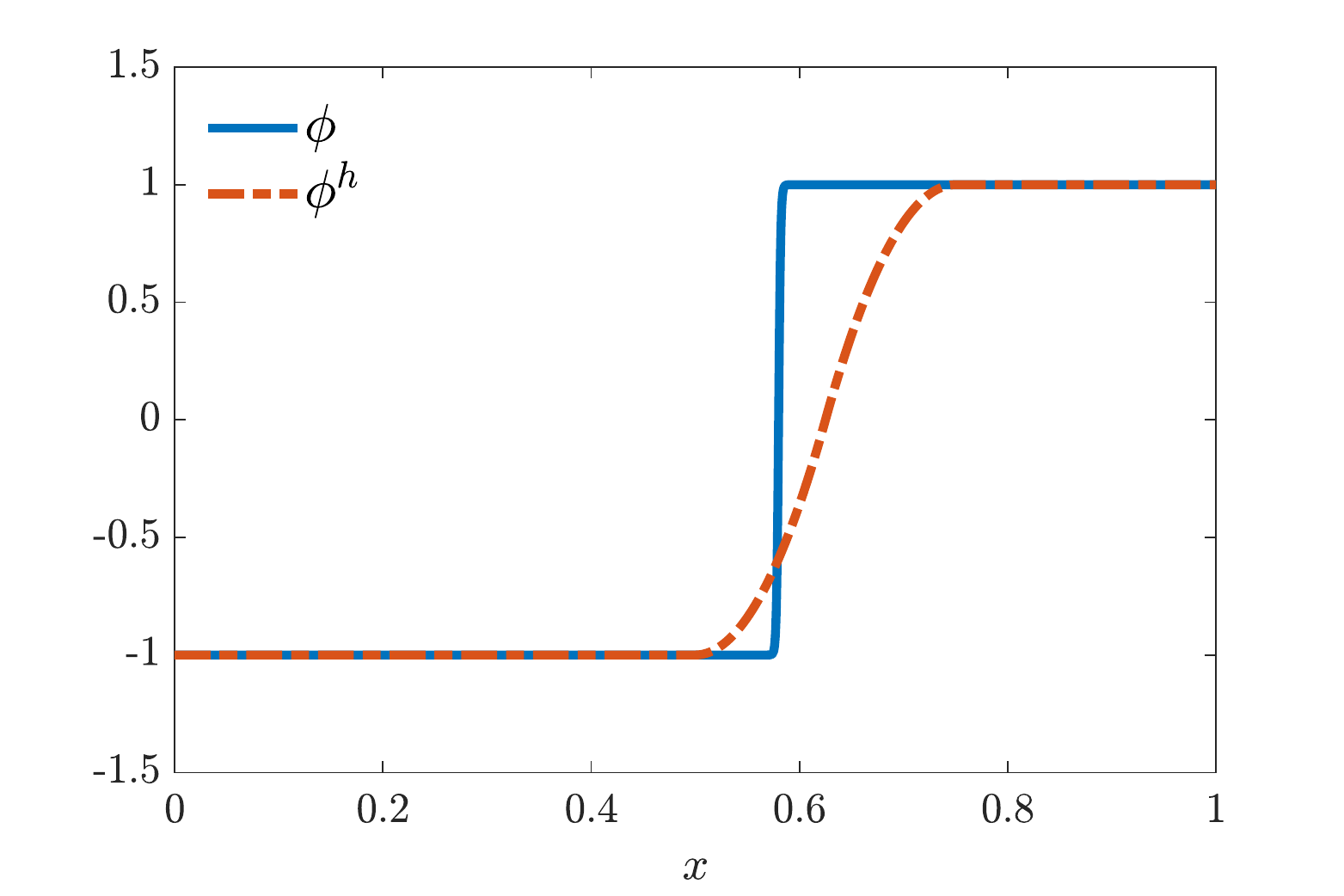}
\caption{Subdivision 2.}
\end{subfigure}
\caption{The $H_0^1$-best approximation $\phi^h \in \mathcal{V}_{D;2,1}^h$, with $n_{\rm el} =8$, of the smooth step function $\phi = \phi_{0.58}$, subject to Gibbs constraints on each of the subdomains illustrated in \cref{fig: subdivisions 1D}.}
\label{fig: H01 Gibbs free 1D p2m1}
\end{figure}

In \cref{fig: H01 Gibbs free 1D p2m1}, we visualize the $H_0^1$-best approximation $\phi^h \in \mathcal{V}_{D;2,1}^h$ of the smooth step function subject to Gibbs constraints on the subdomains of both subdivisions. We see that both approximations are completely free of over- and undershoots. In general, the existence of a feasible solution of the constrained best approximation problem with $\alpha \geq 1$ is an open question. 

\section{Eliminating the Gibbs phenomenon in higher dimensions}\label{sec: Gibbs higher D}

In this section, we extend the strategy to eliminate the Gibbs phenomenon presented in \cref{sec: Gibbs constraints 1D} to higher dimensions. We first introduce the Gibbs constraints for the elimination of the Gibbs phenomenon in general function spaces in \cref{subsec: Gibbs higher D: Gibbs constraints}. After that, in \cref{subsec: Gibbs higher D: Gibbs fem}, we discuss the Gibbs constraints in the context of best approximations in finite element spaces.

\subsection{Gibbs constraints}\label{subsec: Gibbs higher D: Gibbs constraints}

We present the construction of Gibbs constraints on a multidimensional given subdomain $\omega$, by building upon the one-dimensional framework presented in \cref{sec: Gibbs constraints 1D}. The main distinguishing feature in the design of Gibbs constraints in higher dimensions is the notion of directionality. As a consequence, the occurrence of over- and undershoots now depends on the point of view, meaning that there exists no precise mathematical definition of monotonicity in higher dimensions. To construct a set of constraints
, we instead aim to preclude the Gibbs phenomenon in a particular direction. We start by assuming an arbitrary direction, $\mathbf{e}$, and afterwards suggest a particular closed form expression for $\mathbf{e}$.

\begin{definition}[Directional Gibbs functional]\label{def: Gibbs nD}
The directional Gibbs functional of the function $\phi^* \in H^1(\tilde{\Omega})$, with respect to the given function $\phi \in H^1(\Omega)$ and on $\omega \subset \Omega$, is defined as:
\begin{align}\label{eq: large g nD dir}
    \mathscr{G}_{\phi,\omega}(\phi^*;\mathbf{e}) := \displaystyle\int_{\omega} g_{\phi}(\phi^*,\mathbf{e})~{\rm d}\Omega,
\end{align}
where the functional $g_\phi$ is defined as:
\begin{align}\label{eq: small g nD dir}
    g_{\phi}(\phi^*,\mathbf{e}):= |\mathbf{e}\cdot \nabla \phi^*|- {\rm sgn}\left(\mathbf{e}\cdot \nabla \phi^* \right)\mathbf{e}\cdot \nabla \phi= - {\rm sgn}\left(\mathbf{e}\cdot \nabla \phi^*\right)\mathbf{e}\cdot \nabla \phi',
\end{align}
and $\mathbf{e} \in \mathbb{R}^d$ is a unit vector (i.e. $\|\mathbf{e}\|_2=1$).
\end{definition}
In the following remark we provide a motivation for this generalization of the one-dimensional Gibbs functional.
\begin{remark}[Pointwise motivation Gibbs constraint] Define the projection operator $\mathcal{P}_{\mathbf{e}}:\mathbb{R}^d \rightarrow \mathbb{R}^d$ in the direction $\mathbf{e}$ by $\mathcal{P}_{\mathbf{e}} \mathbf{v} := \mathbf{P}_{\mathbf{e}}\mathbf{v}$ 
with projection matrix $\mathbf{P}_{\mathbf{e}}=\mathbf{e} \otimes \mathbf{e}$. The projections of $\nabla \phi^*$ and $\nabla \phi$ in direction $\mathbf{e}$ are thus given by $\mathcal{P}_{\mathbf{e}}\nabla \phi^*$ and $\mathcal{P}_{\mathbf{e}}\nabla \phi$, respectively. The form of the local functional $g_{\phi}(\phi^*,\mathbf{e})$ is now motivated by the following equivalence:
   \begin{align}
      &\left. \begin{matrix}   \|\mathcal{P}_{\mathbf{e}}\nabla \phi^*\|_2- \|\mathcal{P}_{\mathbf{e}}\nabla \phi\|_2 \leq 0 &~\Leftrightarrow& |\mathbf{e}\cdot \nabla \phi^*|\leq |\mathbf{e}\cdot \nabla \phi|,\\
       {\rm sgn}\left(\mathcal{P}_{\mathbf{e}}\nabla \phi^*\cdot \mathbf{e}\right) = {\rm sgn}\left(\mathcal{P}_{\mathbf{e}}\nabla \phi \cdot \mathbf{e}\right) &~\Leftrightarrow & {\rm sgn}\left( \mathbf{e}\cdot \nabla \phi^*\right) ={\rm sgn}\left( \mathbf{e}\cdot \nabla \phi \right) \end{matrix} \right\}\nn\\
       &\Leftrightarrow g_{\phi}(\phi^*,\mathbf{e})\leq 0.
   \end{align}
\end{remark}
We have the following simple, but important, property.
\begin{proposition}[Invariance of directional Gibbs functional opposite direction]
  The directional Gibbs functional is invariant with respect to flipping the sign of the direction $\mathbf{e}$:
\begin{align}
    \mathscr{G}_{\phi,\omega}(\phi^*; -\mathbf{e}) = \mathscr{G}_{\phi,\omega}(\phi^*; \mathbf{e}).
\end{align}
\end{proposition}

\begin{figure}[ht!]
\begin{center}
\begin{tikzpicture}[>=stealth,scale=0.80]
\draw[->, line width=0.5mm, blue1 ] (0,0) -- (\xx1,\yy1);
\draw[-, line width=0.5mm, black ] ({\t*cos(atan(\yy3/\xx3)-90)},{\t*sin(atan(\yy3/\xx3)-90)}) -- ({\t*cos(atan(\yy3/\xx3)-90)-\t*sin(atan(\yy3/\xx3)-90)},{\t*sin(atan(\yy3/\xx3)-90)+\t*cos(atan(\yy3/\xx3)-90)});
\draw[-, line width=0.5mm, black ] ({-\t*sin(atan(\yy3/\xx3)-90)},{\t*cos(atan(\yy3/\xx3)-90)}) -- ({\t*cos(atan(\yy3/\xx3)-90)-\t*sin(atan(\yy3/\xx3)-90)},{\t*sin(atan(\yy3/\xx3)-90)+\t*cos(atan(\yy3/\xx3)-90)});
\draw[-, line width=0.5mm, black ] ({\r*cos(atan(\yy2/\xx2)-90)},{\r*sin(atan(\yy2/\xx2)-90)}) -- ({\r*cos(atan(\yy2/\xx2)-90)-\r*sin(atan(\yy2/\xx2)-90)},{\r*sin(atan(\yy2/\xx2)-90)+\r*cos(atan(\yy2/\xx2)-90)});
\draw[-, line width=0.5mm, black ] ({-\r*sin(atan(\yy2/\xx2)-90)},{\r*cos(atan(\yy2/\xx2)-90)}) -- ({\r*cos(atan(\yy2/\xx2)-90)-\r*sin(atan(\yy2/\xx2)-90)},{\r*sin(atan(\yy2/\xx2)-90)+\r*cos(atan(\yy2/\xx2)-90)});
\draw[densely dashdotted, ->, line width=0.5mm, red1 ] (0,0) -- (\xx2,\yy2);
\draw[dotted, ->, line width=0.5mm, green1 ] (0,0)-- (\xx3,\yy3);
\node[text width=0.7cm] at (0.6,2.8) {$\nabla \phi'$};
\node[text width=0.7cm] at (3.0,1.8) {$\nabla \phi^*$};
\node[text width=0.5cm] at (3.5,4.0) {$\nabla \phi$};
\draw[dotted, -, line width=0.5mm, black ] (0,0)-- ({-3*\yy2/sqrt(\xx2*\xx2+\yy2*\yy2)},{3*\xx2/sqrt(\xx2*\xx2+\yy2*\yy2)});
\draw[dotted, -, line width=0.5mm, black ] (0,0)-- ({3*\yy2/sqrt(\xx2*\xx2+\yy2*\yy2)},{-3*\xx2/sqrt(\xx2*\xx2+\yy2*\yy2)});
\draw[dotted, -, line width=0.5mm, black ] (0,0)-- ({-3*\yy3/sqrt(\xx3*\xx3+\yy3*\yy3)},{3*\xx3/sqrt(\xx3*\xx3+\yy3*\yy3)});
\draw[dotted, -, line width=0.5mm, black ] (0,0)-- ({3*\yy3/sqrt(\xx3*\xx3+\yy3*\yy3)},{-3*\xx3/sqrt(\xx3*\xx3+\yy3*\yy3)});
\draw [red,line width=0.5mm,domain={atan(\yy3/\xx3)+90}:{atan(\yy2/\xx2)+90},<->] plot ({cos(\x)}, {sin(\x)});
\draw [red,line width=0.5mm,domain={atan(\yy2/\xx2)-90}:{atan(\yy3/\xx3)-90},<->] plot ({cos(\x)}, {sin(\x)});
\draw [cadmiumgreen,line width=0.5mm,domain={atan(\yy3/\xx3)-90}:{atan(\yy2/\xx2)+90},<->] plot ({cos(\x)}, {sin(\x)});
\draw [cadmiumgreen,line width=0.5mm,domain={atan(\yy3/\xx3)+90}:{atan(\yy2/\xx2)+270},<->] plot ({cos(\x)}, {sin(\x)});
\node[cadmiumgreen,text width=1.5cm] at (-2.0,-1.0) {$g_{\phi} < 0$};
\node[red,text width=1.5cm] at (-2.0,1.5) {$g_{\phi} > 0$};
\node[red,text width=1.5cm] at (2.5,-1.5) {$g_{\phi} > 0$};
\node[cadmiumgreen,text width=1.5cm] at (3.0,0.7) {$g_{\phi} < 0$};
\node[text width=1.5cm] at (4.0,-0.5) {$g_{\phi} = 0$};
\node[text width=1.5cm] at (-3.5,0.5) {$g_{\phi} = 0$};
\node[text width=1.5cm] at (1.8,-2.9) {$g_{\phi} = 0$};
\node[text width=1.5cm] at (-1.3,2.8) {$g_{\phi} = 0$};
\end{tikzpicture}
    \caption{Two-dimensional visualization of the sign of $g_{\phi}=g_{\phi}(\phi^*,\mathbf{e})$ for a specific case. A unit vector $\mathbf{e}$ in the green region corresponds to $g_{\phi}(\phi^*,\mathbf{e})<0$, in the red region to $g_{\phi}(\phi^*,\mathbf{e})>0$, and at the intersections to $g_{\phi}(\phi^*,\mathbf{e})=0$.}
    \label{fig: sign_g}
    \end{center}
\end{figure}
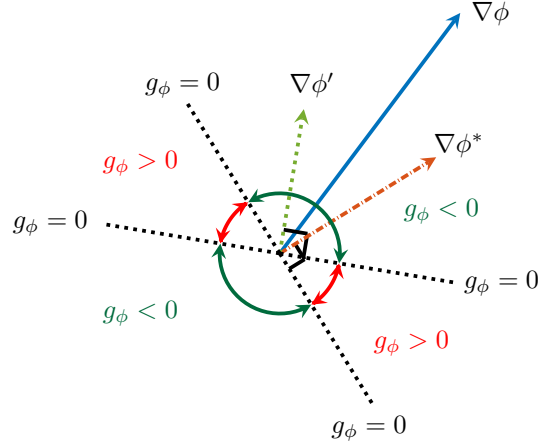

To proceed, we provide a short study on the functional $g_{\phi}=g_{\phi}(\phi^*,\mathbf{e})$.
In the following proposition we comment on the sign thereof, 
depending on the unit vector $\mathbf{e}$.
\begin{proposition}[Sign of directional $g_{\phi}$]\label{prop: sign gK}
  The sign of the functional $g_{\phi}(\phi^*,\mathbf{e})$ depends in the following way on the unit vector $\mathbf{e}$:
  \begin{subequations}
    \begin{align}
        g_{\phi}(\phi^*,\mathbf{e}) < 0 \quad\Leftrightarrow\quad & \left\{\begin{matrix}
        \mathbf{e}\cdot \nabla \phi' > 0 \,\,\text{ and }\,\, \mathbf{e}\cdot \nabla \phi^* >0 & \text{ if }\,\,\mathbf{e}\cdot \nabla \phi > 0\\
        \mathbf{e}\cdot \nabla \phi' < 0  \,\,\text{ and }\,\, \mathbf{e}\cdot \nabla \phi^* <0 & \,\text{ if }\,\,\mathbf{e}\cdot \nabla \phi < 0,\\
        \end{matrix}\right.\\
        g_{\phi}(\phi^*,\mathbf{e}) = 0 \quad\Leftrightarrow\quad & ~\mathbf{e}\cdot\nabla \phi^* = 0 \,\,\text{ or }\,\,\mathbf{e}\cdot \nabla \phi' = 0,\\
        g_{\phi}(\phi^*,\mathbf{e}) > 0 \quad\Leftrightarrow\quad & \left\{\begin{matrix}
        \mathbf{e}\cdot \nabla \phi' < 0 \,\,\text{ or }\,\, \mathbf{e}\cdot\nabla \phi^* < 0& \,\,\,\text{ if }\,\,\mathbf{e}\cdot \nabla \phi > 0\\
        \mathbf{e}\cdot \nabla \phi' > 0 \,\,\text{ or }\,\, \mathbf{e}\cdot\nabla \phi^* > 0& \,\,\,\,\text{ if }\,\,\mathbf{e}\cdot \nabla \phi < 0.
        \end{matrix}\right.
    \end{align}
  \end{subequations}
\end{proposition}
We provide a visualization of the sign of $g_{\phi}(\phi^*,\mathbf{e})$ for varying direction $\mathbf{e}$ for a two-dimensional scenario in \cref{fig: sign_g}. 

In general, there always exists a direction $\mathbf{e}$ for which $g_{\phi}(\phi^*,\mathbf{e}) > 0$, as expressed in 
the following proposition.
\begin{proposition}[Supremum of $g_{\phi}(\phi^*,\mathbf{e})$]\label{prop: sup gK}
  The supremum of $g_{\phi}(\phi^*,\mathbf{e})$ is given by
  \begin{align}
      \sup_{\|\mathbf{e}\|_2=1} g_{\phi}(\phi^*,\mathbf{e}) = \left\{\begin{matrix}
        \mathbf{b}\cdot \nabla \phi' & \text{ if }\,\,\nabla \phi^*\cdot\nabla \phi' \geq 0\\[0.15cm]
        \|\nabla\phi'\|_2 & \,\text{ if }\,\,\nabla \phi^*\cdot\nabla \phi' < 0,
      \end{matrix}\right.
  \end{align}
  where $\mathbf{b} \in \left\{ \mathbf{a} \in \mathbb{R}^d: \|\mathbf{a}\|_2 = 1,\, \mathbf{a}\cdot\nabla \phi^* = 0\right\}$.
\end{proposition}
\begin{proof}

Distinguish the two cases: (i) $\nabla \phi^*\cdot\nabla \phi'<0$ and (ii) $\nabla \phi^*\cdot\nabla \phi'\geq0$.
\begin{enumerate}
    \item[(i)] If $\mathbf{e}\cdot \nabla \phi^*\geq0$, then the supremum of $g_{\phi}$ is $\|\nabla \phi'\|_2$, which is attained at $\mathbf{e} = -\nabla \phi'/\|\nabla \phi'\|_2$. Similarly, if $\mathbf{e}\cdot \nabla \phi^*<0$, then the supremum of $g_{\phi}$ is again $\|\nabla \phi'\|_2$, which is then attained at $\mathbf{e} = \nabla \phi'/\|\nabla \phi'\|_2$.
    \item[(ii)] If $\mathbf{e}\cdot \nabla \phi^*\geq0$, then $g_{\phi}$ is a decreasing function in $\mathbf{e}\cdot \nabla\phi^*$. Therefore, $g_{\phi}$ attains its supremum when $\mathbf{e}\cdot \nabla \phi^* = 0$. Similarly, if $\mathbf{e}\cdot \nabla \phi^*<0$, then $g_{\phi}$ is an increasing function in $\mathbf{e}\cdot \nabla\phi^*$. Again, $g_{\phi}$ attains its supremum when $\mathbf{e}\cdot \nabla \phi^* = 0$.
\end{enumerate} 
\end{proof}

\cref{prop: sign gK,prop: sup gK} show that demanding the elimination of the Gibbs functional in each possible direction, which might intuitively yield the most desirable result, is 
not prudent: selecting the direction $\mathbf{e}$ to maximize $g_{\phi}=g_{\phi}(\phi^*,\mathbf{e})$ and subsequently using this in the Gibbs constraint would only yield (i) the zero approximation $\phi^*=0$ and (ii) the perfect approximation $\phi^*=\phi$ as feasible solutions. To proceed, we must thus select a direction $\mathbf{e}$ that leads to both suitable and practical constraints. 
The relevant directions to consider are functions of the (approximate) solution gradients: (i) $\mathbf{e}=\nabla \phi/\|\nabla \phi\|_2$, (ii) $\mathbf{e}=\nabla \phi^*/\|\nabla \phi^*\|_2$ and (iii) $\mathbf{e}=\nabla \phi'/\|\nabla \phi'\|_2$ leading to the expressions:
\begin{subequations}\label{eq: dir Gibbs indicator}
  \begin{align}
      g_{\phi}\left(\phi^*, \dfrac{\nabla \phi}{\|\nabla \phi\|_2}\right) =&~-|\nabla \phi\cdot\nabla \phi^*|^{-1}\|\nabla \phi\|_2^{-1} (\nabla \phi\cdot\nabla \phi)(\nabla \phi\cdot\nabla \phi'),\\
      g_{\phi}\left(\phi^*, \dfrac{\nabla \phi^*}{\|\nabla \phi^*\|_2}\right) =&~- \|\nabla \phi^*\|_2^{-1}\nabla \phi^*\cdot\nabla \phi',\\
      g_{\phi}\left(\phi^*, \dfrac{\nabla \phi'}{\|\nabla \phi'\|_2}\right) =&~- \|\nabla \phi'\||\nabla \phi^*\cdot\nabla \phi'|^{-1}\nabla \phi^*\cdot\nabla \phi'.
  \end{align}
\end{subequations}

Insisting compatibility with the one-dimensional case requires the choice of direction $\mathbf{e}=\nabla \phi^*/\|\nabla \phi^*\|_2$, which we focus on exclusively in the following.

\begin{definition}[Gibbs functional and constraint]\label{def: Gibbs functional and constraint nD}
The Gibbs functional of the function $\phi^* \in H^1(\tilde{\Omega})$, with respect to the given function $\phi \in H^1(\Omega)$ and on $\omega \subset \Omega$, is defined as:
\begin{align}\label{eq: large g nD}
    \mathscr{G}_{\phi,\omega}(\phi^*) := \displaystyle\int_{\omega} g_{\phi}(\phi^*)~{\rm d}\Omega,
\end{align}
where the functional $g_\phi$ is defined as:
\begin{align}\label{eq: small g nD}
    g_{\phi}(\phi^*) :=&~ \left\{ \begin{matrix}
      - \|\nabla \phi^*\|_2^{-1}\nabla \phi^*\cdot\nabla \phi' & \text{if }\,\, \nabla\phi^* \neq 0,\\[0.15cm]
      0 & \text{if }\,\, \nabla\phi^* = 0.\end{matrix}\right.
\end{align}
The associated Gibbs constraint reads:
\begin{align}\label{eq: Gibbs constraint nD}
    \mathscr{G}_{\phi,\omega}(\phi^*) \leq 0.
\end{align}
\end{definition}

Note that the Gibbs functional is frame-invariant, and that the properties of \cref{prop: perfect approx} and \cref{prop: lower bound G} are inherited from the one-dimensional case.

\subsection{Gibbs constraints for finite element best approximations}\label{subsec: Gibbs higher D: Gibbs fem}

In this subsection, we apply the Gibbs constraints to finite element best approximations. Again, we consider the following general form of the best approximation problem:\\

\textit{find $\phi^h \in \mathcal{V}^h_{D;p,\alpha}$ such that:}
\begin{subequations}
\begin{align}\label{eq: optimization problem fem 2D}
    \phi^h =&~ \underset{\theta^h \in \mathcal{K}_{p,\alpha}}{\rm arginf} \|\phi-\theta^h\|_{\HH},
\end{align}
\indent \textit{where the feasible set $\mathcal{K}_{p,\alpha}$ 
is defined as:}
\begin{align}\label{eq: constraint set fem 2D}
  \mathcal{K}_{p,\alpha} := \left\{ \phi^h \in \mathcal{V}^h_{D;p,\alpha} :~\right.&\left. \mathscr{G}_{\phi,\omega_j}(\phi^h) \leq 0, j = 1,...,J\right\}.
\end{align}
\end{subequations}

\begin{figure}[!ht]
\begin{subfigure}{0.49\textwidth}
\includegraphics[width=0.95\textwidth]{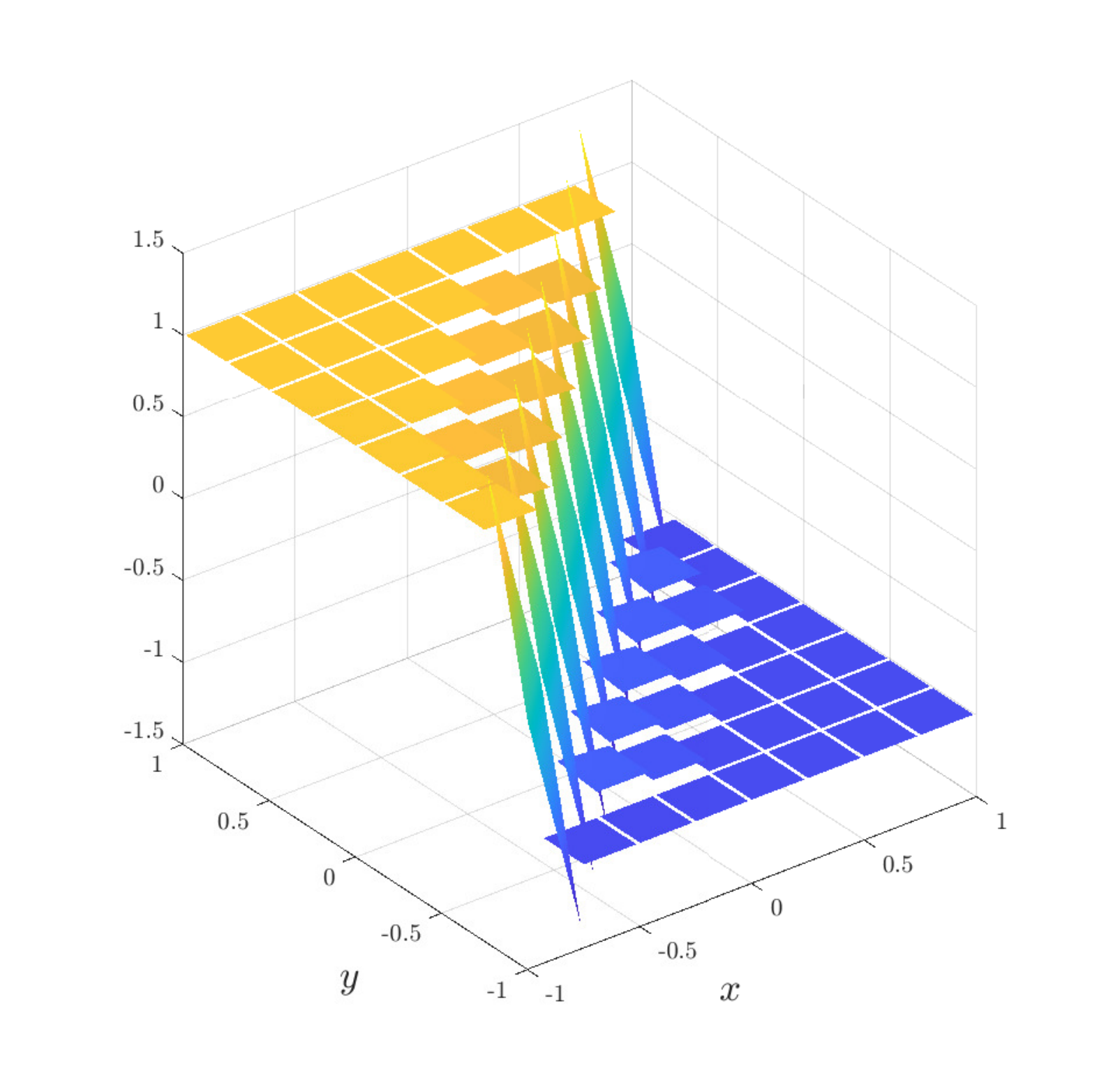}
\caption{$p=1$}
\end{subfigure}
\begin{subfigure}{0.49\textwidth}
\includegraphics[width=0.95\textwidth]{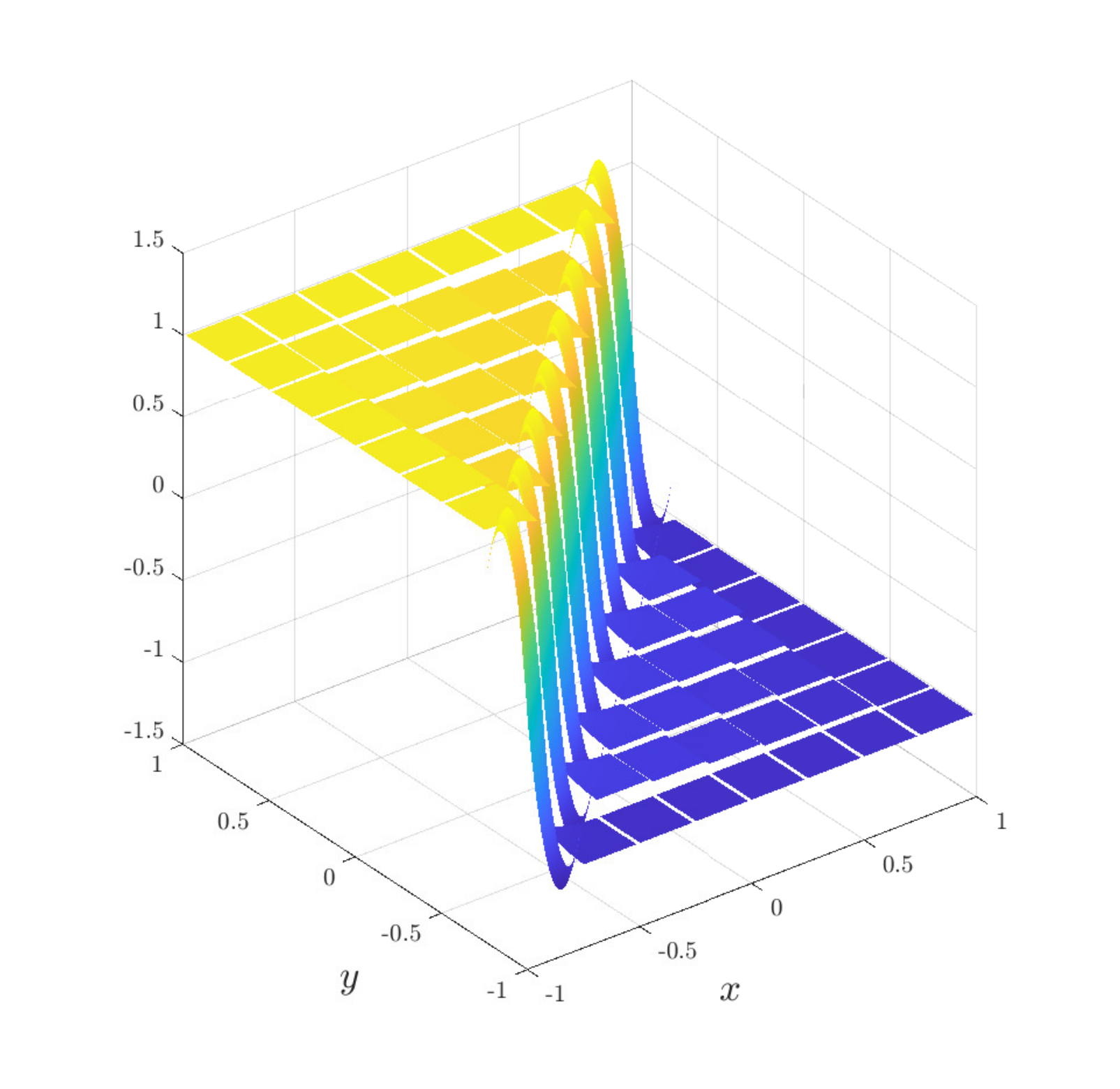}
\caption{$p=2$}
\end{subfigure}
\caption{The interior penalty-best approximation $\phi^h \in \mathcal{V}_{D;p,-1}^h$ of the smooth step function, with $n_{\rm el} =8\times 8$, subject to element-wise Gibbs constraints, for different $p$.}
\label{fig: IP Gibbs 2D}
\end{figure}

We first consider discontinuous finite element approximation spaces ($\alpha = -1$). Analogous to the one-dimensional case, we consider \eqref{eq: optimization problem fem 2D}-\eqref{eq: constraint set fem 2D} with $\mathcal{H}={\rm IP}$, i.e. the interior penalty best approximation. We select the subdomains $\omega_j$ as the finite elements $K_j, j = 1, \dots, n_{{\rm el}}$. Feasibility of the constrained optimization problem is a consequence of the existence of an elementwise constant solution $\phi^h$, for which $\mathscr{G}_{\phi,K_j}(\phi^h)=0$. Again, we may analyze the solution properties via the KKT conditions \eqref{eq: KKT fem}. In case of homogeneous boundary conditions, we have $(\phi^h,\phi')_{\mathcal{H}}\geq 0$ (\cref{lem: positive inner product}). In \cref{fig: IP Gibbs 2D}, we visualize the interior penalty-best approximation $\phi^h \in \mathcal{V}_{D;p,-1}^h$, $p=1, 2$ of the profile of example 1, subject to the element-wise Gibbs constraints. We observe that the constrained interior penalty-best approximations display a significant reduction of over- and undershoots compared to the non-constrained solutions depicted in \cref{sec: review Gibbs 2D}.

\begin{figure}[ht]
\begin{subfigure}{0.55\textwidth}
\includegraphics[scale = 0.55]{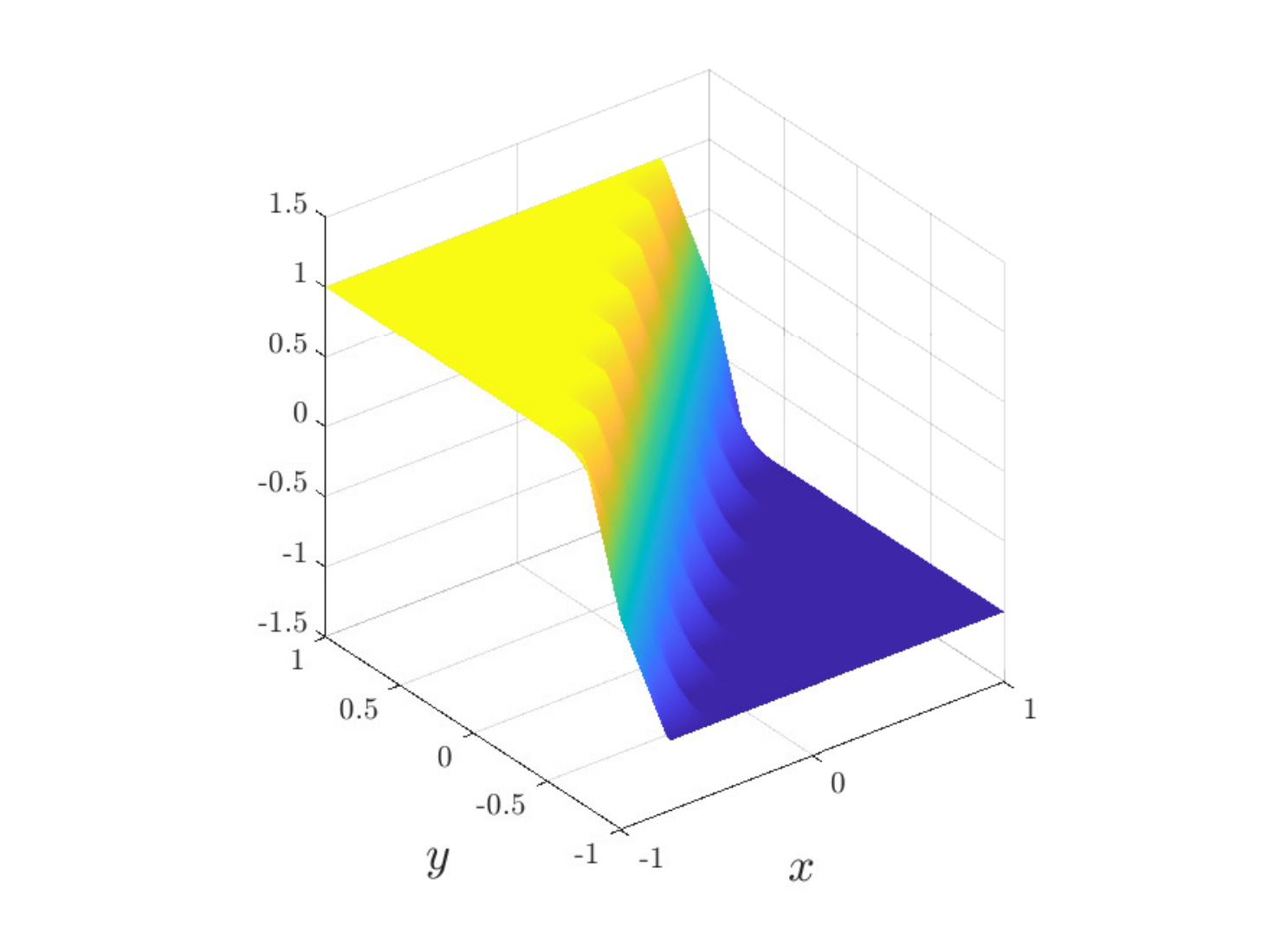}\\
    \caption{The approximation $\phi^h \in \mathcal{V}^h_{1,0}$.}
\end{subfigure}
\begin{subfigure}{0.44\textwidth}
\begin{center}
\scalebox{0.60}{\begin{tikzpicture}[>=stealth,
axis/.style={densely dashed,font=\small}]
\foreach \i in {0,1,2,3,4,5,6,7,8} {
    \draw [very thin,gray] (\i,\yMin) -- (\i,\yMax);
}
\foreach \i in {\yMin,...,\yMax} {
    \draw [very thin,gray] (\xMin,\i) -- (\xMax,\i);
}
\node[text width=2.5cm, green1] at (1.3,1.5) {\normalsize$-0.33$};
\node[text width=2.5cm, green1] at (2.3,2.5) {\normalsize$-0.33$};
\node[text width=2.5cm, green1] at (3.3,3.5) {\normalsize$-0.33$};
\node[text width=2.5cm, green1] at (4.3,4.5) {\normalsize$-0.33$};
\node[text width=2.5cm, green1] at (5.3,5.5) {\normalsize$-0.33$};
\node[text width=2.5cm, green1] at (6.3,6.5) {\normalsize$-0.33$};
\node[text width=2.5cm, green1] at (7.3,7.5) {\normalsize$-0.33$};
\node[text width=2.5cm, green1] at (8.3,8.5) {\normalsize$-0.33$};

\node[text width=2.5cm, carnelian] at (2.45,1.5) {\normalsize$0.18$};
\node[text width=2.5cm, carnelian] at (1.45,2.5) {\normalsize$0.18$};
\node[text width=2.5cm, carnelian] at (3.45,2.5) {\normalsize$0.18$};
\node[text width=2.5cm, carnelian] at (2.45,3.5) {\normalsize$0.18$};
\node[text width=2.5cm, carnelian] at (4.45,3.5) {\normalsize$0.18$};
\node[text width=2.5cm, carnelian] at (3.45,4.5) {\normalsize$0.18$};
\node[text width=2.5cm, carnelian] at (5.45,4.5) {\normalsize$0.18$};
\node[text width=2.5cm, carnelian] at (4.45,5.5) {\normalsize$0.18$};
\node[text width=2.5cm, carnelian] at (6.45,5.5) {\normalsize$0.18$};
\node[text width=2.5cm, carnelian] at (5.45,6.5) {\normalsize$0.18$};
\node[text width=2.5cm, carnelian] at (7.45,6.5) {\normalsize$0.18$};
\node[text width=2.5cm, carnelian] at (6.45,7.5) {\normalsize$0.18$};
\node[text width=2.5cm, carnelian] at (8.45,7.5) {\normalsize$0.18$};
\node[text width=2.5cm, carnelian] at (7.45,8.5) {\normalsize$0.18$};
\node[text width=2.5cm, white] at (1.3-1,1.5-2) {\normalsize$-0.33$};
\end{tikzpicture}}\\
\caption{Element-wise Gibbs functionals.}    
\end{center}
\end{subfigure}
\caption{The sharpest possible approximation $\phi^h \in \mathcal{V}^h_{1,0}$ (a) and the corresponding element-wise values of the Gibbs functionals $\mathscr{G}_{\phi,K_i}(\phi^h)$ (b).}
    \label{fig: 2D linears sharp}
\end{figure}

\newpage
Next, we study the case with regularity $\alpha = 0$. Consider again the example with the sharp layer skew to the mesh. In \cref{fig: 2D linears sharp} we visualize the sharpest possible approximation $\phi^h \in \mathcal{V}^h_{D,1,0}$, and plot the element-wise values of the Gibbs functionals $\mathscr{G}_{\phi,K_i}(\phi^h)$. We note that the sharpest approximation is interpolatory. Nevertheless, the Gibbs constraints are not fulfilled element-wise. Similar to the one-dimensional case, one could define $\omega_j$ as the collection of several elements. However, a careful examination reveals that the only possible choice is the collection of all elements $\omega = \cup_{j=1}^{n_{{\rm el}}} K_j$, which would constitute a rather weak condition. It appears impossible to obtain a feasible solution when using subdomains $\omega_j$ that consist of a few elements. This is a consequence of the support of the basis functions, which contains multiple elements. The Gibbs constraints for any collection of elements $K_i$ on the off-diagonal with $x\geq y ~(x\leq y)$ would enforce $\phi^h_{K_i}$ to be close to $1$ (and $-1$). This is in conflict with the continuity requirement of the finite element approximation space $V^h_{D,1,0}$, and persists when using higher order polynomials.

\begin{figure}
\begin{center}
\scalebox{0.65}{\input{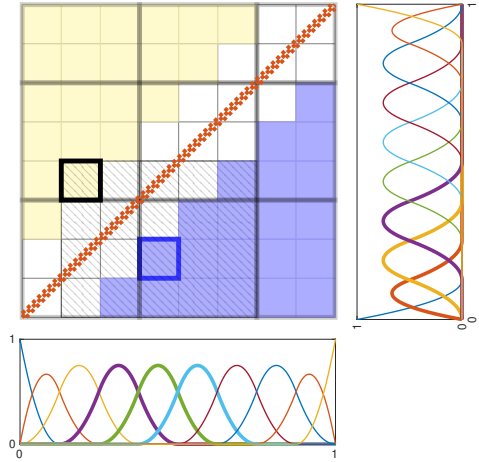}}
    \caption{Subdivision of an $8 \times 8$ element domain into groups of maximum $3 \times 3$ elements with approximation space $\mathcal{V}^h_{D,2,1}$. The right and bottom inset figures display the univariate basis functions along the corresponding axes.}
    \label{fig: subdiv 2D Bsplines}\end{center}
\end{figure}

\newpage
Lastly, we consider higher-order smooth finite element approximation spaces ($\alpha \geq 1$). The above discussed infeasibility of the optimization problem also applies to higher-order smooth approximations. We illustrate this for the smooth step function with a quadratic B-spline approximation space $\mathcal{V}^h_{D,2,1}$. Consider a subdivision into groups of maximum $3 \times 3$ elements, as visualized in \cref{fig: subdiv 2D Bsplines}.
In order to satisfy the Gibbs constraint on the element group on the middle-bottom we must have $\phi^h = -1$ in the blue boxed element. Similarly, in the black boxed element we require $\phi^h = 1$. These two conditions are incompatible.

\section{Conclusions}
\label{sec: conclusions}
In this article, we constructed a set of integral constraints with the aim of eliminating the Gibbs phenomenon in finite element best approximations. We first provided an overview of the Gibbs phenomenon for best approximations in finite element spaces. We illustrated with computational examples that spurious oscillations occur in one and two dimensions for standard projections onto finite element spaces of arbitrary degree and regularity (with the exception of the one-dimensional $H_0^1$-projection with linear continuous finite elements). The proposed constraints 
build onto the concept of total variation. In this regard, we established in one dimension the interrelation between the Gibbs constraint and interpolatory and monotonic approximations, as well as the maximum principle. Furthermore, we displayed in one dimension that the proposed constraints may be applied element-wise when the finite element space is either discontinuous or $\mathcal{C}^0$-continuous. For higher regularity finite element spaces, the integration domains of the constraints depend on multiple elements. We showed that enforcing the constraints removes over- and undershoots for continuous finite element spaces, and suppresses them for discontinuous finite element spaces. In higher dimensions, the constraints act in the direction of the solution gradient. The applicability of the constraints is then limited to discontinuous finite element spaces. We demonstrated that also in this case over- and undershoots are severely reduced. 

We recognize two open problems that require further investigation. The first is linked to the last observation, namely the extension of the set of constraints 
to continuous finite element spaces in higher dimensions. The second open problem is the construction of a finite element method that incorporates these constraints in practical computations. We conjecture that a possible resolution lies in the variational multiscale framework, in particular through the work of Evans et al. \cite{evans2009enforcement} and of the first author \cite{ten2019variation}. 

\section*{Acknowledgments}
The authors are grateful to Thomas J.R. Hughes for insightful discussions on the topic. MtE was supported by the German Research Foundation (Deutsche Forschungsgemeinschaft DFG) via the Walter Benjamin project EI 1210/1-1. DS gratefully acknowledges support from the German Research Foundation (Deutsche Forschungsgemeinschaft DFG) via the Emmy Noether Award SCH 1249/2-1.


\bibliographystyle{unsrtnat}
\bibliography{references}

\end{document}